%% file: conjecture.tex
\documentclass[noaddress, tikz, sepfignums]{nmd/article}
\ifnmd
\fi
\usepackage{microtype}
\usepackage{bm}
\usepackage{enumitem}
\usepackage{tikz-cd}
\usepackage{colortbl}
\usetikzlibrary{arrows}
\usetikzlibrary{shapes}
\usetikzlibrary{intersections}
\tikzset{
  commutative diagrams/.cd,
  arrow style=tikz,
  diagrams={>=stealth', line width=0.7pt},
}

\input tikz_preamble

\pgfkeys{/matplotlibfigure,
  default/.style = {
    width = 0.5\textwidth,
    font = \fontsize{9}{10.7}\selectfont}
}

\tikzset{plot legend/.style = {
    below left,
    align=right,
    font=\scriptsize,
    fill=nmddark!10,
    line width=0.4pt,
    rounded corners=3
  }
}

\title{Floer homology, \\ group orderability, \\
  and taut foliations of \\ hyperbolic 3-manifolds}

\author{Nathan M. Dunfield}
\givenname{Nathan}
\surname{Dunfield}
\address{ Dept.~of Math., MC-382 \\
          University of Illinois \\
          1409 W. Green St. \\
          Urbana, IL 61801 \\ 
          USA
}
\email{nathan@dunfield.info}
\urladdr{http://dunfield.info}

\arxivreference{}
\arxivpassword{}

\renewcommand{\Im}{\operatorname{Im}}
\newcommand{\QHSpheres}{$\mathcal{Y}$}
\newcommand{\mathQHSpheres}{\mathcal{Y}}
\newcommand{\ab}{\mathrm{ab}}

\newcommand{\cFtil}{\widetilde{\cF}}
\newcommand{\cTtil}{\widetilde{\cT}}
\newcommand{\IR}{\mathrm{I}\hspace{0.05em}\mathbb{R}}
\newcommand{\IC}{\mathrm{I}\mathbb{C}}
\newcommand{\IG}{\mathrm{I}G}
\newcommand{\IGtil}{\mathrm{I}\Gtil}

\renewcommand{\a}{\bm{a}}
\renewcommand{\b}{\bm{b}}
\renewcommand{\c}{\bm{c}}
\renewcommand{\d}{\bm{d}}
\renewcommand{\t}{\bm{t}}
\newcommand{\g}{\bm{g}}
\newcommand{\z}{\bm{z}}
\newcommand{\x}{\bm{x}}
\newcommand{\brho}{\bm{\rho}}
\newcommand{\brhotil}{\widetilde{\bm{\rho}}}
\newcommand{\slopes}{\mathit{Sl}}
\newcommand{\Bpunc}{{B^\circ}}
\newcommand{\TF}{{T\cF}}
\newcommand{\TB}{{T\!B}}
\newcommand{\UTB}{{U\!T\!B}}

\newcommand{\UTBpunc}{{U\!T\!\Bpunc}}
\newcommand{\iotainvDtau}[1]{\iota^{-1}\left(\cD^\tau_{>0}(#1)\right)}
\DeclareMathOperator{\interior}{int}
\DeclareMathOperator{\mixed}{mixed}
\DeclareMathOperator{\degen}{degen}
\DeclareMathOperator{\Core}{Core}
\DeclareMathOperator{\FreeGroup}{FreeGroup}

\usepackage[subfigure]{tocloft}
\advance\cftsecnumwidth 0.2em\relax
\advance\cftsubsecindent 0.2em\relax
\advance\cftsubsecnumwidth 0.2em\relax

\begin{document}

\begin{abstract} 
  This paper explores the conjecture that the following are equivalent
  for irreducible rational homology 3-spheres: having left-orderable
  fundamental group, having non-minimal Heegaard Floer homology, and
  admitting a co-orientable taut foliation.  In particular, it adds
  further evidence in favor of this conjecture by studying these three
  properties for more than 300,000 hyperbolic rational homology
  3-spheres.  New or much improved methods for studying each of these
  properties form the bulk of the paper, including a new combinatorial
  criterion, called a foliar orientation, for showing that a
  3-manifold has a taut foliation.
\end{abstract}
\maketitle

\tableofcontents



\input intro

\input floer

\input disorder
\input foliations
\input representations

{\RaggedRight \bibliographystyle{nmd/math} \small
  \bibliography{\jobname} }
\end{document}

%% file: tikz_preamble.tex
\usetikzlibrary{hobby}
\usetikzlibrary{patterns}

\tikzset{
      hatch distance/.store in=\hatchdistance,
      hatch distance=10pt,
      hatch thickness/.store in=\hatchthickness,
      hatch thickness=2pt,
      hatch color/.store in=\hatchcolor,
      hatch color=black,
  }

\pgfdeclarepatternformonly[%
      \hatchdistance,%
      \hatchthickness,%
      \hatchcolor%
          ]{flexible hatch north east}
  {\pgfqpoint{0pt}{0pt}}
  {\pgfqpoint{\hatchdistance}{\hatchdistance}}
  {\pgfpoint{\hatchdistance-1pt}{\hatchdistance-1pt}}%
  {
      \pgfsetlinewidth{\hatchthickness}
      \pgfpathmoveto{\pgfqpoint{0pt}{0pt}}
      \pgfpathlineto{\pgfqpoint{\hatchdistance}{\hatchdistance}}
      \pgfsetstrokecolor{\hatchcolor}%
      \pgfusepath{stroke}
  }

\pgfdeclarepatternformonly[%
      \hatchdistance,%
      \hatchthickness,%
      \hatchcolor%
          ]{flexible hatch north west}
  {\pgfqpoint{0pt}{0pt}}
  {\pgfqpoint{\hatchdistance}{-\hatchdistance}}
  {\pgfpoint{\hatchdistance-1pt}{-\hatchdistance+1pt}}%
  {
      \pgfsetlinewidth{\hatchthickness}
      \pgfpathmoveto{\pgfqpoint{0pt}{0pt}}
      \pgfpathlineto{\pgfqpoint{\hatchdistance}{-\hatchdistance}}
      \pgfsetstrokecolor{\hatchcolor}%
      \pgfusepath{stroke}
}

\makeatletter
\def\IfTikzLibraryLoaded#1{%
  \ifcsname tikz@library@#1@loaded\endcsname
    \expandafter\@firstoftwo
  \else
    \expandafter\@secondoftwo
  \fi
}
\makeatother    

%% file: intro.tex
\section{Introduction}

\subsection{The motivating conjecture}

\begin{figure}
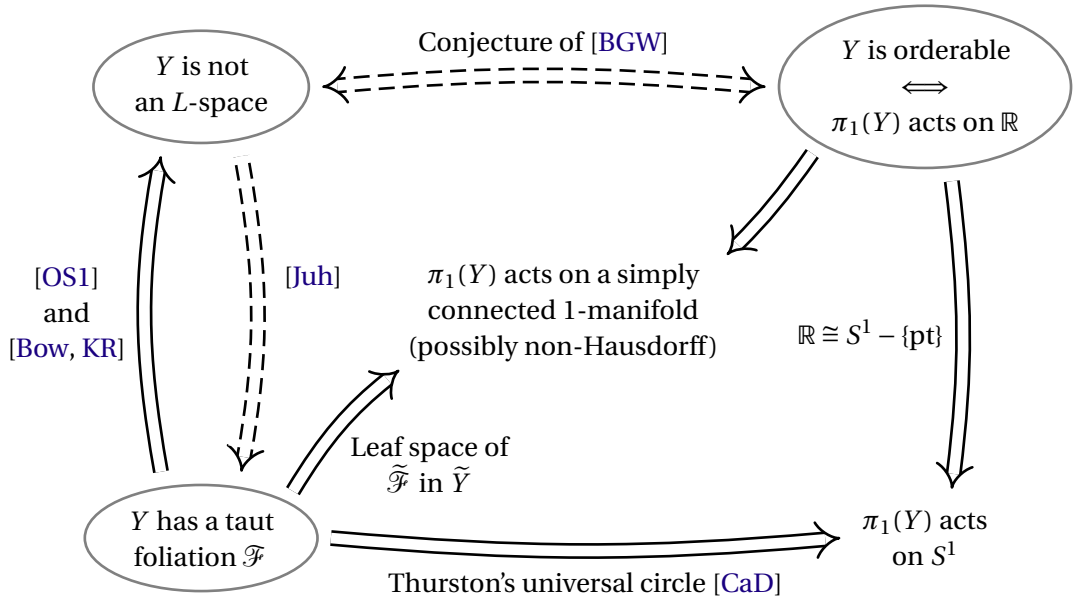

  \begin{center}
    \input plots/diagram
  \end{center}
  \caption{Some results related to Conjecture~\ref{BGWconjecture},
    which asserts the equivalence of the three circled conditions.
    Here $Y$ is an irreducible $\Q$-homology 3-sphere, all
    foliations are co-orientable, and all actions are nontrivial,
    faithful, and orientation preserving; the solid arrows are
    theorems and dotted ones conjectures.  See
    \cite{BoyerGordonWatson2013} for a complete discussion. This
    figure is copied from \cite{CullerDunfield2018}.
  }
  \label{fig:conj}
\end{figure}

Throughout this introduction, please see Section~\ref{sec:terms} for
precise definitions and conventions, which include that all
\3-manifolds are orientable and all foliations are co-orientable.
This paper explores the following:
\begin{conjecture}\label{BGWconjecture}
  For an irreducible $\Q$-homology \3-sphere $Y$, the following
  are equivalent:
  \begin{enumerate}
  \item \label{item:orderable}
    $Y$ is orderable, i.e.~its fundamental group $\pi_1(Y)$ is left-orderable;
  \item \label{item:notL} $Y$ is not an L-space, i.e.~its Heegaard
    Floer homology is not minimal;
  \item \label{item:taut}
    $Y$ admits a taut foliation. 
  \end{enumerate}
\end{conjecture}
The equivalence of (\ref{item:orderable}) and (\ref{item:notL}) was
boldly postulated by Boyer, Gordon, and Watson in
\cite{BoyerGordonWatson2013}, which includes a detailed discussion of
this conjecture.  The equivalence of (\ref{item:notL}) and
(\ref{item:taut}) was formulated as a question by Ozsv\'ath and
Szab\'o after they proved that (\ref{item:taut}) implies
(\ref{item:notL}) \cite{OSgenusbounds2004, KazezRoberts2017,
  Bowden2016}, and upgraded to a conjecture in \cite{Juhasz2015}. On
its face, Conjecture \ref{BGWconjecture} is quite surprising given the
disparate nature of these three conditions, but there are actually a
number of interconnections between them summarized in
Figure~\ref{fig:conj}.  Despite much initial skepticism, substantial
evidence has accumulated in favor of Conjecture~\ref{BGWconjecture}.
For example, it holds for \emph{all} graph manifolds
\cite{HanselmanRasmussenRasmussenWatson2015, Rasmussen2017,
  BoyerClay2014} and many branched covers of knots in the \3-sphere
\cite{GordonLidman2014}, as well as for certain families of Dehn
surgeries on a fixed manifold \cite{CullerDunfield2018}.  Here,
despite my best efforts to disprove this conjecture, I add to this
evidence by exploring these properties for more than 300{,}000
hyperbolic rational homology \3-spheres.

This was challenging in part because the property in
(\ref{item:orderable}) is not known to be algorithmically decidable
(and it is undecidable in the broader category of all finitely
presented groups), and while property (\ref{item:taut}) is known to be
algorithmically decidable, the current algorithm is believed by its
authors to be ``nearly impossible to implement on a computer''
\cite{AgolLi2003}.  The bulk of this paper is devoted to giving new or
much improved methods for exploring all three of these properties; see
Sections~\ref{intro:fol}--\ref{intro:floer} for an overview.  However,
let me first describe Theorem~\ref{thm:main}, which is the main result
here supporting Conjecture~\ref{BGWconjecture}.

\subsection{A few rational homology 3-spheres}\label{sec:3spheres}

Here I consider a census, denoted \QHSpheres, of some 307{,}301
rational homology 3-spheres which are described in
Section~\ref{sec:QHSpheres}. Each manifold in \QHSpheres\ is a Dehn
filling of a 1-cusped hyperbolic 3-manifold that can be triangulated
with at most 9 ideal tetrahedra; the latter were enumerated by Burton
\cite{Burton2014}.

\begin{theorem}\label{thm:arehyp}
  The rational homology 3-spheres in \QHSpheres are all hyperbolic and
  distinct.
\end{theorem}
Additionally, I have strong numerical evidence that the systole, that
is, the length of the shortest closed geodesic, is at least 0.2 for
all manifolds in \QHSpheres.  In fact, I conjecture that \QHSpheres\
is precisely the set of all hyperbolic rational homology 3-spheres
that are Dehn fillings on 1-cusped manifolds from \cite{Burton2014}
where the systole is at least 0.2.

For comparison, the Hodgson-Weeks census consists of $11{,}031$ closed
hyperbolic 3-manifolds with systole at least 0.3 and that are Dehn
fillings on manifolds triangulated by at most 7 ideal tetrahedra
\cite{HodgsonWeeks1994}.  In particular, the census \QHSpheres\
contains the 10{,}903 rational homology \3-spheres in the
Hodgson-Weeks census, and those make up $3.5\%$ of its total.

\subsection{Overall results}

\begin{figure}
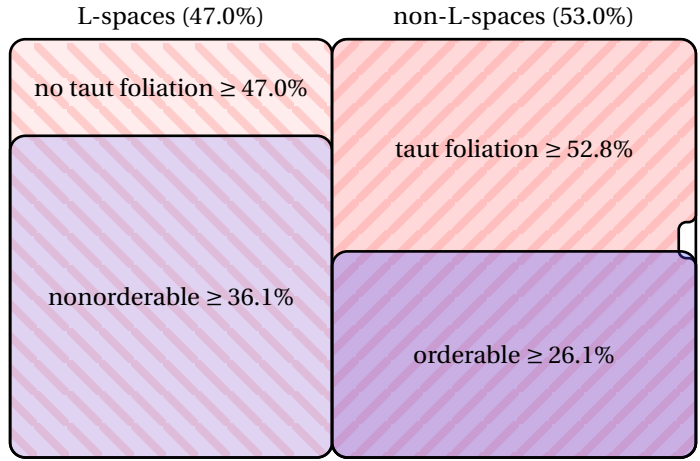

  \begin{center}
      \input plots/summary
    \end{center}
    \caption{A visual summary of what Theorem~\ref{thm:main} says
      about the rational homology 3-spheres in $\cY$.  In particular,
      Conjecture~\ref{BGWconjecture} holds in full for the 62.2\% of
      $\cY$ corresponding to the lowest two regions above.  Moreover,
      the equivalence of parts (\ref{item:notL}) and (\ref{item:taut})
      in the conjecture holds for 99.8\% of $\cY$, namely everything
      except the notch on the middle of the right side.}
    \label{fig:main}
\end{figure}

The main result of this paper supporting
Conjecture~\ref{BGWconjecture} is the following, which is summarized
in Figure~\ref{fig:main}.

\begin{theorem}\label{thm:main}
  Of the 307{,}301 hyperbolic rational homology 3-spheres in \QHSpheres:
  \begin{enumerate}
    \item \label{item:Lsp}
      Exactly 144{,}298 (47.0\%) are L-spaces and 163{,}003
      (53.0\%) are non-L-spaces.

    \item \label{item:fol} At least 162{,}341 (52.8\%) of these
      manifolds admit taut foliations; this is 99.6\% of the
      non-L-spaces.

    \item \label{item:euler0} At least 80{,}236 (26.1\%) of these
      manifolds are orderable; all of the known orderable manifolds
      are non-L-spaces.

    \item \label{item:nonorder} At least 110{,}940 (36.1\%) of these
      manifolds are not orderable; all of the known nonorderable
      manifolds are L-spaces.
  \end{enumerate}
  Overall, Conjecture~\ref{BGWconjecture} holds for at least
  191{,}089 (62.2\%) of these manifolds.
\end{theorem}
I now turn to summarizing the techniques used to prove this theorem,
which form the real heart of the paper.

\subsection{Constructing foliations}
\label{intro:fol}

In Section~\ref{sec:foliar}, I give a new purely combinatorial
technique for constructing a taut foliation on a closed \3-manifold
$Y$: a \emph{foliar orientation} of the edges of a triangulation $\cT$
for $Y$.  This notion is a strengthing of the local orientation of
Calegari \cite{Calegari2000}, and a foliar orientation has an
associated branched surface which is a canonical smoothing of the
2-skeleton of the dual cell complex to $\cT$. Li's theory of laminar
branched surfaces \cite{Li2002} turns out to apply to this branched
surface, showing that it carries a lamination that can be extended to
a taut foliation of $Y$; see Theorems~\ref{thm:foliar} and
\ref{thm:foliarmulti} for more.  Such foliar orientations turn out to
be extraordinarily common, occurring for more than 160,000 of the
manifolds in $\cY$ by Theorem~\ref{thm:usefoliar} and providing the
bulk of the proof of Theorem~\ref{thm:main}(\ref{item:fol}).  It is
unclear whether every taut foliation arises from a foliar orientation
on some triangulation; see Remark~\ref{rem:allfoliar} for some
possible approaches to this question.

The closely related notation of a \emph{persistently foliar
  orientation} on an ideal triangulation of a compact \3-manifold $M$
whose boundary is a torus is introduced in Section~\ref{sec:persist}.
Theorem~\ref{thm:persist} shows that having a persistently foliar
orientation means that all but at most one Dehn filling on $M$ has a
taut foliation. In addition to being used in the proof of
Theorem~\ref{thm:main}(\ref{item:fol}), persistently foliar
orientations are ubiquitous on the exteriors of knots in $S^3$:
\begin{restatable*}{theorem}{persistentlynonalternating}
\label{thm:nonalt}
  Among the 1{,}210{,}608 nonalternating prime knots with at most 16
  crossings, there are exactly 12 that are L-space knots.  All the
  others have ideal triangulations of their exteriors with
  persistently foliar orientations; in particular, any nontrivial Dehn
  surgery on one of these knots admits a co-orientable taut foliation.
\end{restatable*}
\noindent
Motivated by this and the work of Delman and Roberts
\cite{DelmanRobertsTBD} in the case of alternating knots, I posit in
Conjecture~\ref{conj:nonLpersist} that the exterior of a non-L-space
knot in $S^3$ always has a persistently foliar branched surface.

\subsection{Nonorderability and the word problem}

For the orderability of $Y$, I used three separate techniques.  The
first of these, described in Section~\ref{sec:disorder}, is a method
for showing that $\pi_1(Y)$ is \emph{not} left-orderable.  The basic
approach follows \cite[\S 8--10]{CalegariDunfield2003}, but with the
key change being how the word problem is solved in $\pi_1(Y)$.  Rather
than using the theory of automatic groups, I use a new approach
specific to the fundamental group of a finite-volume hyperbolic
3-manifold $Y$.  The idea is to use a numerical approximation of the
holonomy representation $\pi_1(Y) \to \PSL{2}{\C}$ that has been
rigorously verified as correct to within some small tolerance via the
interval analysis method of \cite{HIKMOT2016}.  See
Section~\ref{sec:wordprob} for details.

\subsection{Orderability via foliations and $\PSLRtilde$} 

The remaining techniques for studying orderability were a pair of
independent methods for showing $Y$ is orderable.  The first uses taut
foliations, specifically the fact that if $Y$ has a taut foliation
$\cF$ whose Euler class $e(\cF) \in H^2(Y;\Z)$ vanishes then $Y$ is
orderable \cite[Theorem 8.1]{BoyerHu2018}.  Section~\ref{sec:foleuler}
explains how to calculate $e(\cF)$ when $\cF$ comes from a
foliar orientation, and then uses this to show some 32{,}347 of the
manifolds in $\cY$ are orderable in Theorem~\ref{thm:tauteuler0}.

The second technique for showing orderability is the much-used method
of finding a nontrivial representation
$\rhotil \maps \pi_1(Y) \to \PSLRtilde$; see \cite{CullerDunfield2018}
and the references therein for many prior examples of this.  The new
feature is that I prove the existence of $\rhotil$ using interval
analysis in the same spirit as \cite{HIKMOT2016}.  The fact that
$\PSLRtilde$ is nonlinear complicates matters somewhat as you will see
in Section~\ref{sec:prove_order}, but in the end I successfully
applied it to $64{,}180$ manifolds in Theorem~\ref{thm:tilreps}.

The starting point for Theorem~\ref{thm:tilreps} was a numerical study
of representations of $\pi_1(Y)$ to $\SL{2}{\C}$ for all the
$Y$ in $\cY$ which is described in Section~\ref{sec:numreps}.  There,
using Ptolemy coordinates and numerical algebraic geometry, I found
compelling evidence of 27.8~million such representations, summarized
in Table~\ref{table:reps}.  One interesting observation, given the
importance of $\PSLRtilde$-representations in earlier work on
Conjecture~\ref{BGWconjecture}, is that the L-spaces actually had
more representations to $\SL{2}{\R}$ than the non-L-spaces.  In fact,
if the Euler classes of the $\SL{2}{\R}$-representations found were
simply random elements of $H^2(Y; \Z)$, you would expect $\cY$
to contain about 6{,}000 counterexamples to
Conjecture~\ref{BGWconjecture} (see Remark~\ref{rem:crazy}).  This
allows us to reject the hypothesis that these Euler classes are random
with $p = 10^{-2{,}700}$, providing yet more evidence for
Conjecture~\ref{BGWconjecture}.

\subsection{Computing Floer homology}
\label{intro:floer}

The property of being an L-space is algorithmically decidable by
\cite{SarkarWang2010}.  Moreover, the bordered Heegaard Floer theory
of \cite{LipshitzOzsvathThurston2011, LipshitzOzsvathThurston2014}
provides powerful and effective computational tools for determining
this.  However, rather than attacking the problem head on, I chose to
use a bootstrapping procedure that exploited the structure of the big
Dehn Surgery graph \cite{HoffmanWalsh2015} via the results of
Rasmussen and Rasmussen \cite{Rasmussens2015} on L-space Dehn
fillings.

Suppose $M$ is a compact \3-manifold with $\partial M$ a torus.  When
$M$ has two Dehn fillings that are L-spaces it is \emph{Floer simple},
and \cite{Rasmussens2015} shows how to completely determine which Dehn
fillings are L-spaces by using essentially only the ordinary Alexander
polynomial.  By definition, the manifolds in $\cY$ are Dehn fillings
on a collection $\cC$ of manifolds with torus boundary, and many
manifolds in $\cY$ have multiple such descriptions.  Starting with the
complete list of all exceptional Dehn fillings on $\cC$ provided by
\cite{DunfieldExceptional}, I applied the definition to see that
almost 20\% of the manifolds in $\cC$ are Floer simple.  Then
\cite{Rasmussens2015} determines the L-space status of all Dehn
fillings on those manifolds, in particular showing that almost 8\% of
$\cY$ are L-spaces.  This in turn shows that even more manifolds in
$\cC$ are Floer simple.  Repeating this and several related
deductions, I eventually recovered the complete picture of which
manifolds in $\cY$ are L-spaces.  See Section~\ref{sec:compute-HF} for
details and Table~\ref{table:bootstrap} for a summary.

\subsection{Code and data}

The proof of Theorem~\ref{thm:main} above is of course heavily
computer-assisted.  Moreover, discovering it took many CPU-decades of
computational time, quite possibly several CPU-centuries, using a
computer cluster with a few hundred processor cores.  However,
checking the final proof is much faster.  For example, it is quick to
check that a saved edge orientation of a particular triangulation is
foliar, but much time can be spent searching through triangulations
and orientations in hopes of finding such an object in the first
place.  Complete code and all associated data has been archived at
\cite{PaperData}, see Section~\ref{sec:code} for details.

\subsection{Open questions and next steps}

For interesting specific examples, open questions, and avenues for
further research, see Section~\ref{sec:QHSpheres},
Remarks~\ref{rem:closed}, \ref{rem:floersimple},
\ref{rem:moredisorder}, \ref{rem:allfoliar}, \ref{rem:hardtofol},
\ref{rem:tauteuler}, and \ref{rem:gauge},
Conjecture~\ref{conj:nonLpersist}, and Question~\ref{qu:persist}.

\subsection{History and acknowledgements} 

I began the first iteration of this project in 2004 not long after
\cite{KronheimerMrowkaOzsvathSzabo2007} appeared on the arXiv and have
worked on it on and off since, slowly increasing both the size of the
sample and the range of techniques, always seeking a counterexample to
what is now Conjecture~\ref{BGWconjecture}. Thus I cannot give here a
complete accounting of the debts I owe both individuals and
institutions on this project. Certainly, I gratefully thank Ian Agol,
John Berge, Danny Calegari, Marc Culler, Jake Rasmussen, Rachel
Roberts, Saul Schleimer, and Liam Watson for many helpful
conversations and ideas.  This work was done at Caltech, the
University of Illinois, ICERM, the University of Melbourne, and IAS,
and was funded in part by the Sloan Foundation, the Simons Foundation,
and the US National Science Foundation, the latter most recently by
the GEAR Network (DMS-1107452), DMS-1510204, and DMS-1811156.
Finally, I thank the referees for their helpful comments.

\section{Terminology and conventions} \label{sec:terms}

In this paper, all \3-manifolds will be orientable, all foliations
co-orientable, and all group actions on manifolds will be orientation
preserving. A $\Q$-homology \3-sphere is a closed \3-manifold whose
rational homology is the same as that of $S^3$.  A $\Q$-homology solid
torus is a compact \3-manifold $M$ with boundary a torus where
$H_*(M; \Q) \cong H_*(D^2 \times S^1; \Q)$; this is equivalent to $M$
being the exterior of a knot in some $\Q$-homology \3-sphere. One
defines $\Z$-homology \3-spheres and solid tori analogously.

Suppose $M$ is a compact \3-manifold with $\partial M$ a torus.  A
\emph{slope} on $\partial M$ is an unoriented isotopy class of simple
closed curve, or equivalently a primitive element of
$H_1(\partial M; \Z)$ modulo sign.  The set of all slopes will be
denoted $\slopes(M)$, which can be viewed as the rational points in
the projective line $P^1\big(H_1(\partial M; \R)\big) \cong
P^1(\R)$. The Dehn fillings of $M$ are parameterized by
$\alpha \in \slopes(M)$, with $M(\alpha)$ being the Dehn filling where
$\alpha$ bounds a disk in the attached solid torus.  When the interior
of $M$ is hyperbolic, Thurston showed that all but finitely many
$M(\alpha)$ are also hyperbolic \cite{ThurstonsNotes}.  The
nonhyperbolic Dehn fillings are called \emph{exceptional}, and the
corresponding slopes the \emph{exceptional slopes}.

A group is called \emph{left-orderable} when it admits a total
ordering that is invariant under left multiplication (see
\cite{ClayRolfsen2015} for an introduction to the role of orderable
groups in topology).  We will say that a closed $3$-manifold $Y$ is
\emph{orderable} when $\pi_1(Y)$ is left-orderable. By convention, the
trivial group is not left-orderable, and so $S^3$ is not orderable.

Heegaard Floer homology will always have coefficients in
$\F_2 = \Z/2\Z$. Recall from \cite{OSLensSpace2005} that an
\emph{L-space} is a $\Q$-homology 3-sphere with minimal Heegaard Floer
homology, specifically one where $\dim \HFhat(Y) = |H_1(Y;\Z)|$.

\section{Details on the sample}

\subsection{Some rational homology solid tori} 

Burton proved there are precisely 61{,}911 cusped finite-volume
hyperbolic \3-manifolds that have ideal triangulations with at most 9
tetrahedra \cite{Burton2014}.  Here, the equivalence is up to your
choice of homeomorphism, diffeomorphism, or isometry, provided
orientation reversing maps are allowed.  Each such manifold is the
interior of a compact manifold whose boundary is a nonempty union of
tori.  Of these compact manifolds, some 59{,}068 (95.4\%) are
$\Q$-homology solid tori, and I will denote this collection of
manifolds as $\cC$.

\begin{figure}
  \begin{center}
    \begin{tikzpicture}
      \matrix [matrix of nodes, column sep=-1em, row sep=-1em]{
 \pgfkeys{/matplotlibfigure, default, }%
 \input{plots/closed_volume}
 &
 \pgfkeys{/matplotlibfigure, default, }%
 \input{plots/closed_systole}
 \\
 \pgfkeys{/matplotlibfigure, default, }%
 \input{plots/closed_homology}
 & 
 \pgfkeys{/matplotlibfigure, default, }%
 \input{plots/closed_CS}
 \\
      };
    \end{tikzpicture}
  \end{center}
  \caption{Some basic geometric and topological statistics about the
    manifolds in $\cY$ given as histograms.  For similar plots about
    $\cC$, see \cite{DunfieldExceptional}.}
  \label{fig:basicstats}
\end{figure}
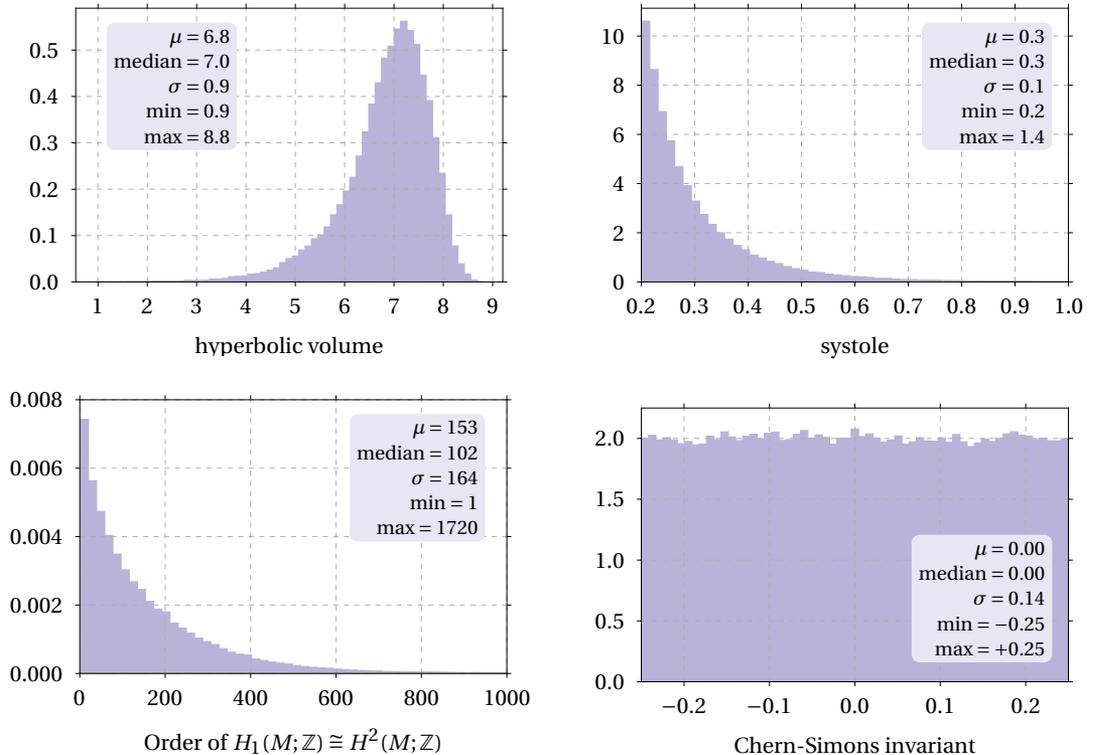

\subsection{The sample of rational homology spheres}
\label{sec:QHSpheres}
The manifolds in $\cY$ in Theorem~\ref{thm:main} are certain
$\Q$-homology \3-spheres that are Dehn fillings on $\cC$. While I
conjecture that they are precisely the Dehn fillings on $\cC$ that
give hyperbolic $\Q$-homology spheres whose systole is at least 0.2,
here I only establish the weaker Theorem~\ref{thm:arehyp} that they
are all hyperbolic and distinct.  (A complete list of all
non-hyperbolic fillings on $\cC$ is known, see
\cite{DunfieldExceptional}.)  Some basic statistics about the
manifolds of $\cY$ are shown in Figure~\ref{fig:basicstats}, and I now
turn to the proof of Theorem~\ref{thm:arehyp}.

\begin{proof}[Proof of Theorem~\ref{thm:arehyp}]
To prove each manifold in \QHSpheres is hyperbolic, I used the method
of \cite{HIKMOT2016}, as reimplemented by Goerner in \cite{SnapPy};
see Section~\ref{sec:intervals} of this paper for an overview of this
technique.  For this, as in the proof of Theorem 5.2 of
\cite{HIKMOT2016}, it was sometimes necessarily to search around for a
triangulation that could be used to certify the existence of a
hyperbolic structure.

To prove that the manifolds of \QHSpheres\ are distinct, I looked at
finite quotients of their fundamental groups.  Specifically, for
$Y \in \mathQHSpheres$, consider $G = \pi_1(M)$.  For a subgroup
$H \leq G$ of finite index, associate the tuple
\[
\big([G:H],\  [G:N],\  H^{\ab},\ N^\ab \big)
\]
where $N = \Core(H)$ is the largest normal subgroup of $G$ contained
in $H$, and $H^\ab$ is the abelianization of $H$.  The set of such
tuples associated to subgroups of index at most $n$ is an invariant of
$G$ and hence of $Y$.  I used coset enumeration \cite[Chapter
5]{HoltEickOBrien2005} to compute this invariant for $n = 6$, which
sufficed to distinguish all but 408 pairs of the manifolds in
\QHSpheres.  Those remaining pairs were separated by looking at the
corresponding invariant for all $H \triangleleft G$ where $G/H$ is
simple and $[G:H] < 10{,}000$.  All of these group-theoretic
calculations were done with Magma \cite{Magma}.
\end{proof}

\begin{remark}
  The above proof of distinctness implies that all manifolds in $\cY$
  have distinct profinite completions, which is an important open
  question for hyperbolic \3-manifolds generally, see \cite[Question
  1]{Agol2014}.  Also, Gardam \cite{Gardam2018} used a similar
  approach to the above to distinguish various census manifolds,
  including the 10,903 manifolds in $\cY$ from the original
  Hodgson-Weeks census.  One difference from \cite{Gardam2018} is that
  I looked at the abelianization of the core ($N^\ab$) as well as that
  of the subgroup itself ($H^\ab)$, allowing me to use $H$ of smaller
  index than \cite{Gardam2018}.  The expensive part of this technique
  is finding the subgroups $H$ rather than computing their homology,
  and, as taking the core is cheap, including $N^\ab$ provides a major
  speedup.
\end{remark}

\begin{remark}
  \label{rem:closed}
  For finite-volume hyperbolic 3-manifolds with cusps, one can
  rigorously compute the canonical Epstein-Penner decomposition and
  use that to prove that two manifolds are \emph{not} isometric
  \cite{DunfieldHoffmanLicata2015}.  However, for closed manifolds, to
  rigorously implement the procedure of \cite{HodgsonWeeks1994} one
  needs to provably find and drill out the shortest closed
  geodesic, thus reducing the problem to the cusped case. An important
  open question is whether one can use the certified hyperbolic
  structure produced via \cite{HIKMOT2016} to produce a provably
  correct shortest geodesic; see \cite{Trnkova2017} for some work in
  this direction.
\end{remark}

%% file: plots/diagram.tex
\begin{tikzpicture}
      [every node/.style={align=center},
       font=\small,
       main/.style={ellipse, draw=black!50, line width=1pt, outer sep=5pt},
       other/.style={inner sep=8pt}
      ]
      \node[main] (taut) at (0,0) {$Y$ has a taut \\ foliation $\cF$};
      \node[main] (orderable) at (9.5, 6)
          {$Y$ is orderable \\ $\Longleftrightarrow$ 
           \\ $\pi_1(Y)$ acts on $\R$}; 
      \node[main] (nonL) at (0, 6) {$Y$ is not \\ an $L$-space}; 

      \node[other] (circle) at (9.5, 0) {$\pi_1(Y)$ acts \\ on $S^1$};

      \node[other, inner sep=3pt] (leaf) at (4.75, 3) 
         {$\pi_1(Y)$ acts on a simply \\ connected 
        1-manifold \\ (possibly non-Hausdorff)};
    
      \begin{scope}[
        line width=1pt,
        double distance=4pt,
      ]
      \draw[-implies, double, bend right=5]    
           (taut.east) to node[below=4pt] 
           {Thurston's universal circle \cite{CalegariDunfield2003}}
           (circle.west);
      \draw[-implies, double] ([xshift=-15pt]taut.north) 
         to [bend left=10] 
         node[left=3pt] {
           \cite{OSgenusbounds2004} \\ 
             and\\ 
           \cite{Bowden2016, KazezRoberts2017}}
         ([xshift=-15pt]nonL.south);
      \draw[-implies, double, bend left=10] 
         (taut.north east) to 
         node[below right] {Leaf space of \\  $\cFtil$ in $\Ytil$}
         (leaf.south west); 

      \draw[-implies, double, bend left=10] ([xshift=15pt]nonL.south)
         to node[above right=8pt]{\cite{Juhasz2015}}
      ([xshift=15pt]taut.north);

      \draw[line width=10pt, loosely dashed, color=white, bend left=10]
         ([xshift=15pt]nonL.south) to ([xshift=15pt, yshift=9pt]taut.north);
      
      \draw[implies-implies, double, bend left=5] 
          (nonL.east) to node[above=3pt] 
          {Conjecture of \cite{BoyerGordonWatson2013}}
          (orderable.west);
     \draw[line width=8pt, loosely dashed, color=white, bend left=5] 
         ([xshift=9pt]nonL.east) to ([xshift=-9pt]orderable.west);

      \draw[-implies, double, bend left=5] (orderable.south west) to
         (leaf.north east);

      \draw[-implies, double, bend left=7] 
         ([xshift=10pt]orderable.south) to 
         node[left=3pt] {$\R \cong S^1 - \{\mbox{pt}\}$} 
         ([xshift=10pt]circle.north);
      \end{scope}
\end{tikzpicture}

%% file: plots/summary.tex
\begin{tikzpicture}[nmdstd, scale=0.9]
  \coordinate (O) at (0, 0);
  \coordinate (X) at (10, 0);
  \coordinate (Y) at (0, 6.18033988749895);
  \coordinate (XY) at ($(X) + (Y)$);
  \coordinate (Lsp) at (4.69599470141002, 0);
  \coordinate (NOrd) at (0, 4.75329532814139);
  \coordinate (Ord) at (0, 3.04547521627703);
  \coordinate (FgapX) at (-0.25, 0);
  \coordinate (FgapY) at (0, 0.531965581656767);
  
  \begin{scope}[
    rounded corners=0.2cm,
    line width=1pt,
    hatch distance=13pt,
    hatch thickness=3pt,
    ]
    \draw[preaction={fill=red!7.5},
          pattern=flexible hatch north west,
          hatch color=red!20]
          (O) rectangle ($(Lsp) + (Y)$);

    \draw (Lsp) rectangle (XY);

    \draw[fill=blue!40, fill opacity=0.3] (O) rectangle ($(Lsp) + (NOrd)$);
    
    \draw[preaction={fill=red!15},
          pattern=flexible hatch north east,
          hatch color=red!25]
           (Lsp) -- (X)[rounded corners=0.1cm] -- ++($(Ord) - (0, 0.1)$)
           -- ++(FgapX) -- ++(FgapY) --
           ++($-1*(FgapX)$) [rounded corners=0.2cm]
           -- (XY) -- ($(Lsp) + (Y)$) -- cycle;

    \draw[fill=blue!70, fill opacity=0.3]  (Lsp) rectangle ($(X) + (Ord)$); 
  \end{scope}

  
  \node[above] at ($(O)!0.5!(Lsp) + (Y)$) {L-spaces (47.0\%)};
  \node at ($(O)!0.5!(Lsp) + 0.5*(NOrd)$) {nonorderable $\geq 36.1$\%};
  \node at ($(O)!0.5!(Lsp) + ($(NOrd)!0.5!(Y)$)$) {no taut foliation $\geq 47.0$\%};
  \node[above] at ($(Lsp)!0.5!(X) + (Y)$) {non-L-spaces (53.0\%)};
  \node at ($(Lsp)!0.5!(X) + 0.5*(Ord)$) {orderable $\geq 26.1$\%};
  \node at ($(Lsp)!0.5!(X) + 1.5*(Ord)$) {taut foliation $\geq 52.8$\%};
  
\end{tikzpicture}

%% file: plots/closed_volume.tex
\begin{tikzoverlay}[width=\matplotlibfigurewidth]{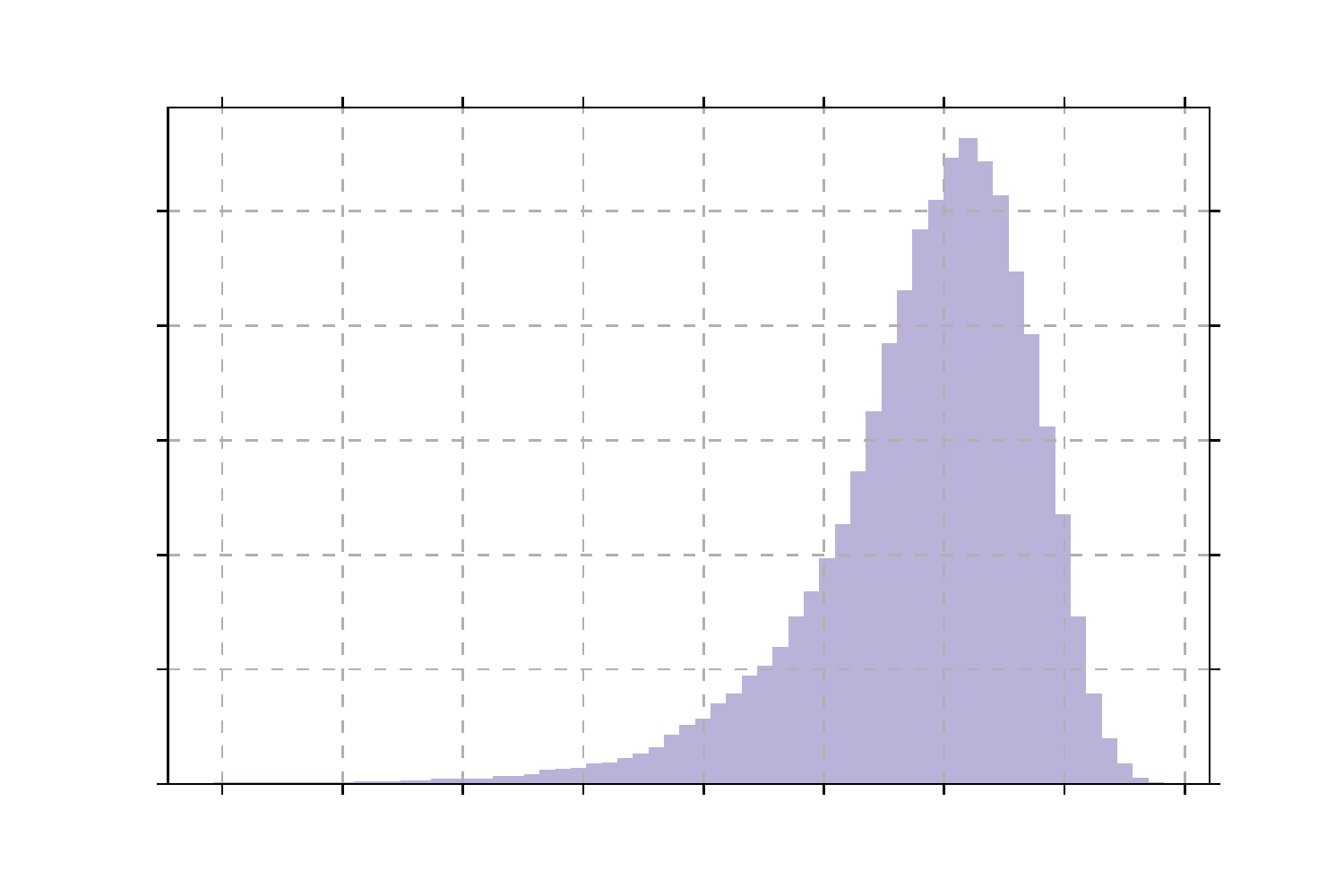}[\matplotlibfigurefont]
  \draw (41.875000, 56.150000)
  node[plot legend]
   {$\mu=6.8$ \\
    $\mathrm{median}=7.0$ \\
    $\sigma=0.9$ \\
    $\min=0.9$ \\
    $\max=8.8$};
  \draw (51.250000, 3.548177) node[below=1em] {hyperbolic volume};
  \draw (16.535651, 6.712963) node[below] {$1$};
  \draw (25.488349, 6.712963) node[below] {$2$};
  \draw (34.441048, 6.712963) node[below] {$3$};
  \draw (43.393746, 6.712963) node[below] {$4$};
  \draw (52.346444, 6.712963) node[below] {$5$};
  \draw (61.299143, 6.712963) node[below] {$6$};
  \draw (70.251841, 6.712963) node[below] {$7$};
  \draw (79.204539, 6.712963) node[below] {$8$};
  \draw (88.157238, 6.712963) node[below] {$9$};
  \draw (10.879630, 8.333333) node[left] {$0.0$};
  \draw (10.879630, 16.859337) node[left] {$0.1$};
  \draw (10.879630, 25.385341) node[left] {$0.2$};
  \draw (10.879630, 33.911345) node[left] {$0.3$};
  \draw (10.879630, 42.437349) node[left] {$0.4$};
  \draw (10.879630, 50.963353) node[left] {$0.5$};
  \begin{scope}[shift={(7.58295262, 8.33333333)},
                xscale=8.95269835, yscale=85.26003941]
  \end{scope}
\end{tikzoverlay}

%% file: plots/closed_systole.tex
\begin{tikzoverlay}[width=\matplotlibfigurewidth]{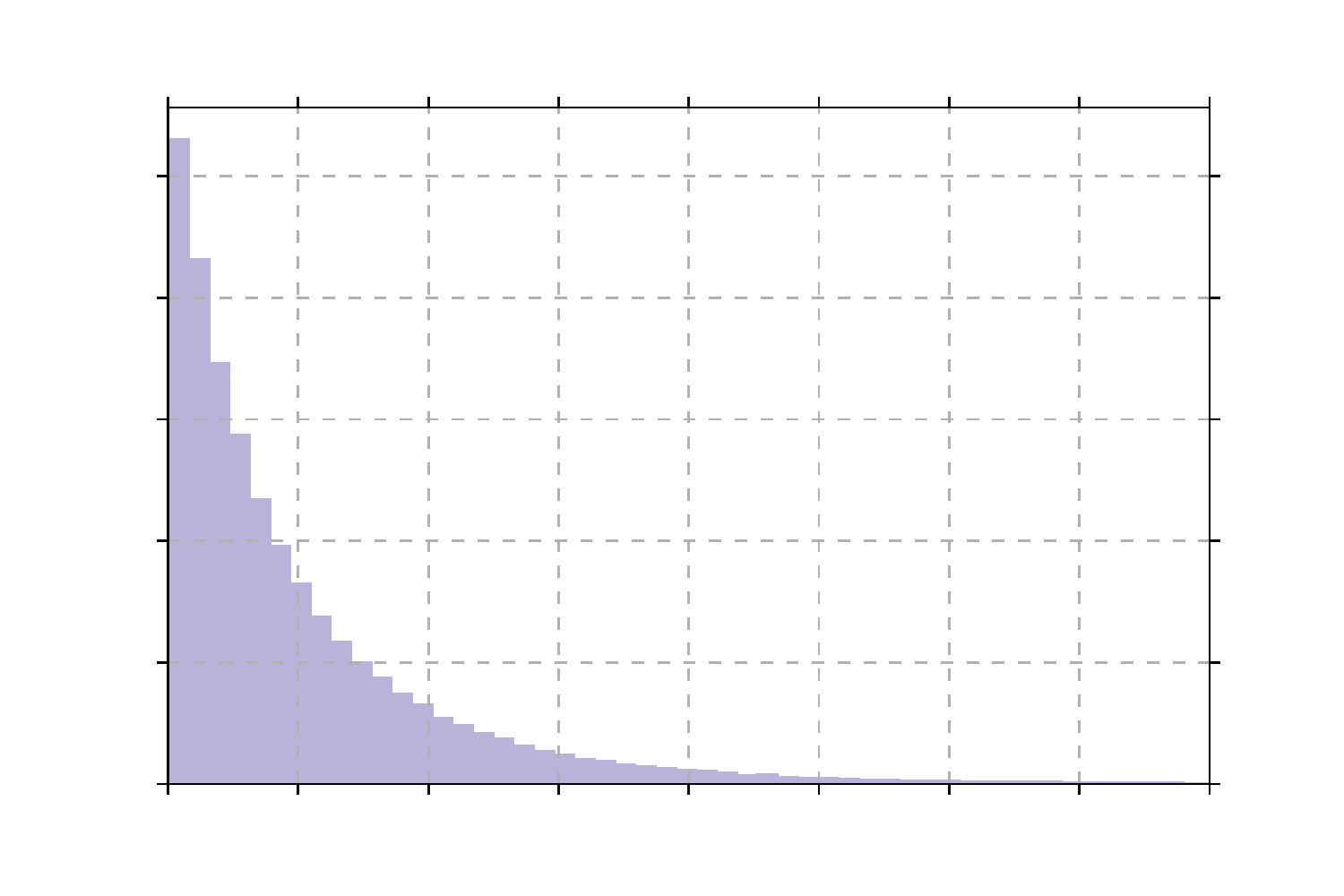}[\matplotlibfigurefont]
  \draw (87.125000, 56.150000) node[plot legend]
  {$\mu=0.3$ \\
   $\mathrm{median}=0.3$ \\
   $\sigma=0.1$ \\
   $\min=0.2$ \\
   $\max=1.4$};
  \draw (51.250000, 3.548177) node[below=1em] {systole};
  \draw (12.500000, 6.712963) node[below] {$0.2$};
  \draw (22.187500, 6.712963) node[below] {$0.3$};
  \draw (31.875000, 6.712963) node[below] {$0.4$};
  \draw (41.562500, 6.712963) node[below] {$0.5$};
  \draw (51.250000, 6.712963) node[below] {$0.6$};
  \draw (60.937500, 6.712963) node[below] {$0.7$};
  \draw (70.625000, 6.712963) node[below] {$0.8$};
  \draw (80.312500, 6.712963) node[below] {$0.9$};
  \draw (90.000000, 6.712963) node[below] {$1.0$};
  \draw (10.879630, 8.333333) node[left] {$0$};
  \draw (10.879630, 17.377214) node[left] {$2$};
  \draw (10.879630, 26.421094) node[left] {$4$};
  \draw (10.879630, 35.464975) node[left] {$6$};
  \draw (10.879630, 44.508855) node[left] {$8$};
  \draw (10.879630, 53.552736) node[left] {$10$};
  \begin{scope}[shift={(-6.87500000, 8.33333333)},
                xscale=96.87500000, yscale=4.52194026]
  \end{scope}
\end{tikzoverlay}

%% file: plots/closed_homology.tex
\begin{tikzoverlay}[width=\matplotlibfigurewidth]{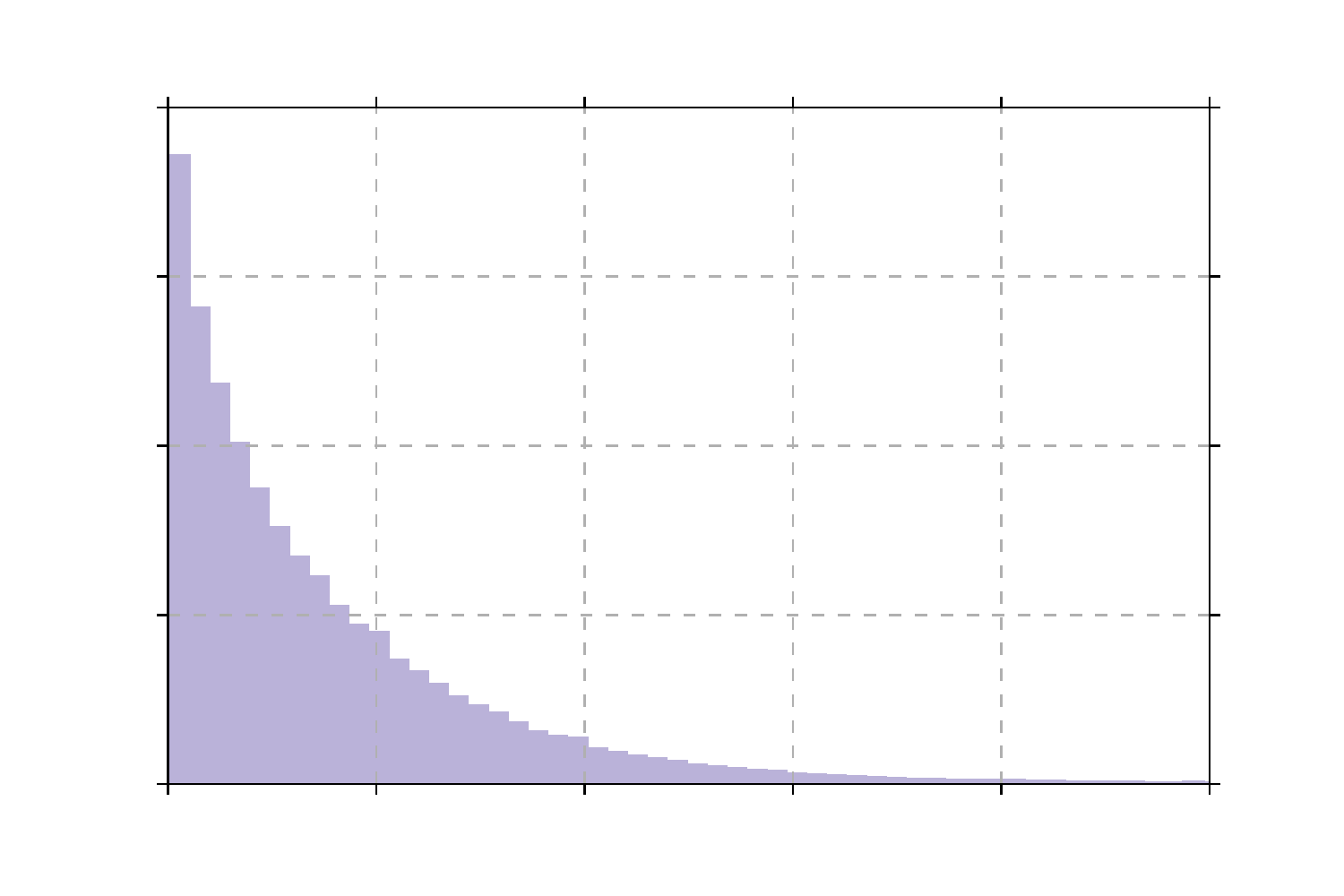}[\matplotlibfigurefont]
  \draw (86.125000, 56.150000) node[plot legend]
  {$\mu=153$ \\
   $\mathrm{median}=102$ \\
   $\sigma=164$ \\
   $\min=1$ \\
   $\max=1720$};
  \draw (51.250000, 3.548177) node[below=1em] {Order of $H_1(M; \Z)
    \cong H^2(M; \Z)$};
  \draw (12.500000, 6.712963) node[below] {$0$};
  \draw (28.000000, 6.712963) node[below] {$200$};
  \draw (43.500000, 6.712963) node[below] {$400$};
  \draw (59.000000, 6.712963) node[below] {$600$};
  \draw (74.500000, 6.712963) node[below] {$800$};
  \draw (90.000000, 6.712963) node[below] {$1000$};
  \draw (10.879630, 8.333333) node[left] {$0.000$};
  \draw (10.879630, 20.916667) node[left] {$0.002$};
  \draw (10.879630, 33.500000) node[left] {$0.004$};
  \draw (10.879630, 46.083333) node[left] {$0.006$};
  \draw (10.879630, 58.666667) node[left] {$0.008$};
  \begin{scope}[shift={(12.50000000, 8.33333333)},
                xscale=0.07750000, yscale=6291.66666667]
  \end{scope}
\end{tikzoverlay}

%% file: plots/closed_CS.tex
\begin{tikzoverlay}[width=\matplotlibfigurewidth]{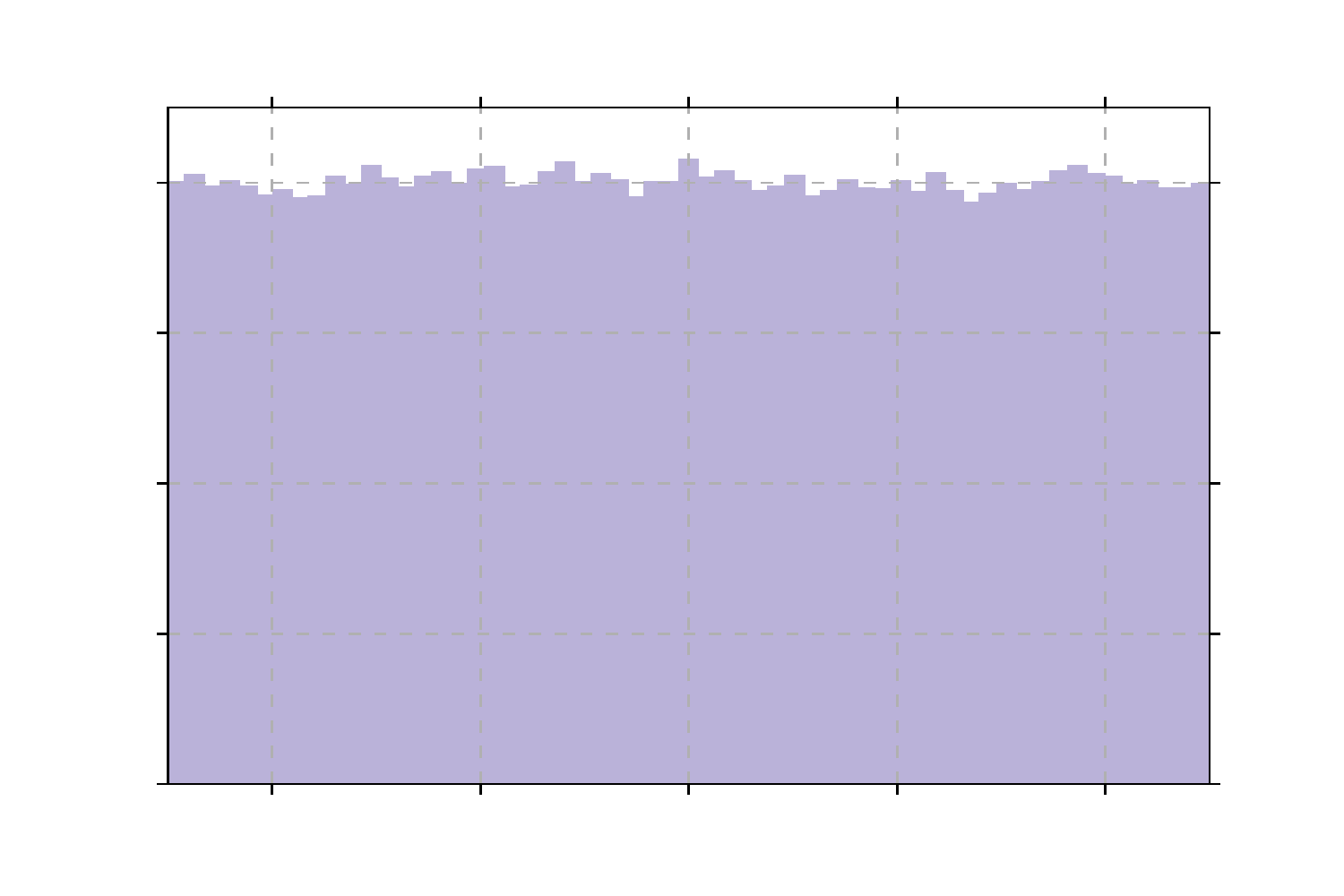}[\matplotlibfigurefont]
  \draw (87.125000, 35.50000) node[plot legend]
  {$\mu=0.00$ \\
   $\mathrm{median}=0.00$ \\
   $\sigma=0.14$ \\
   $\min=-0.25$ \\
   $\max=+0.25$};
  \draw (51.250000, 3.548177) node[below=1em] {Chern-Simons invariant};
  \draw (20.250000, 6.712963) node[below] {$-0.2$};
  \draw (35.750000, 6.712963) node[below] {$-0.1$};
  \draw (51.250000, 6.712963) node[below] {$0.0$};
  \draw (66.750000, 6.712963) node[below] {$0.1$};
  \draw (82.250000, 6.712963) node[below] {$0.2$};
  \draw (10.879630, 8.333333) node[left] {$0.0$};
  \draw (10.879630, 19.518519) node[left] {$0.5$};
  \draw (10.879630, 30.703704) node[left] {$1.0$};
  \draw (10.879630, 41.888889) node[left] {$1.5$};
  \draw (10.879630, 53.074074) node[left] {$2.0$};
  \begin{scope}[shift={(51.25000000, 8.33333333)},
                xscale=155.00000000, yscale=22.37037037]
  \end{scope}
\end{tikzoverlay}

%% file: floer.tex
\section{Finding the L-spaces}
\label{sec:compute-HF}

The simplest Heegaard Floer homology group $\HFhat$ is algorithmically
computable by \cite{SarkarWang2010} and hence the property of being an
L-space is algorithmically decidable.  Moreover, the bordered Heegaard
Floer theory of \cite{LipshitzOzsvathThurston2011,
  LipshitzOzsvathThurston2014} provides powerful and effective
computational tools for computing $\HFhat$.  This theory has been
implemented in e.g.~\cite{BFHPython} and has been successfully applied
to manifolds of the complexity of those in \QHSpheres, though one must
first find a Heegaard splitting of the input manifold specified in
terms of certain preferred generators of the mapping class group.
However, rather than attack the problem of computing $\HFhat$ head on,
I chose to use a bootstrapping procedure that exploited the structure
of the big Dehn Surgery graph \cite{HoffmanWalsh2015} via the results
of Rasmussen and Rasmussen \cite{Rasmussens2015} on L-space Dehn
fillings.  Throughout this section, complete code and manifold lists
are available at \cite{PaperData}.

\subsection{Floer simple manifolds} So that I can explain my
procedure, I first outline some results from
\cite{Rasmussens2015}. Let $M$ be a compact \3-manifold with $\bdry M$
a torus, and recall from Section~\ref{sec:terms} that $\slopes(M)$ is
the set of slopes on $\partial M$.  Since we are interested in which Dehn
fillings are L-spaces, define
\[
\cL(M) = \setdef{\alpha \in \slopes(M)}{\mbox{$M(\alpha)$ is an L-space}}
\]
Such an $M$ is called \emph{Floer simple} when it has \emph{two}
distinct Dehn fillings that are L-spaces, i.e.~$\abs{\cL(M)} > 1$,
see \cite[Proposition 1.3]{Rasmussens2015}.  For Floer simple
manifolds, the set $\cL(M)$ has the following structure; here, the
Tureav torsion of $M$ is a power series $\tau(M) \in \Z[[t]]$ which is
only slightly more complicated than the Alexander polynomial to
compute.
\begin{theorem}[{\cite[Theorem 1.6]{Rasmussens2015}}]
  \label{thm:Linterval}
  Suppose $M$ is a 3-manifold with $\partial M$ a torus. If $M$ is
  Floer simple, then $\cL(M)$ is either a closed interval or consists
  of every slope except the homological longitude. If you know that
  two slopes $\alpha \neq \beta$ are in $\cL(M)$, then $\cL(M)$
  can be explicitly computed from $\alpha$, $\beta$, and the Turaev
  torsion of $M$.
\end{theorem}
I leveraged this
result to determine the L-spaces in $\cY$ by identifying a large
number of manifolds in $\cC$ as either Floer simple or not Floer
simple.  This was done by an inductive bootstrapping procedure that
increased the level of knowledge about $\cY$ and $\cC$ in tandem. I
will describe this in detail below, but first I will explain the
starting point with regards to which manifolds in $\cC$ are Floer
simple.

\subsection{Priming the pump} \label{sec:prime}

Define $M$ to be \emph{Turaev simple} when every coefficient of
$\tau(M) \in \Z[[t]]$ is either 0 or 1. A basic obstruction to being
Floer simple is:
\begin{proposition}[{\cite[Prop.~1.4]{Rasmussens2015}}]
  \label{prop:turaev}
  A Floer simple manifold $M$ is Turaev simple.
\end{proposition}
Computing $\tau(M)$ for all the manifolds in $\cC$ identified
7{,}895 of them as not Turaev simple and hence not Floer simple.
To identify some initial Floer simple manifolds, I looked at
\emph{finite Dehn fillings}, that is slopes $\alpha$ where
$\pi_1\big(M(\alpha)\big)$ is finite. Since any \3-manifold with
finite fundamental group is an L-space, such an $\alpha$ is in
$\cL(M)$.  An immediate consequence of the data in
Theorem~1.2 of \cite{DunfieldExceptional} is thus:
\begin{corollary}\label{cor:finite}
  There are exactly 59{,}200 finite Dehn fillings on manifolds in
  $\cC$, with 78.2\% having at least one such filling. There are
  11{,}594 manifolds in $\cC$ with two or more finite fillings, all of
  which are therefore Floer simple.
\end{corollary}

\input tables/bootstrap.tex

\subsection{Bootstrapping procedure}

For 59.0\% of the manifolds in $\cY$, I am aware of only a single
description as a Dehn filling on something in $\cC$.  However, the
remaining 41.0\% average 3.4 known descriptions, and this
will be the key to my procedure for expanding what we know about
$\cY$ and $\cC$ in tandem.  The scheme is based on the following
allowed deductions, where here ``$M$ is Floer simple'' really means
``$M$ is Floer simple with at least two elements in $\cL(M)$
explicitly known''.
\begin{enumerate}[label=D\arabic*:, ref=D\arabic*]
\item\label{item:def1}
  If $M \in \cC$ has two Dehn fillings that are L-spaces, then
  $M$ is Floer simple. (This is just the definition.)

\item\label{item:RR}
  If $M \in \cC$ is Floer simple, then
  Theorem~\ref{thm:Linterval} determines exactly which of its Dehn
  fillings in $\cY$ are L-spaces.
  
\item\label{item:def2} If $M \in \cC$ is not Floer simple and we know
  $M(\alpha)$ is an L-space, then every other manifold in $\cY$ which
  is a Dehn filling on $M$ is not an L-space.  (This is also just the
  definition.)
\end{enumerate}
Table~\ref{table:bootstrap} shows what happens when using these three
deductions starting from the data in Section~\ref{sec:prime}; after
repeated applications, one arrives at the ``initial fixed point''
where only 7.4\% of the manifolds in \QHSpheres\ have unknown L-space
status.

I now describe some more sophisticated deductions, for which I need to
say a little more about \cite{Rasmussens2015}; for each Tureav simple
$M$, the authors define from $\tau(M)$ a subset $\iotainvDtau{M}$ in
$\slopes(M)$ which is either empty or infinite with a single limit
point, namely the homological longitude $\lambda$.  The precise
statement of \cite[Theorem~1.6]{Rasmussens2015} is as follows. For a
Floer simple $M$, if $\iotainvDtau{M}$ is empty, then
$\cL(M) = \slopes(M) \setminus \{\lambda\}$; otherwise, if
$\alpha \neq \beta$ are in $\cL(M)$, then $\cL(M)$ is the unique
closed interval with consecutive end points in $\iotainvDtau{M}$
containing both $\alpha$ and $\beta$.  This allows for the following
additional deductions when $M$ is Turaev simple but may or may not be
Floer simple.

\begin{enumerate}[label=D\arabic*:, ref=D\arabic*, start=4]
\item\label{item:Dempty}
  If $\iotainvDtau{M}$ is empty and $M$ has a non-L-space filling then
  $M$ is not Floer simple. 

\item\label{item:Fsimp} If $\iotainvDtau{M}$ is nonempty and
  $\alpha \in \cL(M)$, should $M$ be Floer simple there are only one
  or two possibilities for $\cL(M)$; there is one when $\alpha$ is not
  in $\iotainvDtau{M}$ and two when it is.  If all possibilities for
  $\cL(M)$ contain a known non-L-space filling, then we can conclude
  $M$ is not Floer simple.

\item\label{item:trick} As in \ref{item:Fsimp}, suppose
  $\iotainvDtau{M}$ is nonempty and $\alpha \in \cL(M)$.  As in
  \ref{item:Fsimp}, should $M$ be Floer simple there are at most two
  possibilities $P_i$ for $\cL(M)$.  If a $P_i$ contains a known
  non-L-space slope, it can be eliminated. Then any
  $\beta \in \slopes(M)$ not in the union of the remaining $P_i$ must
  be a non-L-space slope, even though we don't know whether or not $M$
  is Floer simple.
\end{enumerate}
Applying all six deductions reduces the number of \QHSpheres\ with
unknown L-space status from 7.4\% to 6.5\%; see
Table~\ref{table:bootstrap}.

\subsection{Endgame}

From \cite{DunfieldExceptional}, we know there are exactly 201{,}798
exceptional Dehn fillings on manifolds in $\cC$ which are
$\Q$-homology 3-spheres, and so far I have only used the 59,200 that
give spherical manifolds. One moreover has:
\begin{theorem}\label{thm:extraexcepts}
  Of the 199{,}662 exceptional $\Q$-homology sphere Dehn fillings on
  $\cC$ that do not have a hyperbolic piece in their JSJ
  decompositions, exactly 181{,}317 are L-spaces and 18{,}345 are not.
\end{theorem}

\begin{proof}
  As per Table 2 and Section 4.5 of \cite{DunfieldExceptional}, the
  201{,}798 exceptional $\Q$-homology sphere fillings consist of
  59,200 spherical manifolds, 4{,}296 connected sums of spherical
  manifolds, 72{,}841 Seifert fibered manifolds with infinite $\pi_1$,
  63{,}325 proper graph manifolds, and 2{,}136 manifolds with a
  non-trivial JSJ decomposition with a hyperbolic piece.  All
  spherical manifolds are L-spaces as are their connected sums
  since the connected sum of two L-spaces is again an L-space.  For
  everything except the ones with a hyperbolic piece, it is possible
  to compute $\HFhat$ directly as follows.  For each manifold, I
  translated from Regina's \cite{Regina} description of the graph
  manifold given in \cite{DunfieldExceptional} over to the weighted
  tree description of \cite{Neumann1981}. In that form, I used
  Hanselman's program \cite{HFhatGraph} associated to
  \cite{Hanselman2016} to compute $\HFhat$.
\end{proof}

Beyond the spherical fillings which I already used,
Theorem~\ref{thm:extraexcepts} provides an additional 140{,}462
fillings on $\cC$ whose L-space status is known.  Some 67{,}612 of
these are fillings on manifolds in $\cC$ that are already known to be
Floer simple, so only 72{,}850 of these fillings provide new
information, though the cases where we have two ways of determining
whether an exceptional Dehn filling is an L-space give a strong check
on the correctness of the computation.  Combining
Theorem~\ref{thm:extraexcepts} with repeated applications of the six
deductions results in only 6{,}437 (2.1\%) of the manifolds in
\QHSpheres\ having unknown L-space status.  It will turn out that only
3 are L-spaces and the other 6{,}434 are non-L-spaces.

By Theorem~\ref{thm:main}(\ref{item:fol}), some 162{,}341 of the
manifolds in \QHSpheres\ have taut foliations and hence are not
L-spaces; using this takes care of all but 3 of the manifolds in
\QHSpheres.  Moreover, the 8{,}115 persistent foliar orientations of
Theorem~\ref{thm:persistex} tell us that 172 additional manifolds in
$\cC$ are not Floer simple.

The three remaining manifolds in \QHSpheres\ are all Dehn fillings on
$o9_{34146}$, which is Turaev simple and has a single known L-space
filling, namely $o9_{34146}(1, 0)$ is the lens space $L(35,11)$.  I
claim that $o9_{34146} (0, 1)$ is also an L-space; in this case, the
manifold $o9_{34146}$ is then Floer simple and we can finish off the
last three manifolds in \QHSpheres.  The filling $o9_{34146} (0, 1)$
is hyperbolic but it is not in \QHSpheres\ because its systole is
$\approx 0.08648$ which is less than $0.2$.  SnapPy confirms that
$o9_{34146}(0, 1)$ is homeomorphic to $m007(11, 3)$.  Now $m007$ is
already known to be Floer simple and using \ref{item:RR} gives that
$m007(11, 3) = o9_{34146} (0, 1)$ is an L-space.  This completes the
proof of Theorem~\ref{thm:main}(\ref{item:Lsp}).  For code and further
details see \cite{PaperData}.

\begin{remark}\label{rem:floersimple}
  In the end, we know at least 50{,}598 (85.7\%) of the manifolds in
  $\cC$ are Floer simple and at least 8{,}352 (14.1\%) are not Floer
  simple but there are 118 (0.2\%) whose status is unknown.  The first
  ten unknown ones are: $t08191$, $t08263$, $o9_{10045}$,
  $o9_{18999}$, $o9_{19314}$, $o9_{19325}$, $o9_{19344}$,
  $o9_{19372}$, $o9_{19424}$, and $o9_{19478}$.
\end{remark}

%% file: tables/bootstrap.tex
\begin{table}
  \begin{center}\small
  \begin{tabular}{rrrrrrl}
    \toprule
    \multicolumn{3}{c}{$\cY$: $\Q$-hom.~\3-spheres} 
    &
    \multicolumn{3}{c}{$\cC$: $\Q$-hom.~solid tori} \\ 
    L-sp & non-L & L-sp? & F-simp & non-F & simp? & \\
    \midrule
    0 & 0 & 100 & 0 & 0 & 100 & initial state \\ 
    0 & 0 & 100 & 0 & 13.4 & 86.6 &  Proposition~\ref{prop:turaev}\\ 
    0 & 24.3 & 75.7 & 0 & 13.4 & 86.6 & \ref{item:def2} with Corollary~\ref{cor:finite}\\ 
    0 & 24.3 & 75.7 & 19.6 & 13.4 & 67.0 & \ref{item:def1} with Corollary~\ref{cor:finite}\\
   
    7.7 & 26.1 & 66.2 & 19.6 & 13.4 & 67.0 &  $\cY \Longleftarrow \cC$
                                             via \ref{item:RR} and \ref{item:def2} \\ 

    7.7 & 26.1 & 66.2 & 45.0 & 13.4 & 41.7 &  $\cY \Longrightarrow
                                             \cC$ via \ref{item:def1}\\
    40.1 & 40.2 & 19.7 & 45.0 & 13.4 & 41.7 &  $\cY \Longleftarrow
                                              \cC$ via \ref{item:RR} and \ref{item:def2}  \\
    
    40.1 & 40.2 & 19.7 & 51.1 & 13.4 & 35.5 &  $\cY \Longrightarrow
                                              \cC$ via \ref{item:def1}  \\
    46.8 & 45.6 & 7.6 & 51.1 & 13.4 & 35.5 &  $\cY \Longleftarrow \cC$ via \ref{item:RR} and \ref{item:def2}\\

    \midrule

    46.9 & 45.8 & 7.4 & 51.2 & 13.4 & 35.4 &  initial fixed point\\

    46.9 & 46.6 & 6.5 & 51.2 & 13.7 & 35.1 &  fixed point of
                                             \ref{item:def1}--\ref{item:trick} \\

    \midrule
    47.0 & 46.6 & 6.4 & 85.7 & 13.7 & 0.6 & \ref{item:def1} and
                                            \ref{item:RR} with Thm.~\ref{thm:extraexcepts} \\

    47.0 & 50.9 & 2.1 & 85.7 & 13.8 & 0.5 & fixed point of
                                             \ref{item:def1}--\ref{item:trick} \\

    47.0 & 53.0 & 0.001 & 85.7 & 14.1 & 0.2 & foliations \\

    \midrule
    47.0 & 53.0 & 0 & 85.7 & 14.1 & 0.2 & Last trick and final
                                          answer\\

    \bottomrule
  \end{tabular}
\end{center}
\caption{This table illustrates the steps in the proof of
  Theorem~\ref{thm:main}(\ref{item:Lsp}).  The first three columns
  record the percentages of the manifolds of \QHSpheres\ that are
  known to be L-spaces, known not to be L-spaces, and whose L-space
  status is unknown, respectively. The next three columns similarly
  record what is known about the manifolds in $\cC$ being Floer simple.
  At the beginning, we are completely ignorant about which manifolds
  in $\cY$ are L-spaces and which manifolds in $\cC$ are Floer
  simple. Then we apply the tools indicated in the rightmost column to
  learn more and more about both collections of manifolds. It takes
  five applications of each of \ref{item:def1}, \ref{item:RR}, and
  \ref{item:def2} before the data stabilizes and becomes
  self-consistent, arriving at the row labelled ``initial fixed
  point''; only the first three applications are shown individually as
  the others are indistinguishable at the level of rounding used in
  this table. }
  \label{table:bootstrap}
\end{table}

%% file: disorder.tex
\section{Proving groups are not orderable}
\label{sec:disorder}

This section is devoted to proving
Theorem~\ref{thm:main}(\ref{item:nonorder}), which is restated here:
\begin{theorem}\label{thm:nonorder}
  At least 110{,}940 of the manifolds in \QHSpheres\ are not orderable.
\end{theorem}
The basic approach follows \cite[\S 8--10]{CalegariDunfield2003}, with
the most significant change being how the word problem is solved in
the relevant \3-manifold groups.  This is discussed in
Section~\ref{sec:wordprob} below and is the key for going from
proving 44 manifolds are not orderable in \cite{CalegariDunfield2003}
to more than 100{,}000 manifolds here. The precise list of the
manifolds in Theorem~\ref{thm:nonorder} is available at
\cite{PaperData} along with the code and
additional data needed to prove it.

\begin{figure}
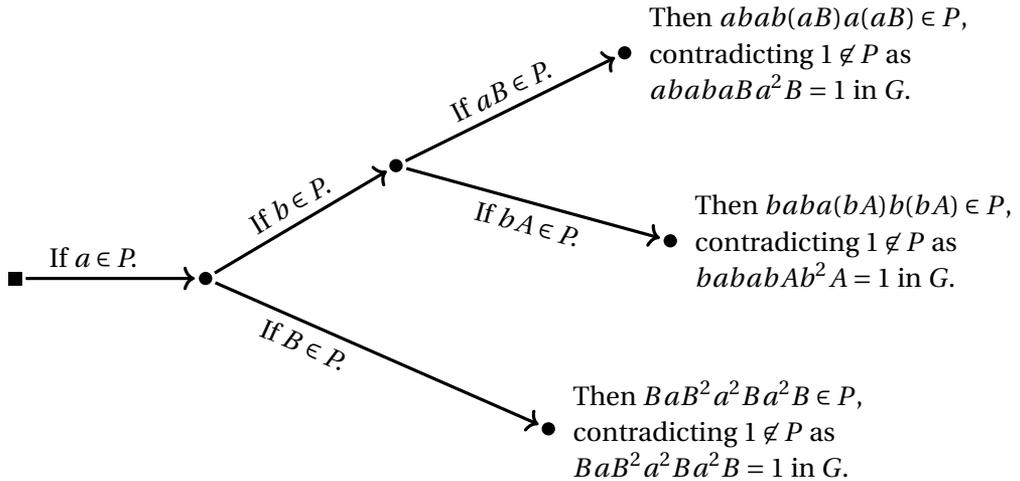

  \begin{center}
    \input plots/prooftree
  \end{center}

  \caption{This figure illustrates the proof of Theorem 9.1 of
    \cite{CalegariDunfield2003}, namely that the group
    $G = \big\langle a, b \ \big| \ ababaBa^2B,\, ababAb^2Ab
    \big\rangle$ is not left-orderable; here $A = a^{-1}$ and
    $B = b^{-1}$. (The group $G$ is the fundamental group of the Weeks
    manifold.) The nonordering tree should be read starting from
    the leftmost vertex and encodes a proof by contradiction: we
    assume that a positive cone $P$ exists, and then consider all
    possibilities for whether certain nontrivial elements of $G$ are
    or are not in $P$.  In each case, the contradiction comes from
    showing that $P \cdot P \subset P$ implies that $1 \in P$.
    Because we can reverse the roles of $P$ and $P^{-1}$, we need not
    consider the case when $A \in P$.}
  \label{fig:weeksproof}
\end{figure}

\subsection{Proof trees} 

For Theorem~\ref{thm:nonorder}, each manifold in the statement
was handled separately, though in a uniform manner as I
will now describe. Specifically, for each one I found a proof of
nonorderability that has a structure implicit in many arguments
that a group is not left-orderable, including those in
\cite{CalegariDunfield2003}.  To formalize this kind of proof, I first
need to fix the general context.  Let $G$ be a group with a fixed
finite generating set $S$.  Suppose further we have a solution to the
word problem for $G$, that is, an algorithm that determines whether or not a
word in $S$ corresponds to the trivial element in $G$.  If we have a
left-order on $G$, we consider its \emph{positive cone}
$P = \setdef{g \in G}{g > 1}$.  This gives a partition of $G$ into
$P \cup P^{-1} \cup \{1\}$ such that $P \cdot P \subset P$.
(Conversely, any such $P$ gives an associated order \cite[\S
1.4]{ClayRolfsen2015}.)  The prototype proof of nonorderability is given in
Figure~\ref{fig:weeksproof}, and the following formal definition is
most easily understood in the context of that example.

\begin{definition}
  A \emph{nonordering tree} for a group $G$ with
  generators $S$ is a finite trivalent tree $T$ with the following
  additional structure:
  \begin{enumerate}
    \item \label{item:root}
      The tree $T$ has a preferred root vertex, and all edges of
      $T$ are oriented pointing away from this root. 
      
    \item Each edge of $T$ is labeled by an element of
      $\FreeGroup(S)$.  At each interior vertex of $T$, the labels on
      the two outgoing edges are inverses in $\FreeGroup(S)$.

    \item Each leaf vertex, other than the root, is labeled by a word
      in the edge labels that appear along the unique directed path
      from the root to the leaf.

    \item \label{item:group}
      Every edge label corresponds to a nontrivial element of $G$,
      but every leaf label corresponds to $1$ in $G$. 
  \end{enumerate}
\end{definition}

Note that a nonordering tree is a finite combinatorial object, and
that if we can solve the word problem in $G$ then we can check
whether a given labeled tree satisfies
(\ref{item:root})--(\ref{item:group}).  It is straightforward to
generalize the thinking behind the example in
Figure~\ref{fig:weeksproof} to show:
\begin{theorem}\label{thm:nonorderingtree}
  Suppose $G$ has a nonordering tree $T$.  Then $G$ is not
  left-orderable. 
\end{theorem}
The converse to Theorem~\ref{thm:nonorderingtree} turns out to be true
as well, but I have no use for this fact here.  I will now outline the proof
of Theorem~\ref{thm:nonorder}.

\begin{proof}[Proof of Theorem~\ref{thm:nonorder}]
  For each of the manifolds in the statement, I used a heuristic
  method (described in Section~\ref{sec:findtree} below) to find a
  likely nonordering tree for its fundamental group.  I then certified
  that each labeled tree was in fact a nonordering tree using the
  solution to the word problem given in Section~\ref{sec:wordprob}.
  The largest tree was for the manifold $o9_{39416}(4, 1)$ which had
  some 20{,}329 leaves and 40{,}655 edges, though the median
  nonordering tree had only 12 leaves and 21 edges.  In total, there
  were 17.1 million instances of the word problem to solve and the
  longest word considered had length 76{,}196.  You can find the
  labeled tree for each manifold at \cite{PaperData}, along with the
  code used to verify them.  All together, the labeled trees weigh in
  at some 919MB of data, making this one very long proof.
\end{proof}

\subsection{Finding nonordering trees}\label{sec:findtree}

To find a nonordering tree for each manifold, I used the procedure
described in \cite[\S 8]{CalegariDunfield2003}, but using a fast,
although nonrigorous, method to ``solve'' the word problem.
Specifically, for a given $G = \pi_1(Y)$ with $Y \in \cY$, I used
SnapPy \cite{SnapPy} to find the images of a generating set $S$ of $G$
under an alleged holonomy representation
$\rho \maps \pi_1(Y) \to \SL{2}{\C}$ as matrices with double-precision
floating-point entries.  To test if a word in $S$ is $1$ in $G$, I
multiplied the corresponding matrices and checked if the result was
close to the identity matrix in $\SL{2}{\C}$.  While errors can and do
accumulate in this situation \cite{FloydWeberWeeks2002}, the words in
question typically had length less than 10, and the fact that we're
approximating a \emph{discrete} subgroup in $\SL{2}{\C}$ makes this quite
numerically robust in practice.  Indeed, not a single one of the
110{,}940 proofs found this way turned out to be incorrect when it was
rigorously verified using the solution to the word
problem of Section~\ref{sec:wordprob}.  I am grateful to Saul
Schleimer for suggesting attacking the word problem in this way,
rather than the original approach in \cite{CalegariDunfield2003} which
used the theory of automatic groups.

The procedure of $\cite[\S 8]{CalegariDunfield2003}$ is to look at a
ball $B_r$ of radius $r$ in the Cayley graph of $G$ with respect $S$,
and try to partition $B_r$ into $P \cup P^{-1} \cup \{1\}$ so that $P$
is closed under multiplications that stay in $B_r$.  (To give a sense
of scale, for most manifolds in Theorem~\ref{thm:nonorder}, I looked
at a $B_r$ with 5{,}000 to 30{,}000 elements.)  Because $G$ has
exponential growth, exactly which elements are in $B_r$ is very
sensitive to the choice of generating set $S$, even if we choose $r$
so that the number of elements in $B_r$ is roughly constant.  One
obvious choice, and the one used in \cite[\S
10]{CalegariDunfield2003}, is to minimize the size of $S$.  However, I
got much better results (that is, found more nonordering trees)
starting with presentations that had more than the minimal number of
generators but relatively short relators.  A typical example is this
presentation for $\pi_1\big(o9_{36382}(5,1)\big)$:
\[
\spandef{a, b, c, d, e, f}{
  \mathit{bdbbd},\, \mathit{adBf},\, \mathit{cAefB},
  \mathit{adFcACd}, \, \mathit{deeC},\, \mathit{aaaaaBcc}}
\]

\begin{remark}
  \label{rem:moredisorder}
  In searching for the proofs that formed Theorem~\ref{thm:nonorder},
  I always used the default Dehn surgery description of $Y$ and group
  presentation produced by SnapPy with the option
  \verb|minimize_number_of_generators=False|.  Working with a greater
  variety of descriptions and presentations would certainly show that
  several thousand additional manifolds are nonorderable.  For
  example, the very simple L-spaces $m006(3, 1)$ and $m016(6, 1)$ are
  not listed as nonorderable in \cite{PaperData}, though in fact both
  are.  The first is excluded because the technique of
  Section~\ref{sec:intervals} does not apply to this description (some
  shape $z$ has $\Im{z} < 0)$, but this can be fixed by using the
  description $m011(3, 1)$.  The second example can easily be handled
  by switching to the description $s002(-3, 2)$ for whatever reason.
\end{remark}

\section{Solving the word problem}\label{sec:wordprob}

Studying hyperbolic structures on \3-manifolds numerically has a
history going back 40 years to Riley \cite{Riley2013}, with much work
being done using Thurston's perspective of ideal triangulations and
gluing equations (see \cite{ThurstonsNotes, Weeks2005}) via the
program SnapPea and its successors \cite{SnapPy}.  While these methods
provide robust and compelling numerical evidence for the existence and
particulars of a hyperbolic structure for many a \3-manifold, they do
not prove that a given manifold is hyperbolic much less guarantee that
any hyperbolic invariants computed (e.g.~the volume) are approximately
correct. However, the breakthrough paper \cite{HIKMOT2016} shows how
many of these numerical computations can be made fully rigorous with
remarkably few changes using the framework of interval arithmetic and
interval analysis.  I now recall their method and explain how to use
it to rigorously solve the word problem for the fundamental group of a
given finite-volume hyperbolic \3-manifold.  While there are several
other ways to solve the word problem here, for example using the
Knuth-Bendix procedure to find and certify a short-lex automatic
structure \cite{WordProcessing1992}, the present approach is by far
the most efficient I know of and the only one fast enough for proving
Theorem~\ref{thm:nonorder}.

\subsection{Interval analysis}\label{sec:intervals}
In interval analysis, a number $z\in\C$ is partially specified by
giving a closed rectangle $\z$ with vertices in $\Q(i)$ that contains
$z$.  Because the vertices are rational, such intervals can be exactly
stored on a computer and rigorously combined by the operations
$+,-,\times,/$ to create other such intervals; see
e.g.~\cite{MooreKearfottCloud2009, Rump2010} for general background on
interval analysis and its applications.  One cost is that the sizes of
the rectangles grow with the number of operations. I will use $\IC$ to
denote the set of such complex intervals. Two elements of $\IC$ are
not viewed as equal even when all their vertices agree since each
represents some unknown point inside the rectangle.  In contrast, when
two elements of $\IC$ are disjoint they are considered not equal.

The key tools of interval analysis needed here are effective versions
of the Inverse Function Theorem.  For a suitable smooth function
$f \maps \C^n \to \C^n$ these results provide simple tests for proving
that $f$ has a zero in some interval vector $\z \in \IC^n$.  Suppose
$Y$ is a closed 3-manifold specified by Dehn fillings on a topological
ideal triangulation $\cT$.  If $\cT$ has $n$ tetrahedra, the
logarithmic form of Thurston's gluing equations give a smooth map
$f \maps \C^n \to \C^n$ where a zero $z$ of $f$ gives a hyperbolic
structure on $Y$ provided the imaginary parts of all $z_i$ are
positive; here each $z_i$ specfies the geometric shape in $\H^3$ of
one of the tetrahedra in $\cT$.  In favorable circumstances, one can
apply the interval Newton's method or Krawczyk's test to prove that an
interval vector $\z \in \IC^n$ contains a zero of $f$ of the needed
type.  In this case, the hyperbolic structure on $Y$ is known in terms
of the shapes of the tetrahedra of $\cT$ to within a tolerance
determined by the sizes of the rectangles $\z_i$. A $\z \in \IC^n$
that has been proven to contain such a zero of $f$ will be
called an \emph{approximate hyperbolic structure}.  The authors of
\cite{HIKMOT2016} demonstrated this approach works fantastically well
in practice, in particular certifying the existence of hyperbolic
structures on more than 27{,}000 manifolds.

For the computations here, I used Goerner's implementation of this
approach included in SnapPy \cite{SnapPy} when used inside SageMath
\cite{SageMath}.  The original code of \cite{HIKMOT2016} used elements
of $\IC$ whose vertices are machine precision floating point numbers.
This provides maximum speed but limits how small one can make the
rectangles $\z_i$ in an approximate hyperbolic structure; in contrast,
SnapPy allows arbitrary precision vertices.  While this is
significantly slower, for the manifolds here it still takes less than
a second to use the interval Newton's method to certify an approximate
hyperbolic structure where each $\z_i$ has diameter less than
$10^{-1000}$.  Using the interval implementation of the dilogarithm
function provided by Arb
\cite{Johansson2017}, such an approximate
hyperbolic structure allows SnapPy to give the volume of one of these
manifolds to 1{,}000 provably correct digits, though it takes an
additional 5--10 seconds to do so.

\subsection{Approximate holonomy representations}
\label{sec:approxholo}
  
Let $M_2(\IC)$ denote 2-by-2 matrices with entries in $\IC$.  The
arithmetic operations on $\IC$ naturally give addition and
multiplication operations on $M_2(\IC)$.  Let $\GL{2}{(\IC)}$ be the
elements of $M_2(\IC)$ whose determinant, itself an element of $\IC$,
does not contain $0$. While many products of elements of
$\GL{2}{(\IC)}$ are again in $\GL{2}{(\IC)}$, it is not actually
closed under multiplication.  Suppose $G$ is the fundamental group of
a hyperbolic \3-manifold $Y$ with a finite symmetric generating set
$S$. Any holonomy representation from $G$ to
$\Isom^+(\H^3) = \PSL{2}{\C}$ of the hyperbolic structure on $Y$ lifts
to a representation $G \to \SL{2}{\C}$ which I will still call a
holonomy representation.  An \emph{approximate holonomy
  representation} is a function $\brho \maps S \to \GL{2}{(\IC)}$ such
that there is a holonomy representation
$\rho_0 \maps G \to \SL{2}{\C}$ so that
$\rho_0(\gamma) \in \brho(\gamma)$ for all $\gamma \in S$. Given a
word $w = s_1 s_2 \cdots s_n$ in $S$, define
$\brho(w) = \rho(s_1) \rho(s_2) \cdots \rho(s_n)$ as you would expect.
Regardless of whether $\brho(w)$ is in $\GL{2}{(\IC)}$ or just
$M_2(\IC)$, we are guaranteed that $\rho_0(w) \in \brho(w)$.

In the setup of Section~\ref{sec:intervals}, suppose a finite-volume
hyperbolic \3-manifold $Y$ is specified by an approximate hyperbolic
structure $\z$ on an ideal triangulation $\cT$.  Since we know the
approximate shapes of the geometric ideal tetrahedra in $\cT$, we can
start tiling out to build an approximation of the developing map
between the universal cover of $\cT$ and $\H^3$ and use that to
compute an approximate holonomy representation; see e.g.~\cite[Lemma
3.5]{DunfieldGaroufalidis2012} for details.

\subsection{Solving the word problem} Suppose $w$ is a word in our
generators $S$ of $G = \pi_1(Y)$ and we wish to know if $w = 1$ in
$G$.  If $\brho$ is an approximate holonomy representation and
$\brho(w)$ does not contain the identity matrix $I$, then we have
proved that $w \neq 1$ in $G$ since then there is an actual holonomy
representation $\rho_0$ with $\rho_0(w) \neq I$. However, if instead
$I \in \brho(w)$, we cannot immediately conclude $w = 1$ since
$\brho(w)$ might contain other elements of $\rho_0(G)$ as well.
Fortunately, the subgroup $\rho_0(G)$ is discrete in $\SL{2}{\C}$, so
if the entries of $\brho(w)$ have small enough diameter then
$I \in \brho(w)$ does imply $\rho_0(w) = I$ and hence $w = 1$.  This
can be made effective by exploiting J\o rgensen's inequality:

\begin{theorem}\label{thm:word}
  Suppose $\brho$ is an approximate holonomy representation for a
  hyperbolic 3-manifold $Y$ where $s_1$ and $s_2$ are in $S$ with $2
  \notin \tr\big(\brho([s_1, s_2])\big)$.  If $w$ is a word in $S$ with
  \begin{equation}\label{eqn:Jorg}
    \abs{\tr^2(\brho(w))
 - 4} 
    + \abs{\tr\big([\brho(s_i), \ \brho(w)]\big) - 2} < 1
  \end{equation}
  for each $i = 1, 2$ then $w = 1$ in $G = \pi_1(Y)$. 
\end{theorem}
Here, the lefthand side of (\ref{eqn:Jorg}) is a real interval and the
equation should be read as saying both of its endpoints are less than
$1$.

\begin{proof}[Proof of Theorem~\ref{thm:word}]
Let $\rho_0$ be a holonomy representation compatible with $\brho$. For
this proof, we will focus on the action of elements of $G$ on $\H^3$,
and so regard $\rho_0$ as a representation $G \to \PSL{2}{\C}$.  By
hypothesis, we have $\rho_0\big([s_1, s_2]\big) \neq I$ and so both
$\rho_0(s_i) \neq I$ and moreover $\pair{\rho_0(s_1), \rho_0(s_2)}$ is
a nonelementary torsion-free discrete subgroup of $\PSL{2}{\C}$.

For each $i$, the equation (\ref{eqn:Jorg}) implies that J\o rgensen's
inequality (see e.g.~\cite[Theorem 2.17]{MatsuzakiTaniguchi1998}) is
violated for the pair $\rho_0(w), \rho_0(s_i)$; in particular,
$\Gamma_i = \pair{\rho_0(w), \rho_0(s_i)}$ cannot be a nonelementary
discrete subgroup of $\PSL{2}{\C}$.  Since $\Gamma_i$ is discrete,
this means it must be elementary.  By
e.g.~\cite[Proposition~2.2]{MatsuzakiTaniguchi1998}, as $\Gamma_i$ is
torsion-free the nontrivial elements of $\Gamma_i$ are either all
parabolic with a common fixed point on $\partial \H^3$ or all
hyperbolic with a common axis in $\H^3$.  So if $\rho_0(w) \neq I$,
then $\pair{\rho_0(s_1), \rho_0(s_2)}$ would also be elementary, a
contradiction.  So $\rho_0(w) = I$ and hence $w = 1$ in $G$ as needed.
\end{proof}

\subsection{Practical considerations}
\label{sec:praccon}

Any given approximate holonomy representation $\brho$ will not solve
the word problem for all $w$ as the intervals making up $\brho(w)$
frequently ``smear out'' to the point where $I \in \brho(w)$ but
Theorem~\ref{thm:word} does not apply.  Thus to determine if a given
$w = 1$ in $G$, you may need to refine the approximate hyperbolic
structure to get one accurate enough so that the associated $\brho$
can determine whether $\rho_0(w) = 1$ or not.  In the proof of
Theorem~\ref{thm:nonorder}, for $78.0\%$ of the manifolds it sufficed
to use approximate hyperbolic structures that were correct to about
100 bits of precision (about 30 decimal digits, or roughly twice that
of an IEEE-754 double-precision floating point number), though
$22.0\%$ of the manifolds needed 200 bits for some words, and 11
manifolds ($0.01\%$) needed 400 bits. As mentioned above, the longest
word $w$ that was considered in Theorem~\ref{thm:nonorder} had length
$\abs{w} = 76{,}196$, which makes these precision requirements seem
modest by comparison.  In fact, this particular $w$ only required
$\brho$ to be computed to 200 bits, which resulted in the entries of
the $\brho(s)$ for $s \in S$ having diameter about $2^{-175}$, whereas
the entries of $\brho(w)$ had diameter about $2^{-31}$.

%% file: plots/prooftree.tex
\begin{tikzpicture}
  [font=\small,
   sloped, 
   main/.style={shape=circle, fill,
     inner sep=0pt, outer sep=1.5pt, 
     minimum size=5pt},
   line width=1.2pt,
   ]

  \node[main, shape=rectangle] (root) at (0, 0) {};
  \node[main] (a) at (2.5, 0) {};
  \node[main] (b) at (5, 1.5) {};
  \node[main] (c) at (7, -2) {};
  \node[main] (d) at (8.6, 0.5) {};
  \node[main] (e) at (8, 3) {};

  \node[right=5pt, align=left] at (e) {Then $abab(aB)a(aB) \in P$, \\
       contradicting $1 \not\in P$ as \\ $ababaBa^2B = 1$ in $G$.};

  \node[right=5pt, align=left] at (d) {Then $baba(bA)b(bA) \in P$, \\
       contradicting $1 \not\in P$ as \\ $bababAb^2A = 1$ in $G$.};

  \node[right=5pt, align=left] at (c) {Then $BaB^2a^2Ba^2B \in P$, \\
       contradicting $1 \not\in P$ as \\ $BaB^2a^2Ba^2B = 1$ in $G$.};

  \draw[->, above=4pt, pos=0.4] (root) to node {If $a \in P.$} (a); 
  \draw[->, below=4pt, pos=0.3] (a) to node {If $B \in P.$} (c); 
  \draw[->, above=4pt] (a) to node {If $b \in P.$} (b); 
  \draw[->, below=4pt] (b) to node {If $bA \in P.$} (d); 
  \draw[->, above=4pt] (b) to node {If $aB \in P.$} (e); 

\end{tikzpicture}

%% file: foliations.tex
\section{Foliar orientations}
\label{sec:foliar}

A \emph{triangulation} of a closed 3-manifold is a cell complex made
out of finitely many tetrahedra with all their \2-dimensional faces
identified in pairs via affine maps so that the link of every vertex
is a 2-sphere.  (For such face gluings, the link condition is
equivalent to the complex being a closed 3-manifold, see
e.g.~\cite[Prop.~3.2.7]{Thurston1997}.) In particular, a triangulation
is not necessarily a simplicial complex, but rather what is sometimes
called a semi\hyp simplicial or pseudo\hyp triangulation.

\begin{figure}
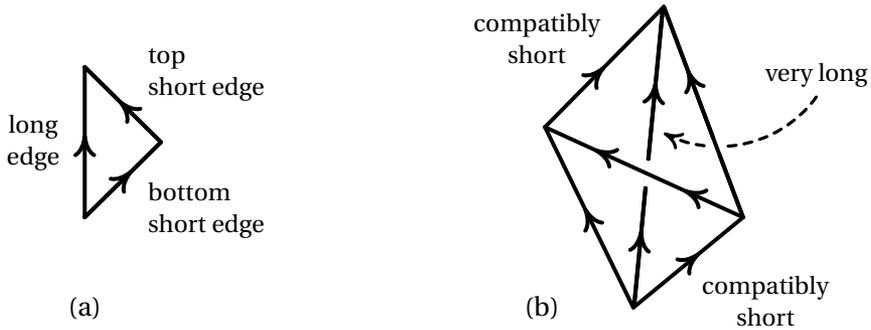

\begin{center}
  \input plots/local_edge_orientations
\end{center}
\vspace{-2ex}
\caption{Local possibilities for acyclic edge orientations.}
\label{fig:edges}
\end{figure}

An \emph{edge orientation} of a 3-manifold triangulation $\cT$ is a
choice of direction for each edge of the \1-skeleton $\cT^{1}$. An
edge orientation $\mu$ of $\cT$ is \emph{acyclic} when there is no
face of $\cT^2$ where $\mu$ orients its edges as a directed
cycle. Thus, in an acyclic edge orientation, each face of $\cT$ has
the orientation pattern shown in Figure~\ref{fig:edges}(a), and the
edges of the face are called \emph{long} or \emph{short} as indicated,
with the short edges further subdivided into \emph{bottom} and
\emph{top}.  A simple case check shows that any acyclic edge
orientation of a tetrahedron is as in Figure~\ref{fig:edges}(b), up to
a possibly orientation reversing homeomorphism.  (An acyclic edge
orientation gives an ordering of the vertices of each simplex of $\cT$
making it into a $\Delta$-complex as defined in \cite[\S
2.1]{Hatcher2002}, and indeed in our context these two notions are
equivalent.)  In each such tetrahedron, the unique edge that is long
in both of its adjacent faces is called \emph{very long}.  The two
edges labelled \emph{compatibly short} are those that are the same
kind of short edge in both of their adjacent faces.


Now suppose $\mu$ is an acyclic edge orientation of $\cT$. First, an
edge of $\cT$ is a \emph{sink edge} with respect to $\mu$ when it is
very long in every tetrahedron for which it is adjacent.  Second, the
\emph{face relation} of $\mu$ is the equivalence relation on the faces
of $\cT^{2}$ generated by the rule that for each compatibly short edge
in a tetrahedron the two faces adjacent to it are equivalent.  To
visualize the equivalence classes of the face relation, consider the
dual cellulation $\cD$ to $\cT$, where in particular the edges of
$\cD^1$ correspond to the faces of $\cT^2$. In each tetrahedron
$\sigma$ of $\cT$, slice $\cD^1 \cap \sigma$ into two pieces at the
vertex in $\cD^0$ via a quadrilateral that separates the two
compatibly short edges and is otherwise disjoint from $\cD^1$.  This
cuts $\cD^1$ into a collection of closed loops which correspond
precisely to the equivalence classes of the face relation.

\begin{figure}
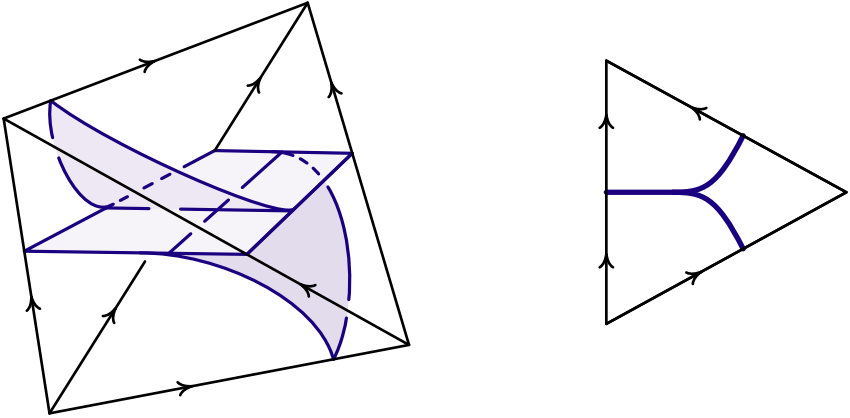

\begin{center}
  \input plots/branched
\end{center}
\vspace{-1ex}
\caption{The branched surface $B(\mu)$ in one tetrahedron.  Note that
  each face of the tetrahedron looks like the figure at right, and
  hence these local pictures glue up to give a branched surface in
  $\cT$.}
\label{fig:branched}
\end{figure}

I now focus on the case when $\cT$ has a single vertex. Putting these
concepts together, a \emph{foliar orientation} of $\cT$ is an
acyclic edge orientation that has no sink edges, and whose face
relation has a single equivalence class. The key result of this
section is:
\begin{theorem}\label{thm:foliar}
  Suppose $\cT$ is a \1-vertex triangulation of a closed orientable
  \3-manifold $Y$ that admits a foliar orientation $\mu$. Then
  $Y$ has a co-orientable taut foliation $\cF$ that is transverse to
  $\cT^1$ and induces the edge orientation $\mu$.
\end{theorem}
Having a foliar orientation is a strengthing of Calegari's notion of a
local orientation for a \3-manifold triangulation \cite{Calegari2000}.
Specifically, for a \1-vertex triangulation, it turns out that
Calegari's local orientation condition is equivalent to an acyclic
edge orientation whose face relation has a single equivalence class,
and thus the added restriction is that there are no sink edges.  The
basic idea for Theorem~\ref{thm:foliar} is to build a branched surface
that is ``dual'' to the foliar orientation as shown in
Figure~\ref{fig:branched}, and then apply Li's foundational work
on laminar branched surfaces \cite{Li2002} to prove the existence of a
related foliation.

\subsection{Branched surfaces}

Figure~\ref{fig:branched} shows how to associate a branched surface
$B(\mu)$ to any acyclic edge orientation $\mu$; note here that $\mu$
gives a global co-orientation to $B(\mu)$; since $Y$ is orientable,
this gives an global orientation to the branches of $B(\mu)$.  The key
property is that when $\mu$ is foliar then $B(\mu)$ is a laminar
branched surface as defined by Li \cite{Li2002}.

\begin{proof}[Proof of Theorem~\ref{thm:foliar}]

Within this proof, I will use concepts and notation from \cite[\S
0-1]{Li2002} without further comment and also abbreviate $B(\mu)$ to
$B$. To start, let us understand the topology of the guts
$G = Y \setminus \interior(N(B))$ where $N(B)$ is a regular
neighborhood of $B$.  As $\cT$ has a single vertex, we see from
Figure~\ref{fig:branched} that $G$ is a single \3-ball and to complete
the picture we need to understand how the horizontal boundary
$\partial_h N(B)$ and vertical boundary $\partial_v N(B)$ decompose
the sphere $\partial G$.  Each tetrahedron of $\cT$ contributes two
rectangles to $\partial_v N(B)$, and the gluing pattern of these
rectangles matches the face relation of $\mu$; since the face relation
has a single equivalence class, the vertical boundary
$\partial_v N(B)$ is a single annulus and so the horizontal boundary
$\partial_h N(B)$ is two discs.  Hence $G$ is just a
$D^2 \times I$ region.

I next show that $B$ is laminar.  As per Definition 1.4 of
\cite{Li2002}, there are four conditions to check.  The first two
conditions are (1) and (2) from Proposition 1.1 of \cite{Li2002},
namely that $\partial_h N(B)$ is incompressible in $G$, no component
of $\partial_h N(B)$ is a sphere, that $G$ is irreducible, and finally
that $G$ has no monogons; all of these are immediate since $G$ is a
$D^2 \times I$ region.  The third condition is (3) of Proposition 1.1
of \cite{Li2002}, namely that $B$ does not carry a torus that bounds a
solid torus in $Y$.  To see this, first note that any surface $S$
carried by $B$ inherits an orientation from $B$ and moreover $S$ has
strictly positive algebraic intersection number with each
$\mu$-oriented edge of $\cT$ that it meets.  Since $\cT$ has one
vertex, each edge is a closed loop, and so we conclude $S$ is
non-separating and in particular does not bound a solid torus as
required. The fourth and final condition to check is that, after
collapsing all trivial bubbles, the branched surface $B$ has no sink
discs where all cusps point inward.  As G is connected, there are no
trivial bubbles; looking back at Figure~\ref{fig:branched}, the only
local branch that could contribute to a sink disk for $B$ is the
square that meets the very long edge.  The hypothesis that $\cT$ has
no sink edges with respect to $\mu$ thus gives exactly that $B$ has no
sink discs.  Therefore, we have shown that $B$ is laminar.

Here is how to build the foliation promised in the theorem.  By
Theorem 1(a) of \cite{Li2002}, as $B$ is laminar it fully carries an
essential lamination. Since the guts of $B$ is a $D^2 \times I$
region, we can extend this lamination to a foliation $\cF$ of $Y$
that is transverse to $\cT^1$ and induces the orientation $\mu$
there. The foliation $\cF$ is taut since the edges of $\cT$ form a
collection of transverse closed loops meeting every leaf.
\end{proof}

\subsection{Examples} Foliar orientations are extremely common in
practice, and I was able to use them to prove 98.6\% of 
Theorem~\ref{thm:main}(\ref{item:fol}):
\begin{theorem}\label{thm:usefoliar}
  At least 160{,}003 of the manifolds in \QHSpheres\ have
  1-vertex triangulations admitting foliar orientations.
\end{theorem}
The rest of Theorem~\ref{thm:main}(\ref{item:fol}) will be handled
using a related technique introduced in Section~\ref{sec:persist}.  I
will now outline the proof of Theorem~\ref{thm:usefoliar}.

\begin{proof}[Proof of Theorem~\ref{thm:usefoliar}]
  For each manifold $Y$ in \QHSpheres, a variety of 1-vertex
  triangulations were generated by doing random Pachner moves starting
  from the manifold's various descriptions as Dehn fillings on
  manifolds in $\cC$.  For each such triangulation, an exhaustive
  search for foliar orientations was done using the method of
  Section~\ref{sec:foliarsearch}. The first such triangulation found
  and all of its foliar orientations were saved and can be found at
  \cite{PaperData} together with the relevant code.  The saved
  triangulations averaged 21.8 distinct foliar orientations each
  (modulo a complete reversal of signs); there were 1,215
  triangulations with a unique foliar orientation, and largest number
  of foliar orientations for a single triangulation was 884.
\end{proof}

\subsection{Triangulations with more vertices}

Here is one way to generalize Theorem~\ref{thm:foliar} to
triangulations with more than one vertex. Throughout, let $\cT$ be a
triangulation of a closed orientable \3-manifold $Y$.  First, an edge
orientation $\mu$ of $\cT$ is \emph{strongly connected} when the
directed graph $(\cT^1, \mu)$ is strongly connected, i.e.~every pair
of vertices $v_0$ and $v_1$ in $\cT^0$ are joined by a $\mu$-directed
path $\cT^1$ that starts at $v_0$ and ends at $v_1$. Second, I will
say that $\cT$ has \emph{short loops} if for every vertex $v$ of
$\cT^0$ there is an edge in $\cT^1$ which joins $v$ to itself.
Finally, an acyclic edge orientation $\mu$ of $\cT$ is \emph{foliar}
when $\cT$ has short loops and $\mu$ is strongly connected, has no sink
edges, and where the number of equivalence classes of its face
relation is equal to the size of $\cT^0$.  Note that when $\cT$ has
only one vertex this is equivalent to the original definition of
foliar.  Again, we have:
\begin{theorem}\label{thm:foliarmulti}
  Suppose $\cT$ is a triangulation of a closed orientable \3-manifold
  $Y$ that admits a foliar orientation $\mu$. Then $Y$ has
  a co-orientable taut foliation $\cF$ that is transverse to $\cT^1$ and
  induces the edge orientation $\mu$.
\end{theorem}
\begin{proof}[Proof of Theorem~\ref{thm:foliarmulti}]
  The proof is very similar to that of Theorem~\ref{thm:foliar}.  Let
  $B = B(\mu)$ be the branched surface associated to $\mu$. The
  condition on the face relation of $\mu$ means that the guts
  $Y \setminus \interior(N(B))$ is a union of $D^2 \times I$ regions,
  one for each vertex of $\cT$. Moreover, since $\cT$ has short loops
  none of these are trivial bubbles, which means there is nothing to
  collapse before checking for sink disks, and $B$ of course has none
  of these since $\mu$ lacks sink edges.  Finally, as $\mu$ is
  strongly connected, the argument from before shows that $B$ cannot
  carry a separating surface.  Thus $B$ is laminar and fully carries
  an essential lamination, which we can again fill in to get a
  foliation $\cF$ which will be taut since the needed family of
  transverse loops can be constructed using that $\mu$ is strongly
  connected.
\end{proof}

For the 3{,}000 non-L-spaces that are not foliated by
Theorem~\ref{thm:usefoliar}, I searched for 2-vertex triangulations to
which Theorem~\ref{thm:foliarmulti} applied, without a great deal of
success:
\begin{theorem}\label{thm:foliartwovert}
  There at least 174 manifolds in \QHSpheres\ not covered by
  Theorem~\ref{thm:usefoliar} with 2-vertex triangulations that admit
  a foliar orientation.
\end{theorem}
As always, the specific manifolds and source code are at
\cite{PaperData}. 

\begin{remark}\label{rem:allfoliar}
  I do not know if every co-orientable taut foliation $\cF$ of $Y$
  arises from a foliar orientation. Here are two approaches for trying
  to show this is the case.
  \begin{enumerate}
  \item Blow some air into $\cF$ to get an essential lamination of $Y$
    with nonempty complement, and then apply \cite{Li2002} to get a
    laminar branched surface $B$ that carries it. The guts of $B$ must
    be product regions, and if they are not $D^2 \times I$ then add
    some extra branches to $B$ to cut them down into such regions.
    After collapsing any trivial bubbles, try to build a cell complex
    that is dual to $B$ and modify it to one that is a triangulation.
  \item Alternatively, you could start with the kind of triangulation
    considered in \cite{Gabai2000} which is transverse to $\cF$ and
    supports a combinatorial volume preserving flow; these have some
    thematic similarities with foliar orientations.
  \end{enumerate}
  Unfortunately, it seems unlikely that either approach would give an
  algorithm for determining whether $Y$ has such a foliation, much
  less one that is actually implementable and practical.
\end{remark}

\subsection{Searching for foliar orientations}
\label{sec:foliarsearch}

Since the Euler characteristic of a closed \3-manifold is $0$, a
1-vertex triangulation $\cT$ with $n$ tetrahedra has $n+1$ edges.
Fixing the orientation of the first edge arbitrarily, there are $2^n$
edge orientations to search through when looking for a foliar
orientation.  As I considered triangulations with a median of $19$ and
a maximum of $27$ tetrahedra, this makes a completely naive search
infeasible.  I dealt with this by reframing the task as a boolean
satisfiability (SAT) problem and used a general-purpose SAT solver to
enumerate all acyclic edge orientations, from which the foliar ones
were easily selected. Although SAT is $\NP$-complete and enumerating
all solutions to a SAT problem is $\#\P$-complete, modern heuristic
solvers make short work of problems of the size and type that arise
here.

Here is the basic idea. Given a \1-vertex triangulation $\cT$,
arbitrarily fix some edge orientation $\mu_0$.  Any other edge
orientation $\mu$ is then encoded by using one boolean variable per
edge to record whether $\mu$ agrees with $\mu_0$ there. For each face
$\sigma$ of $\cT^2$ there is a boolean clause equivalent to the
orientation being acyclic there as follows.  Let $\epsilon_i$,
$\epsilon_j$, and $\epsilon_k$ be the variables associated to the
three edges of $\sigma$. If $\mu_0$ orients $\partial \sigma$ as a
cycle, then $\mu$ is acyclic if and only if
\[
  (\epsilon_i \lor \epsilon_j \lor \epsilon_k) \land 
  (\lnot \epsilon_i \lor \lnot\epsilon_j \lor \lnot\epsilon_k)
\]
If $\mu_0$ does not orient $\partial \sigma$ as a cycle, one can
always create a cycle by reversing one of the edges and use this to
formulate the correct clause.  For example, if reversing the
orientation on the edge associated to $\epsilon_i$ creates a cycle,
then $\mu$ is acyclic if and only if
\[
  (\lnot \epsilon_i \lor \epsilon_j \lor \epsilon_k) \land 
  (\epsilon_i \lor \lnot\epsilon_j \lor \lnot\epsilon_k)
\]
Thus finding all acyclic edge orientations is equivalent to finding
all solutions to a certain SAT problem, specifically a 3-SAT problem
with $n$ variables and $4n$ clauses in conjunctive normal form. 

I used CryptoMiniSat \cite{CryptoMiniSat} for the enumeration.  On a
sample of about $2{,}900$ triangulations with 22 tetrahedra each, the
average time to enumerate all acyclic edge orientations was $0.02$
seconds per triangulation, and the largest number of acyclic edge
orientations found was $1{,}138$, which is quite small compared to
$2^{22} \approx 4.2$ million. 

\begin{remark}
  The number of acyclic edge orientations can in fact be exponentially
  large in the size of $\cT$, which precludes a universal
  polynomial-time algorithm for enumerating them.  Here is one easy if
  somewhat degenerate way to see this.  Start with any $\cT$ that has
  an acyclic edge orientation $\mu$.  Pick any tetrahedron in $\cT$
  and do a $0 \to 2$ move on the pair of faces adjacent to its very
  long edge; that is, unglue those faces from their neighbors and
  insert a ``pillow'' consisting of two new tetrahedra that form the
  link of a new valence two edge $\tau$.  If $\cT'$ is the new
  triangulation, then $\mu$ can be extended to $\cT'$ by \emph{either}
  orientation of $\tau$.  Repeating $n$ times gives a triangulation
  with $2n + c$ tetrahedra and at least $2^n$ acyclic edge
  orientations.
\end{remark}

\section{Persistently foliar orientations}
\label{sec:persist}

I now describe an analog of foliar orientations for a manifold with
torus boundary which gives taut foliations on all but one of its Dehn
fillings. The approach is motivated by the constructions of taut
foliations on all non-trivial Dehn surgeries on certain knots in
$S^3$; see Section~\ref{sec:knots} below.  Specifically, let $M$ be an
orientable \3-manifold with $\partial M$ a single torus. An
\emph{ideal triangulation} of $M$ is a cell complex $\cT$ obtained
from finitely many tetrahedra by identifying their faces in pairs so
that $\cT^0$ is a single point whose complement is homeomorphic to
$M \setminus \partial M$.  The last condition is equivalent to the
complement of a small open regular neighborhood of $\cT^0$ being
homeomorphic to $M$.  Henceforth, I will view $M$ as embedded in $\cT$
in this way, and as such $M$ inherits a cell structure as a union of
truncated tetrahedra.

Just as in the closed case, an acyclic edge orientation $\mu$ for
$\cT$ has an associated branched surface $B(\tau)$ defined by
Figure~\ref{fig:branched}, and I will always choose $B(\tau)$ so that
it lies inside $M$.  If $N(B(\tau))$ is the standard regular
neighborhood of $B(\tau)$ as in \cite[\S0--1]{Li2002}, then
$M \setminus \interior(N(B(\tau)))$ is homeomorphic to
$\partial M \times [0, 1]$ and the vertical boundary
$\partial_vN(B(\tau))$ forms a collection of annuli in
$\partial M \times \{1\}$.  The acyclic edge orientation $\mu$ is
\emph{persistently foliar} when all of these annuli are essential in
$\partial M$ and $\mu$ has no sink edges.  When $\mu$ is persistently
foliar, the annuli making up $\partial_vN(B(\tau))$ must all be
parallel in $\partial M \times \{1\}$ and the common unoriented
isotopy class of their core curves is called the \emph{degeneracy
  slope} of $\mu$ and denoted $\degen(\mu)$.  The main result of this
section is:

\begin{theorem}\label{thm:persist}
  Suppose $M$ is a compact orientable \3-manifold with $\partial M$ a
  torus with an ideal triangulation $\cT$.  If $\mu$ is a persistently
  foliar orientation for $\cT$, then every Dehn filling $M(\alpha)$ of
  $M$ with $\alpha \neq \degen(\mu)$ has a co-orientable taut
  foliation.
\end{theorem}

\begin{proof}[Proof of Theorem~\ref{thm:persist}]
The argument will closely parallel the one given for
Theorem~\ref{thm:foliar}.  From now on, view $B = B(\mu)$ as a
branched surface in $M(\alpha)$.  I will first show that $B$ is
laminar.  The guts $G = M(\alpha) \setminus \interior(N(B))$ are a
solid torus with both $\partial_hN(B)$ and $\partial_vN(B)$ being
parallel annuli on $\partial G$ which are essential as
$\alpha \neq \degen(\mu)$. This immediately gives the first two
criteria for $B$ to be laminar, namely conditions (1) and (2) from
Proposition 1.1 of \cite{Li2002}. The reason that $B$ does not carry a
surface that bounds a solid torus is essentially the same as in
Theorem~\ref{thm:foliar}: every branch of $B$ has nonzero algebraic
intersection number with a closed loop in $M$, namely a loop made by
closing off a component of $\cT^1 \cap M$ by an arc in $\partial
M$. Finally, just as before the branched surface $B$ has no sink discs
since $\mu$ has no sink edges. So $B$ is laminar.

As $B$ is laminar, it carries an essential lamination by
\cite{Li2002}.  We can extend this lamination to a co-orientable
foliation $\cF$ of all of $M(\alpha)$ by using the usual ``stack of
monkey saddles'' or ``stack of chairs'' construction; see e.g.~Example
4.22 of \cite{Calegari2007}, and note we can extend the monkey saddles
into $N(B)$ to fill up the remaining space since $B$ has no discs of
contact by Corollary 2.3 of \cite{Li2002}. It remains to check that
$\cF$ is taut. It is not hard to extend each arc in $\cT^1 \cap M$
through the stack of monkey saddles to form a closed transverse loop,
and together with the core curve of $G$ we now have a collection of
closed transverse loops meeting every leaf.  Thus $M(\alpha)$ admits
the taut foliation claimed.
\end{proof}

\subsection{Exteriors of knots in the 3-sphere}
\label{sec:knots}

\begin{table}
  \begin{center}
    \small
    \begin{tabular}{rl}
      \toprule
      Knot & Description\\
      \midrule
      $8n3$ & Torus knot $T(3,4)$ \\ 
      $10n21$ & Torus knot $T(3, 5)$ \\
      $12n242$ & $(-2, 3, 7)$ pretzel knot \\
      $13n4587$ & $(2, 7)$ cable on $T(2, 3)$\\
      $13n4639$ & $(2, 5)$ cable on $T(2, 3)$\\
      $14n6022$ &  $(-2, 3, 9)$ pretzel knot\\
      $14n21881$ & Torus knot $T(3, 7)$\\
      $15n40211$ & $(2, 9)$ cable on $T(2, 3)$\\
      $15n41185$ & Torus knot $T(4, 5)$\\ 
      $15n124802$ & $(2, 3)$ cable on $T(2, 3)$\\ 
      $16n184868$ & $(-2, 3, 11)$ pretzel knot\\
      $16n783154$ & Torus knot $T(3, 8)$\\
      \bottomrule
  \end{tabular}
  \end{center}
  \caption{All 12 nonalternating L-space knots with at most 16
    crossings. Nomenclature follows
    \cite{HosteThistlethwaiteWeeks1998}. The above cables
    on the trefoil $T(2, 3)$ are L-space knots by Theorem~1.10 
    of \cite{Hedden2009}.  The pretzel knots are L-space knots by 
    Theorem 1 of \cite{BakerMoore2014}.
  }
  \label{table:lspace}
\end{table}

The idea of constructing an object in a knot exterior $M$ in $S^3$
that induces foliations in all but the trivial Dehn filling goes back
at least to \cite{GabaiKazez1990} and was used in \cite{Delman1995,
  Naimi1997, Brittenham2001, HirasawaKobayashi2001} among others. In
these references, the focus was on an essential lamination in $M$,
called a \emph{persistent lamination}, that remains essential in all
nontrivial Dehn fillings.  However, as in the proof of
Theorem~\ref{thm:persist}, in each Dehn filling one can fill up the
complement of the lamination with monkey saddles to get a foliation.
In the context of Theorem~\ref{thm:persist}, the persistent lamination
is carried by the branched surface $B(\tau)$.  More generally, one can
consider a certain class of \emph{persistently foliar branched
  surfaces}, of which these $B(\tau)$ are examples, that give taut
foliation on all nontrivial Dehn fillings of $M$.

A knot $K$ in $S^3$ has a nontrivial Dehn surgery that is an L-space
if and only if its exterior is Floer simple; knots with such surgeries
are called L-space knots.  Recently, Delman and Roberts have announced
that such persistently foliar branched surfaces exist for all
alternating knots and all Montesinos knots that have no $L$-space
surgeries \cite{DelmanRobertsTBD}.  Motivated by this work, I tried to
apply Theorem~\ref{thm:persist} to some nonalternating knots and
found:
\persistentlynonalternating
\noindent
The 12 exceptions are listed in Table~\ref{table:lspace}.  This may
seem like very few, but only 492 of these knots have Alexander
polynomials with the form required by \cite[Cor.~1.3]{OSLensSpace2005}
for being an L-space knot. The ideal triangulations used in
Theorem~\ref{thm:nonalt} had between 5 and 39 tetrahedra, with a mean
of 25.6; you can get all $1{,}210{,}596$ ideal triangulations with
persistent foliar orientations at \cite{PaperData} together with the
relevant source code.

Combined with the work of Delman and Roberts \cite{DelmanRobertsTBD},
Theorem~\ref{thm:nonalt} makes it safe to posit:
\begin{conjecture}\label{conj:nonLpersist}
  The exterior of a non-L-space knot in $S^3$ has a persistently
  foliar branched surface.  
\end{conjecture}

\subsection{Examples} As with the closed case, persistently foliar
orientations are quite common:
\begin{theorem}\label{thm:persistex}
  At least 8{,}115 of the manifolds in $\cC$ have ideal triangulations
  that admit persistently foliar orientations.
\end{theorem}
To put this in context, when $M$ has a persistently foliar
orientation, by Theorem~\ref{thm:persist} at most one Dehn filling of
$M$ can be an L-space and hence $M$ is not Floer simple; by
Remark~\ref{rem:floersimple}, there are between 8{,}352 and 8{,}470
manifolds in $\cC$ that are not Floer simple, so
Theorem~\ref{thm:persistex} covers at least 95.8\% of them.
Collectively, the foliar orientations in
Theorem~\ref{thm:persistex} give taut foliations on some 99{,}339
manifolds in \QHSpheres.

It is intriguing that every non-Floer simple knot exterior examined in
Theorem~\ref{thm:nonalt} has a persistently foliar orientation, and yet
there are non-Floer simple manifolds in $\cC$ where I am unable to
find such an orientation. This is particularly striking
in light of the fact that the manifolds in Theorem~\ref{thm:nonalt}
typically require many more tetrahedra to triangulate than those in
$\cC$. This suggests the answer to the following question might be
yes:
\begin{question}
  \label{qu:persist}
  Are there non-Floer simple manifolds that do not contain a persistently
  foliar branched surface?
\end{question}
The first few candidates in $\cC$ for such a manifold are $v2347$,
$v2626$, $v3452$, $t03447$, $t04027$, $t06287$, $t06400$, and
$t06953$.

\subsection{The manifold m137}

Gao showed in \cite{Gao2017} that the manifold $m137$ in $\cC$ has
infinitely many $\Z$-homology sphere fillings $Y$ where $\pi_1(Y)$ has
no nontrivial homomorphisms to $\PSL{2}{\R}$.  Thus, one cannot order
$\pi_1(Y)$ using the technique of $\PSLRtilde$ representations.
However, the standard 4-tetrahedra triangulation of $m137$ has a
persistently foliar orientation that shows that every Dehn filling
except the homological longitude has a taut foliation.  (The
longitudinal filling is $S^2 \times S^1$.)  In particular, each
$\Z$-homology sphere filling $Y$ has a taut foliation, and as
$H^2(Y; \Z) = 0$, Theorem~\ref{thm:euler_agree} from the next section
applies to show that $Y$ is orderable. So Gao's examples do fully
satisfy Conjecture~\ref{BGWconjecture}.

\subsection{Collecting foliations} I now turn to:
\begin{proof}[Proof of Theorem~\ref{thm:main}(\ref{item:fol})]
  By Theorems~\ref{thm:usefoliar} and~\ref{thm:foliartwovert}, there
  are at least 160{,}177 manifolds in $\cY$ with taut foliations whose
  existence can be certified using a foliar orientation of the
  individual manifold. In addition, as per the discussion immediately
  after Theorem~\ref{thm:persistex}, some 99{,}339 manifolds in $\cY$
  have taut foliations coming from persistently foliar orientations of
  manifolds in $\cC$.  Combined, these two methods cover 162{,}341
  manifolds in $\cY$, completing the proof; see~\cite{PaperData} for
  details.
\end{proof}

\begin{remark}
\label{rem:hardtofol}
Here is a sample of non-L-spaces where I was unable to find a taut foliation:
\begin{center}\footnotesize
\begin{tabular}{cccccc}
  $s137(5, 4)$    &   $s460(6, 1)$    & $s593(6, 1)$ &
  $s614(5, 1)$    &  $s753(-6, 1)$    & $s956(4, 1)$ \\
  $v1333(-5, 1)$  & $v3045(-4, 1)$    & $t06114(5, 1)$ &
  $t08155(-5, 1)$ & $o9_{12518}(-6, 1)$ & $o9_{12544}(-4, 3)$ \\
  $o9_{13679}(-6, 1)$ & $o9_{14675}(-1, 5)$ & $o9_{15066}(-5, 1)$ &
  $o9_{22743}(7, 1)$ & $o9_{30634}(6, 1)$ & $o9_{36699}(7, 1)$ \\
\end{tabular}
\end{center}
I also do not know whether or not any of these manifolds are
orderable.
\end{remark}
  
\section{Foliar orientations and the Euler class}
\label{sec:foleuler}

\begin{figure}
  \begin{center}
    \begin{minipage}[b]{0.45\textwidth}
      \begin{center}
        \input plots/mixed
      \end{center}
      \caption{The mixed edges of an acyclically oriented tetrahedron.}
      \label{fig:mixed}
    \end{minipage}
    \begin{minipage}[b]{0.5\textwidth}
      \begin{center}
        \input plots/tri_section
      \end{center}
      \caption{The section $\eta$ of $\UTBpunc$ in each face of $\cT$.}
      \label{fig:sect_on_tri}
    \end{minipage}
  \end{center}
\end{figure}

\begin{figure}
  \begin{center}
    \input plots/section
  \end{center}
  \caption{The section $\eta$ of $\UTBpunc$ in each tetrahedron of $\cT$.}
  \label{fig:section}
\end{figure}

\begin{figure}
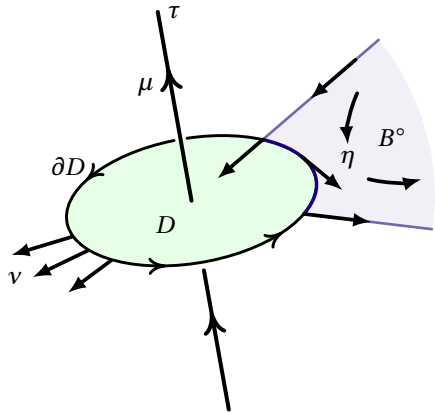

  \begin{center}
    \input plots/orientation
  \end{center}
  \caption{How $\eta$ twists with respect to $\nu$ near a
    mixed edge $\tau$.}
  \label{fig:twisting}
\end{figure}

The following is a consequence of Thurston's Universal Circle
construction:
\begin{theorem}[{\cite[Theorem 8.1]{BoyerHu2018}}]
  \label{thm:euler_agree}
  A $\Q$-homology \3-sphere $Y$ with a taut foliation $\cF$ whose
  Euler class $e(\cF) \in H^2(Y;\Z)$ is zero is orderable.
\end{theorem}
Roughly, Thurston associates to the taut foliation $\cF$ an action of
$\pi_1(Y)$ on a certain circle built from the circles at infinity
of the induced foliation of the universal cover of $Y$; the Euler
class of this action, which is the obstruction to lifting the action
to $\R$, turns out to be $e(\cF)$. See Sections 6--8 of
\cite{BoyerHu2018} for more.

Theorem~\ref{thm:euler_agree} gives a very useful tool for showing a
\3-manifold is orderable, and I now explain how to
compute the Euler class for the foliation $\cF$ associated to a
foliar orientation $\mu$.  In any acyclically oriented tetrahedron,
there are two edges that are short on one face but long on another; I
will call these the \emph{mixed} edges and they are shown in
Figure~\ref{fig:mixed}.  An edge $\tau \in \cT^1$ appears as an edge
in various of the tetrahedra in $\cT^3$, and let $\mixed_\mu \tau$ be
the number of times it appears as a mixed edge with respect to
$\mu$. (As $\cT$ is not simplicial, the edge $\tau$ can appear several
times in the same tetrahedron, and so a single tetrahedron could
contribute as much as $2$ to $\mixed_\mu \tau$.)

Consider the dual cell complex $\cD$ to
the triangulation $\cT$ and orient the 2-cells of $\cD$ by the
orientation $\mu$ gives to the corresponding edge of $\cT$. Define a
cochain $\phi_\mu \in C^2(\cD; \Q)$ by 
\[
\phi_\mu(D) = 1 - \frac{1}{2}\mixed_\mu(\tau)
  \mtext{where $\tau \in \cT^1$ is the edge dual to $D \in \cD^2$.}
\]
I will show:
  
\begin{theorem}\label{thm:cocycle}
  Suppose $\cF$ is the foliation of $Y$ associated to a foliar
  orientation $\mu$ of a 1-vertex triangulation $\cT$.  Then the
  cochain $\phi_\mu$ is a cocycle lying in $C^2(\cD; \Z)$ and
  $[\phi_\mu] = e(\cF)$ in $H^2(\cD; \Z) = H^2(Y; \Z)$.
\end{theorem}

\begin{proof}
First, I will show that $e(\cF)$ is essentially just the Euler class
of the tangent bundle to the branched surface $B = B(\mu)$. If we
relax the condition that $B$ has a well-defined tangent plane along
its singular locus, then as a \2-complex $B$ is isotopic to $\cD^2$.
As $\cT$ has one vertex, the coboundary map
$C^2(\cD;\Z) \to C^3(\cD;\Z)$ is 0, and hence $B \hookrightarrow Y$
induces an isomorphism on $H^i$ for $i \leq 2$. Thus
$e(\cF) \in H^2(Y; \Z)$ is determined by the restriction of $\TF$ to
$B$, which is bundle isomorphic to the tangent bundle $\TB$.  To prove
the theorem, I will compute $e(\TB)$ as the obstruction to the unit
circle bundle $\UTB \subset \TB$ having a section; see
e.g.~\cite[Chapter~4]{CandelConlon2003} for a detailed discussion of
the Euler class along these lines.

Let $\Bpunc$ be $B$ with an open regular neighborhood of $\cT^1$
removed. As $\Bpunc$ deformation retracts to a graph, specifically one
isotopic to $\cD^1$, the restricted bundle $\UTBpunc$ necessarily has
a section, and here is a concrete one.  First, on $\Bpunc \cap \cT^2$
we take the section $\eta$ that is tangent to the faces of $\cT^2$ and
within each face points away from the long edge as shown in
Figure~\ref{fig:sect_on_tri}.  Now extend $\eta$ over all of $\Bpunc$
using the template shown in Figure~\ref{fig:section}, where this
includes using the mirrored picture when needed; the section $\eta$ is
roughly summarized by the rule that it points away from the very long
edge.

Each component of $B \setminus \Bpunc$ is a disc that is a shrunken
copy of a face of $\cD^2$.  For each such disc $D$, fix a section
$\delta$ of $\UTB$ over $D$.  Thus we have two different sections of
$\UTB$ over $\partial D$, namely the restrictions of $\delta$ and
$\eta$.  The usual cocycle for $e(\TB)$ assigns to the face containing
$D$ the integer that expresses the difference between $\delta$ and
$\eta$ on $\partial D$.  To compare these sections, consider a third
section $\nu$ on $\partial D$ that is normal to $\partial D$ and
points directly out of $D$ and hence into $\Bpunc$.  To understand how
$\eta$ compares to $\nu$ look at Figure~\ref{fig:section} which
contains six segments of $\partial \Bpunc$.  On the segment near the
very long edge, the section $\eta$ is equal to $\nu$, and there are
three other segments where $\eta = -\nu$.  The remaining two segments
correspond to the mixed edges of the tetrahedron, and there $\eta$
twists through half a turn with respect to $\nu$.  Now $D$ gets an
orientation from the face of $\cD^2$ that contains it, which we use to
orient $\partial D$ as usual.  Zooming in on the leftmost mixed edge
of Figure~\ref{fig:section}, we get Figure~\ref{fig:twisting} which
shows that $\eta$ does a \emph{clockwise} twist on that segment
with respect to $\nu$.  The same is true for the other mixed edge in
Figure~\ref{fig:section}, as well the two mixed edges in the mirror
image of Figure~\ref{fig:section}.  Thus if $D$ is dual to an edge
$\tau \in \cT^1$, we see that $\eta$ twists through
$\frac{1}{2}\mixed_\mu \tau$ clockwise full turns with respect to
$\nu$ as you go around $\cD$; in particular, the number
$\frac{1}{2}\mixed_\mu \tau$ is in $\Z$.  In contrast, the section
$\delta$ that extends over $D$ does one clockwise twist compared to
$\nu$.  Thus the difference between $\eta$ and $\delta$ on
$\partial D$, expressed in terms of \emph{anticlockwise} twists is:
\[
\phi_\mu(D) = -\frac{1}{2}\mixed_\mu(\tau) + 1
\]
To confirm that $\phi_\mu(D)$ is the obstruction to $\UTB$ having a
section, simply note that if $\phi_\mu(D)$ is the coboundary of some
one chain $\alpha \in C^1(\cD;\Z)$, then $\alpha$ gives instructions
on how to twist $\eta$ along the ``arms'' of $\Bpunc$ (which
correspond to the edges of $\cD^1$) to get a new section that agrees
with our section of $\UTB$ over $B \setminus \Bpunc$ and hence gives
global section for $\UTB$. Finally, to see that my conventions mean
that $[\phi_\mu]$ represents $e(\TB)$ rather than $-e(\TB)$, simply
carry out the corresponding calculation for e.g.~$TS^2$.
\end{proof}

Applying this method to the sample at hand yields:
\begin{theorem}\label{thm:tauteuler0}
  At least 32{,}347 of the manifolds in \QHSpheres\ have taut
  foliations whose Euler class vanishes and hence are orderable.
\end{theorem} 
\begin{proof}
  For Theorem~\ref{thm:usefoliar}, some 160,003 1-vertex
  triangulations were found for manifolds in \QHSpheres\ that each had
  at least one foliar orientation. For each of these 3.5 million
  foliar orientations, I used Theorem~\ref{thm:cocycle} to determine
  whether the associated foliation had $e(\cF) = 0$.  There were
  610{,}326 foliar orientations where $e(\cF) = 0$, showing that some
  32{,}347 manifolds in \QHSpheres\ are orderable.  The triangulations
  and the complete list of these foliar orientations are available at
  \cite{PaperData}, as well as the code for computing the Euler class.
\end{proof}

\begin{remark}\label{rem:tauteuler}
  It is natural to ask to what extent $e(\cF)$ behaves like a random
  element of $H^2(Y ; \Z)$.  Suppose that for each foliar orientation
  one has $e(\cF) = 0$ with probability $1/\abs{H^2(Y ; \Z)}$, and
  that the events $e(\cF) = 0$ are independent for each foliar
  orientation.  Then for the 3.5 million foliar orientations in
  Theorem~\ref{thm:usefoliar}, the expected number with $e(\cF) = 0$
  is 276{,}602.7 and the expected number of $Y$ with at least one
  foliar orientation with $e(\cF) = 0$ is 67{,}595.3.  (In both cases,
  the standard deviation is tiny by comparison: 416.4 and 149.6,
  respectively.) These predictions are off by about a factor of 2 in
  \emph{opposite} directions from what was actually observed:
  610{,}326 and 32{,}347.  I have no explanation for this odd
  phenomenon.
\end{remark}

\begin{remark} \label{rem:gauge}
  It would be very interesting to gauge the extent
  Theorem~\ref{thm:tauteuler0} finds all manifolds in \QHSpheres\ with
  taut foliations whose Euler class vanishes.  One approach here would
  be to refine the method of Section~\ref{sec:compute-HF} so that it
  computes $\HFhat(Y)$ itself together with its
  $\mathrm{Spin}^c$-grading as that provides obstructions to having
  taut foliations with particular Euler classes.  For those $Y$ that
  are Dehn fillings on Floer simple manifolds in $\cC$, the methods of
  \cite{Rasmussens2015} should make this feasible using the Turaev
  torsion computations included in \cite{PaperData}.
\end{remark}

%% file: plots/local_edge_orientations.tex
\begin{tikzpicture}[nmdstd]
  \begin{scope}
    \coordinate (a) at (0, 0); 
    \coordinate (b) at (1, 1); 
    \coordinate (c) at (0, 2); 
    \begin{scope}[line width=1.5pt]
      \draw[mid arrow=0.55] (a) -- (c); 
      \draw[mid arrow=0.62] (a) -- (b);
      \draw[mid arrow=0.55] (b) -- (c);
      \draw (a) -- (b) -- (c);
    \end{scope}
    \node[align=center, left=0.2] at (0, 1) {long \\ edge};
    \node[align=left, right=0.2] at (0.5, 1.9) {top \\ short edge};
    \node[align=left, right=0.2] at (0.5, 0.1) {bottom \\ short edge};
  \end{scope}


  \begin{scope}[scale=0.8, shift={(9, -1.5)}]
    \coordinate (A) at (0, 0); 
    \coordinate (B) at (1.8, 1.5); 
    \coordinate (C) at (-1.45, 3.0); 
    \coordinate (D) at (0.5, 5);

    \begin{scope}[line width=1.5pt]
      \draw[mid arrow=0.55] (A) -- (C);
      \draw[mid arrow=0.5] (C) -- (D);
      \draw[mid arrow=0.6] (A) -- (B);  
      \draw[mid arrow=0.7] (B) -- (D);
      \draw[name path=AD]
            ($(A)+(0.015,0.05)$) -- ($(D)+(-0.016,-0.02)$);
      \draw[mid arrow=0.25, mid arrow=0.75]
            ($(A)+(0.015,0.05)$) -- ($(D)+(-0.016,-0.02)$);
      \path[name path=BC] (B) -- (C);

      \path[name intersections={of=AD and BC, by=X}];
      \path[name path=small circle, red] (X) circle (0.25);
      \path[name intersections={of=AD and small circle, by={U,V}}];
      \draw[{Round Cap[reversed]}-{Round Cap[reversed]}, 
            color=white, line width=1.9pt] (U) -- (V);
      \draw[mid arrow=0.3, mid arrow=0.75]  (B) -- (C);
    \end{scope}

    \coordinate (c) at (0.5, 2.9); 
    \coordinate (b) at (1.7, 2.7); 
    \coordinate (a) at (3.0, 3.5); 
    \draw[->, dashed, name path=indicator, line width=1.2pt] (a) to[curve through={(b)}] (c);
    \node[above] at (a) {very long};
    \path[name path=BD] (B) -- (D);
    \fill[color=white, name intersections={of=BD and indicator, by=x}]
         (x) circle (0.145);
    \draw[line width=1.5pt] (B) -- (D); 

    \node[above left=-3pt, align=center] 
        at ($(C)!0.5!(D)$) {compatibly \\ short};
    \node[below right=2pt, align=center] 
        at ($(A)!0.5!(B)$) {compatibly \\ short};
  \end{scope}

  \begin{scope}[subfig label]
    \node at (0, -1.2) {(a)}; 
    \node at (6, -1.2) {(b)};
  \end{scope}
\end{tikzpicture} 

%% file: plots/branched.tex
\begin{tikzpicture}[scale=0.7, nmdstd]
  \coordinate (A) at (2.06,3.8);
  \coordinate (B) at (8.8,5.1);
  \coordinate (C) at (1.2,9.4);
  \coordinate (D) at (6.9,11.6);

  \coordinate (U) at ($(A)!0.55!(C)$);
  \coordinate (V) at ($(B)!0.40!(C)$);
  \coordinate (W) at ($(B)!0.56!(D)$);
  \coordinate (X) at ($(A)!0.64!(D)$);

  \coordinate (K) at ($(V)!0.43!(W)$);
  \coordinate (L) at ($(U)!0.43!(X)$);
  \coordinate (M) at ($(U)!0.65!(V)$);
  \coordinate (N) at ($(X)!0.49!(W)$);

  \coordinate (P) at ($(C)!0.155!(D)$);
  \coordinate (Q) at ($(A)!0.79!(B)$);

  \def\toptrileft{%
    ($(L)!0.05!(X)$).. controls ($(L)!0.3!(U)$) and ($(P)!0.3!(U)$).. (P)
  }
  \def\toptriright{%
     (P) .. controls ($(P)!0.3!(V)$) and ($(K)!0.25!(V)$) .. ($(K)!0.04!(W)$) 
  }

  \def \bottrileft{%
    (Q) .. controls ($(K)!0.6!(Q)$) and +(1.5,0) .. ($(U)!0.80!(M)$)
  }  
  \def \bottriright{%
    ($(X)!0.85!(N)$)  .. controls +(1.5,0) and ($(W)!0.7!(Q)+(0.5,0)$) .. (Q)
  }

  \fill[color=nmdlight!70] \bottriright -- \bottrileft -- cycle;
  \fill[color=nmdlight!20] (X) -- (W) -- (V) -- (U) -- cycle;
  \fill[color=nmdlight!45] 
        (L) -- \toptrileft -- \toptriright -- (K) -- cycle;

  \begin{scope}[color=nmddark, line width=1.2pt]
    \draw  (X) -- (W) -- (V) -- (U) -- (L);
    \draw[dash pattern=on 3pt off 3pt on 3pt off 7pt
          on 3.5pt off 4 pt on 3.5 pt off 6pt on 10cm]
          (L) -- (X);
    \draw[dash pattern=on 11pt off 7pt on 12pt off 7pt on 10cm]
          (M) -- (N);
    
    \draw \toptriright;
    \draw[dash pattern=on 29.5pt off 8pt] \toptrileft;
    \draw (L) -- ($(L)!0.23!(K)$);
    \draw ($(L)!0.4!(K)$) -- (K);

    \draw \bottrileft;
    \draw[dash pattern=on 8pt off 3pt on 3pt off 3pt on 3pt
          off 6pt on 44pt off 7pt on 10cm] \bottriright;        
    \draw (W) -- (V);
  \end{scope}

  \begin{scope}[line width=1pt]
    \draw[mid arrow=0.4] (A) -- (B);
    \draw[mid arrow=0.4] (A) -- (C);
    \draw[mid arrow=0.7]  (A) -- ($(A)!0.37!(D)$);
    \draw[mid arrow=0.5]  (X) -- (D);
    \draw[mid arrow=0.27] (B) -- (C);
    \draw[mid arrow=0.77] (B) -- (D);
    \draw[mid arrow=0.5] (C) -- (D);
  \end{scope}


  \begin{scope}[shift={(12.5,5.5)}]
    \coordinate (a) at (0, 0); 
    \coordinate (b) at (4.5, 2.5); 
    \coordinate (c) at (0, 5); 
    \coordinate (ac) at ($(a)!0.5!(c)$);
    \coordinate (ab) at ($(a)!0.57!(b)$);
    \coordinate (bc) at ($(c)!0.57!(b)$);
    \coordinate (branch) at (1.25, 2.5);
    \begin{scope}[line width=2pt, color=nmddark]
      \draw (ac) -- (branch);
      \draw (branch) .. controls +(0.55, 0) and +(-0.5, 1).. (ab);
      \draw (branch) .. controls +(0.55, 0) and +(-0.5, -1).. (bc);
    \end{scope}

    \begin{scope}[line width=1pt]
      \draw[mid arrow=0.8, mid arrow=0.27] (a) -- (c); 
      \draw[mid arrow=0.4] (a) -- (b);
      \draw[mid arrow=0.65] (b) -- (c);
      \draw (a) -- (b) -- (c);
    \end{scope}
  \end{scope}

\end{tikzpicture}

%% file: plots/mixed.tex
\begin{tikzpicture}[nmdstd]

  \begin{scope}[scale=0.8, shift={(9, -1.5)}]
    \coordinate (A) at (0, 0); 
    \coordinate (B) at (1.8, 1.5); 
    \coordinate (C) at (-1.45, 3.0); 
    \coordinate (D) at (0.5, 5);

    \begin{scope}[line width=1.5pt]
      \draw[mid arrow=0.55] (A) -- (C);
      \draw[mid arrow=0.5] (C) -- (D);
      \draw[mid arrow=0.6] (A) -- (B);  
      \draw[mid arrow=0.7] (B) -- (D);
      \draw[name path=AD]
            ($(A)+(0.015,0.05)$) -- ($(D)+(-0.016,-0.02)$);
      \draw[mid arrow=0.25, mid arrow=0.75]
            ($(A)+(0.015,0.05)$) -- ($(D)+(-0.016,-0.02)$);
      \path[name path=BC] (B) -- (C);

      \path[name intersections={of=AD and BC, by=X}];
      \path[name path=small circle, red] (X) circle (0.25);
      \path[name intersections={of=AD and small circle, by={U,V}}];
      \draw[{Round Cap[reversed]}-{Round Cap[reversed]}, 
            color=white, line width=1.9pt] (U) -- (V);
      \draw[mid arrow=0.3, mid arrow=0.75]  (B) -- (C);
    \end{scope}

    \node[above left=-3pt, align=center] 
        at ($(C)!0.5!(D)$) {compatibly \\ short};
    \node[below right=2pt, align=center] 
        at ($(A)!0.5!(B)$) {compatibly \\ short};

    \node[left=4pt, align=center] 
        at ($(A)!0.5!(C)$) {mixed};
    \node[right=3pt, align=center] 
        at ($(B)!0.55!(D)$) {mixed};
  \end{scope}
\end{tikzpicture}

%% file: plots/tri_section.tex
\begin{tikzpicture}[nmdstd]

  \coordinate (A) at (0.7, 0.0);
  \coordinate (B) at (4.5, 2.5);
  \coordinate (C) at (0.7, 5.0);

  \coordinate (X) at (barycentric cs:A=0.5,B=0.07,C=0.5);
  \coordinate (Y) at (barycentric cs:A=0.5,B=0.5,C=0.5);
  \coordinate (U) at (barycentric cs:A=0.05,B=0.6,C=0.35);
  \coordinate (V) at (barycentric cs:A=0.35,B=0.6,C=0.05);

  \begin{scope}[line width=0.9pt]
    \draw[mid arrow=0.3, mid arrow=0.7] (A) -- (C);
    \draw[mid arrow=0.5] (A) -- (B);
    \draw[mid arrow=0.55] (B) -- (C);
  \end{scope}

  \def\topbranch{(Y) .. controls +(0:0.7) and +(-105:0.3) .. (U)}
  \def\bottombranch{(Y) .. controls +(0:0.7) and +( 105:0.3) .. (V)}
  
  \begin{scope}[line width=1.2, color=nmddark!50]
    \draw (X) -- (Y);
    \draw \topbranch; 
    \draw \bottombranch;
  \end{scope}

  \begin{scope}[line width=1.9pt, -{Latex[length=7pt]}]
    \path[tangent vector={0.3}{0.7}] (X) -- (Y);
    \path[tangent vector={0.45}{0.75}] \topbranch;
    \path[tangent vector={0.45}{0.75}] \bottombranch;
  \end{scope}

  \node[below=3pt] at ($(Y)+(0.2,0)$) {$B^\circ$};

\end{tikzpicture}

%% file: plots/section.tex
\begin{tikzpicture}[scale=1.7, nmdstd]
  \coordinate (A) at (1.06,0.8);
  \coordinate (B) at (7.8,2.1);
  \coordinate (C) at (0.2,6.4);
  \coordinate (D) at (5.9,8.6);

  \coordinate (U) at ($(A)!0.55!(C)$);
  \coordinate (V) at ($(B)!0.40!(C)$);
  \coordinate (W) at ($(B)!0.56!(D)$);
  \coordinate (X) at ($(A)!0.64!(D)$);

  \coordinate (K) at ($(V)!0.43!(W)$);
  \coordinate (L) at ($(U)!0.43!(X)$);
  \coordinate (M) at ($(U)!0.64!(V)$);
  \coordinate (N) at ($(X)!0.38!(W)$);

  \coordinate (P) at ($(C)!0.20!(D)$);
  \coordinate (Q) at ($(A)!0.69!(B)$);

  \coordinate (U1) at ($(U)!0.32!(V)$);
  \coordinate (U2) at ($(U)!0.22!(X)$);
  \def\Ucorner{%
    (U1) .. controls +(80:0.25) and +(-5:0.2) ..  (U2)
  }

  \coordinate (V1) at ($(V)!0.20!(W)$);
  \coordinate (V2) at ($(V)!0.175!(U)$);
  \def\Vcorner{%
    (V1) .. controls +(170:0.35) and +(60:0.3) ..  (V2)
  }

  \coordinate (W1) at ($(W)!0.35!(X)$);
  \coordinate (W2) at ($(W)!0.2!(V)$);
  \def\Wcorner{%
    (W1) .. controls +(-80:0.2) and +(170:0.2) ..  (W2)
  }

  \coordinate (X1) at ($(X)!0.10!(U)$);
  \coordinate (X2) at ($(X)!0.20!(W)$);
  \def\Xcorner{%
    (X1) .. controls +(-5:0.30) and +(-145:0.34) ..  (X2)
  }

  \def\toptrileft{%
    ($(L)!0.05!(X)$).. controls ($(L)!0.3!(U)$) and ($(P)!0.3!(U)$).. (P)
  }
  \def\toptriright{%
     (P) .. controls ($(P)!0.65!(V)$) and ($(K)!0.1!(V)$) .. ($(K)!0.04!(W)$) 
  }

  \path[add coordinate={0.75}{P1}] \toptrileft;
  \path[add coordinate={0.15}{P2}] \toptriright;
  \def\Pcorner{%
     (P1) .. controls +(-10:0.18) and +(-127:0.19) ..  (P2)
  }
  
  \def\toptrileft{%
    ($(L)!0.028!(X)$) .. controls +(-152:0.52) and +(-88:0.7) .. (P1)
  }
  \def\toptriright{%
      (P2) .. controls +(-40:1.0) and +(180:1.05) .. (K)
  }

  \def \bottrileft{%
    (Q) .. controls ($(K)!0.6!(Q)$) and +(1.5,0) .. ($(U)!0.80!(M)$)
  }  
  \def \bottriright{%
    ($(X)!0.85!(N)$)  .. controls +(1.2,0) and ($(W)!0.7!(Q)+(0.5,0)$) .. (Q)
  }

  \path[add coordinate={0.19}{Q1}] \bottrileft;
  \path[add coordinate={0.78}{Q2}] \bottriright;
  \def\Qcorner{%
     (Q2) .. controls +(170:0.2) and +(53:0.17) .. (Q1)
  }
  \def \bottrileft{%
    (Q1) .. controls +(125:0.8) and +(-5:1.20) .. ($(U)!0.965!(M)$)
  }
  \def \bottriright{%
    ($(X)!0.93!(N)$)  .. controls +(-1:1.1) and +(77:1.0) .. (Q2)
  }

  \fill[color=nmdlight!55] \bottriright -- \Qcorner -- \bottrileft-- cycle;
  \fill[color=nmdlight!15] \Xcorner -- \Wcorner -- 
        \Vcorner -- \Ucorner -- cycle;
  \begin{scope}[color=nmdlight!40] 
    \fill (L) -- \toptrileft -- \Pcorner -- \toptriright -- (K);
    \fill (L) rectangle +(0.1, 0.1);
  \end{scope}
  \begin{scope}[color=nmddark!60, line width=0.7pt]
    \draw (U1) -- (V2);
    \draw (V1) -- (W2);
    \draw (W1) -- (X2);
    \draw (X1) -- ($(X1)!0.14!(U2)$);
    \draw (L) -- (U2);

    \draw[dash pattern=on 1 off 4.2 on 64 off 8 on 100,
          dash phase=2] \toptrileft;
    \draw \toptriright;
    \draw (L) -- ($(L)!0.27!(K)$);
    \draw ($(L)!0.36!(K)$) -- (K);

    \draw \bottrileft;
    \draw[dash pattern=on 11 off 5 on 4.5 off 5 on 4.5 off 5
          on 4.5 off 5 on 4.5 off 5 on 4.2 off 6 
          on 88 off 8 on 100
    ] \bottriright;   
    \draw \Qcorner;
    \draw (W2) -- (V1);

    \draw[dash pattern=on 32 off 8 on 16 off 22 on 100] (M) -- (N);
  \end{scope}

  \begin{scope}[color=nmddark, line width=1.1pt]
    \draw \Pcorner;
    \draw \Qcorner;
    \draw \Xcorner;
    \draw \Ucorner;
    \draw[dash pattern=on 37 off 8] \Vcorner;
    \draw \Wcorner;
  \end{scope}
  
  \begin{scope}[line width=1pt]
    \draw[mid arrow=0.4] (A) -- (B);
    \draw[mid arrow=0.4] (A) -- (C);
    \draw[mid arrow=0.7]  (A) -- ($(A)!0.38!(D)$);
    \draw[mid arrow=0.5]  ($(A)!0.625!(D)$) -- (D);
    \draw[mid arrow=0.15, mid arrow=0.9] (B) -- (C);
    \draw[mid arrow=0.77] (B) -- (D);
    \draw[mid arrow=0.5] (C) -- (D);
  \end{scope}

  

  \begin{scope}[line width=1.6pt, 
    -{Latex[length=7pt]}]
    %

    %
    \draw (U1) -- ($(U1)!-0.5cm!(U)$);
    \draw (U2) -- ($(U2)!0.4cm!(U)$);
    \path[vector along={0.6}{(5:-0.7cm)}] \Ucorner;

    \draw (V1) -- ($(V1)!0.4cm!(V)$);
    \draw (V2) -- ($(V2)!0.4cm!(V)$);

    \draw (X1) -- ($(X1)!-0.3cm!(X)$);
    \draw (X2) -- ($(X2)!-0.35cm!(X)$);
    \path[normal vector={0.5}{-0.5cm}] \Xcorner;
    
    \draw (W1) -- ($(W1)!0.4cm!(W)$);
    \draw (W2) -- ($(W2)!-0.4cm!(W)$);
    \path[vector along={0.3}{(-10:0.6cm)}] \Wcorner;

    \draw (5.5, 5.5) to [out=20, in=165, relative, looseness=1.0] +(-70:0.4cm);
    \draw (5.0, 5.3) to [bend left=5pt] +(-80:0.4cm);
    \path[tangent vector={0.47}{0.8cm}] (W2) -- (V1);
 
    \draw (4.105, 4.61) to  +(-55:0.4cm);
    \coordinate (f1) at (2.5, 4.5);
    \coordinate (g1) at ($(f1)+(-60:0.4cm)$);
    \draw \easybezier{(f1)}{(g1)}{-40}{0.4}{15}{0.2};
    \coordinate (f2) at (4.8, 4.38);
    \coordinate (g2) at ($(f2)+(-105:0.4cm)$);
    \draw \easybezier{(f2)}{(g2)}{15}{0.4}{-7.5}{0.2};
    \draw 
         (3.3, 4.2) .. controls ++(-45:0.2) and ++(-10:-0.15) .. ++(-20:0.5cm);
    \path[tangent vector={0.45}{0.7cm}] (U1) -- (V2);
    \path[tangent vector={0.25}{0.5cm}] (L) -- (U2);

    %
    %
    \draw (P1) -- ($(P1)!0.4cm!5:(P)$);
    \draw (P2) -- ($(P2)!0.4cm!(P)$);

    \path[tangent vector={0.35}{0.8cm}] \toptrileft;
    \path[tangent vector={0.35}{-0.8cm}, 
          tangent vector={0.7}{-0.8cm}]
          \toptriright;

    \draw (2.0, 5.7) to [bend left=10pt, looseness=0.9] +(120:0.4);
    \draw (3.3, 5.0) to [bend left=10pt, looseness=0.9] +(150:0.5);
    \draw (2.3, 4.9) to [bend left=20pt] +(140:0.5);

    %
    %
    \draw (Q1) -- ($(Q1)!0.4cm!12:(Q)$);
    \draw (Q2) -- ($(Q2)!0.4cm!10:(Q)$);

    \path[tangent vector={0.45}{-0.8cm}] \bottrileft;
    \path[tangent vector={0.70}{0.7cm},
          tangent vector={0.45}{0.7cm}]
          \bottriright;

    \draw (5.6, 4) to [out=15, in=170, relative, looseness=1.0] +(-60:0.5);
    \draw (5.3, 3.2) to [out=15, in=170, relative, looseness=0.8] +(-50:0.4);
  \end{scope}


\end{tikzpicture}

%% file: plots/orientation.tex
\begin{tikzpicture}[nmdstd]

  \begin{scope}[rotate=10]
  \coordinate (T0) at (3.66, 0.5);
  \coordinate (T1) at (3.66, 5.87);
  \coordinate (A) at ($(T0)!0.525!(T1)$);
  \coordinate (B) at (6.7,2.4);
  \coordinate (C) at (6.45,5.0);

  \def\bdryD{(A) circle [x radius=1.66, y radius=0.83]} 
  \path[name path=bdryD] \bdryD;
 
  \path[name path=big D]
       (A) circle [x radius=2.5, y radius=1.25, rotate=3];
  \path[name path=R1] (A) -- +(-173:3);
  \path[name path=R2] (A) -- +(-164:3);
  \path[name path=R3] (A) -- +(-153:3);
  \path[name intersections={of=bdryD and R1, by=U1}];
  \path[name intersections={of=bdryD and R2, by=U2}];
  \path[name intersections={of=bdryD and R3, by=U3}];
  \path[name intersections={of=big D and R1, by=V1}];
  \path[name intersections={of=big D and R2, by=V2}];
  \path[name intersections={of=big D and R3, by=V3}];

  \begin{scope}[line width=1.4pt, -{Latex[length=7pt]}]
    \draw (U1) -- (V1);
    \draw (U2) -- (V2);
    \draw (U3) -- (V3);
  \end{scope}

  \fill[green!10] \bdryD;
  \path[mid arrow=0.42, 
        mid arrow=0.69, 
        mid arrow=0.89, 
        line width=1.1,
        ] \bdryD;

  \draw[line width=1.1] ($(A)+(93:1.66 and 0.83)$) arc 
         [x radius=1.66, y radius=0.83, start angle=93, 
          end angle=447];

  \path[name path=AB] (A) -- (B);
  \path[name path=AC] (A) -- (C);
  \path[name intersections={of=bdryD and AB, by=X}];
  \path[name intersections={of=bdryD and AC, by=Y}];

  \def\commoncurve{($(A)+(-31:1.66 and 0.83)$)
       arc [x radius=1.66, y radius=0.83, start angle=-31,
            end angle=50.5]}

  \fill[color=nmdlight!30] 
      (B) -- (X) -- \commoncurve -- (Y) -- (C) 
      .. controls +(-60:0.7) and +(70:1) .. (B);
  
  \draw[color=nmddark!60, line width=1.0pt]
      (B) -- (X) -- \commoncurve -- (Y) -- (C);
  \draw[color=nmddark, line width=1.3pt] \commoncurve;

  
  \begin{scope}[line width=1.4pt, -{Latex[length=7pt]}]
    \path[tangent vector={0.2}{0.8}] (C) -- (Y);
    \path[tangent vector={1.0}{0.8}] (C) -- (Y);
    \path[tangent vector={0.0}{0.9}] (X) -- (B);
    \path[tangent vector={0.6}{-0.7}] \commoncurve;

    \draw (6.05, 4.35) 
      to [bend right=15pt,looseness=0.8]  +(-110:0.7) 
      node[below=-2pt] {$\eta$};
    
    \draw (6.0, 3.2) 
      to [bend right=15pt,looseness=0.8]  +(-10:0.7);
  \end{scope}

  \begin{scope}[line width=1.5pt]
    \draw[mid arrow=0.7] (A) -- (T1);
    \draw[mid arrow=0.7] (T0) -- ($(T0)!0.35!(T1)$);
  \end{scope}


  \node[below left=3pt] at (A) {$D$};
  \node[right] at (T1) {$\tau$};
  \node[left=3pt] at ($(A)!0.6!(T1)$) {$\mu$};

  \node[left] at (V2) {$\nu$};
  \node at (6.4, 3.7) {$B^\circ$};
  \node at (2.15, 4.05) {$\partial D$};

  \end{scope}
\end{tikzpicture}

%% file: representations.tex
\section{Representations to $\PSLRtilde$}\label{sec:numreps}

A very useful technique for showing a 3-manifold $Y$ is orderable is
to find a representation from $\pi_1(Y)$ to $\PSLRtilde$, which is the
universal covering Lie group of $\PSL{2}{\R}$, and also its
universal central extension (see e.g.~\cite[\S 5]{Ghys2001} or \cite[Chapter
2]{Bruggeman1994}):
\[
0 \to \Z \to \PSLRtilde \to \PSL{2}{\R} \to 1
\]
The action of $\PSL{2}{\R}$ on the circle $P^1(\R)$ lifts to an action
of $\PSLRtilde$ on $\R$.  If $Y$ is an irreducible
compact \3-manifold, then by Theorem~1.1 of
\cite{BoyerRolfsenWiest2005} having a \emph{nontrivial} representation
$\rho \maps \pi_1(Y) \to \PSLRtilde$ implies that $Y$ is orderable.
The last piece of Theorem~\ref{thm:main} is:
\begin{theorem}\label{thm:tilreps}
  There are at least $64{,}180$ manifolds $Y \in \cY$ where $\pi_1 Y$
  has a nontrivial homomorphism to $\PSLRtilde$; all these $Y$ are
  therefore orderable.  
\end{theorem}
The manifolds in this theorem are all non-L-spaces, consistent with
Conjecture~\ref{BGWconjecture}.

\subsection{Finding representations numerically}

My starting point for Theorem~\ref{thm:tilreps} was numerical evidence
of representations to $\PSLRtilde$ coming from points on the extended
Ptolemy variety that were computed using numerical algebraic geometry.
I now describe this method in detail. Throughout, let $Y$ be a closed
\3-manifold given as Dehn filling on a manifold $M$ with $\partial M$
a torus, and suppose $\cT$ is a fixed ideal triangulation of $M$.  The
first step is to find many (preferably all) representations
$\pi_1(Y) \to \PSL{2}{\C}$; later, these will be filtred to identify
those that are conjugate into $\PSL{2}{\R}$ and lift to
$\PSLRtilde$. Put another way, we need to compute the
$\PSL{2}{\C}$-character variety $X(Y)$, which is essentially
representations $\pi_1(Y) \to \PSL{2}{\C}$ modulo conjugacy and is an
affine complex algebraic variety.  For a relatively simple
$\Q$-homology sphere like those in $\cY$, one expects that
$\dim_\C X(Y) = 0$ and so there are only finitely many representations
up to conjugacy.

Since $\pi_1(Y)$ is a quotient of $\pi_1(M)$, a representation of
$\pi_1(Y)$ gives one of $\pi_1(M)$.  A $\PSL{2}{\C}$-representation of
$\pi_1(M)$ can frequently be described in terms of the ideal
triangulation $\cT$. Specifically, if $\cTtil$ is the induced ideal
triangulation of the universal cover of $M$, one tries to encode
$\rho \maps \pi_1(M) \to \PSL{2}{\C}$ by an equivariant developing map
$\cTtil \to \H^3$ that takes each topological ideal tetrahedron of
$\cTtil$ to a geodesic ideal tetrahedron in $\H^3$.  (Technical aside:
For some pairs $(\rho, \cT)$ no such developing map exists,
and in practice one does lose some representations by taking this
perspective; for ways around this limitation, which I did not use
here, see \cite{Segerman2012, GoernerZickert2018}.)  There are two
ways of encoding such developing maps, both of which form algebraic
varieties closely related to the character variety of $M$.

First, Thurston used one complex parameter for each tetrahedron in
$\cT$ to describe its geometric shape in $\H^3$, giving his
\emph{gluing equation variety} $D(\cT)$ which is cut out by one
equation for each edge of $\cT$; see e.g.~\cite[\S
3]{DunfieldGaroufalidis2012} for an overview and more on the map
$D(\cT) \to X(M)$.  If $\alpha$ is the slope where $Y = M(\alpha)$,
then by adding the cusp equation corresponding to $\alpha$ one obtains
a subvariety $D(\cT, \alpha)$ of $D(\cT)$.  The subvariety
$D(\cT, \alpha)$ consists of all points mapping to $X(Y)$ together
with those corresponding to representations of $\pi_1(M)$ where
$\alpha$ acts by a nontrivial parabolic.  More formally, if
$X(M, \alpha)$ is the subset of $X(M)$ consisting of characters
$[\rho]$ with $\tr^2 \big(\rho(\alpha)\big) = 4$, then
$D(\cT, \alpha)$ is the preimage of $X(M, \alpha)$.

Second, the \emph{enhanced Ptolemy variety} $P(\cT)$ for $\PSL{2}{\C}$
uses one complex parameter for each edge of $\cT$ as well as two such
parameters for recording the eigenvalues of a basis of
$\pi_1(\partial M)$; there is then one polynomial equation for each
tetrahedron of $\cT$.  The basic Ptolemy variety was introduced in
\cite{GaroufalidisThurstonZickert2015} and extended to non-boundary
parabolic representations in \cite{Zickert2016} in the enhanced
version used here, but I recommend the reader start with
\cite{GaroufalidisGoernerZickert2015a} which focuses on the case of
$\PSL{2}{\C}$ rather than the more general $\SL{n}{\C}$ of
\cite{GaroufalidisThurstonZickert2015}. To study the points of
$P(\cT)$ mapping to $X(M, \alpha)$, one imposes an additional equation
involving the two cusp parameters.

\subsection{The Ptolemy advantage}

The gluing and Ptolemy coordinates are in a certain sense dual
\cite[\S 12]{GaroufalidisGoernerZickert2015b}.  In particular, as
$\cT$ has the same number of edges as tetrahedra, both $D(\cT)$ and
$P(\cT)$ live in $\C^n$ for roughly the same $n$.  However, from the
perspective of computational algebraic geometry there is a key
difference between them: the equations defining $P(\cT)$ are typically
much simpler than those defining $D(\cT)$.  Specifically, they have
much lower degrees. Both symbolic (e.g.~G\"obner bases) and numerical
algorithms in algebraic geometry tend to be highly sensitive to the
degree of the equations defining the variety.  With the numerical
method I used here, working with $P(\cT)$ was typically 10 times
faster than using $D(\cT)$ and sometimes more than 50 times faster.
Moreover, this was using a description of $D(\cT)$ with two variables
per tetrahedron (representing $z_i$ and $1 - z_i$ in the usual
parlance) which was itself 10 to 100 times faster than the standard
description of $D(\cT)$ with just one variable per tetrahedron.

\subsection{Finding points on the Ptolemy variety}

In its simplest form, the Ptolemy variety parameterizes
representations to $\SL{2}{\C}$ rather than $\PSL{2}{\C}$, though the
latter case can be handled as well \cite[\S
3]{GaroufalidisGoernerZickert2015a}.  Since a representation to
$\PSL{2}{\R}$ that lifts to $\PSLRtilde$ necessarily lifts to the
intermediate group $\SL{2}{\R}$, I worked exclusively with the
$\SL{2}{\C}$ Ptolemy varieties. I used the \emph{enhanced Ptolemy
  variety} of \cite{Zickert2016} to study representations of
$\pi_1(M)$ that are not boundary-parabolic, specifically a variant of
that construction due to Goerner (personal communication) that is
easier to implement.  I will use $P(\cT, \alpha)$ to denote the subset
of the extended Ptolemy variety $P(\cT)$ where $\alpha$ acts by a
matrix with eigenvalue 1.  Using the peripheral holonomy map of
\cite[\S 4.2.2]{Zickert2016}, one sees this adds a single equation,
which if we take $\alpha$ to be our preferred meridian for
$\partial M$, is simply the equation $m_1 = 1$ in the notation of
\cite{Zickert2016}.

While Goerner has had tremendous success computing the trace field of
cusped hyperbolic 3-manifolds by applying exact symbolic methods to
the Ptolemy variety \cite{GoernerPtolemyData}, the equations in the
closed case have higher degrees and for the manifolds in $\cY$ this
seems to put the problem beyond the range where one can use Gr\"obner
bases to find the solutions.  Instead, I found many points on
$P(\cT, \alpha)$ via the path homotopy continuation method from
numerical algebraic geometry (see \cite[\S
8.0-8.3]{SommeseVerscheldeWampler2005} for general background),
specifically the software PHCPack \cite{Verschelde1999, PHCPack}.  In
the form used, this method provides compelling numerical evidence of
points on $P(\cT, \alpha)$ though it does not prove that any of the
apparent solutions are close to actual solutions, and it can
definitely miss some points with the parameters that I used.  As
mentioned above, one expects $X(Y)$ to be 0-dimensional, and similarly
for $X(M, \alpha)$ and hence $P(\cT, \alpha)$.  Thus I used PHCPack's
algorithm for finding the isolated points of the input variety since
that should typically be all of them here. The theory behind this is
discussed in \cite[\S 8.1]{SommeseVerscheldeWampler2005}, and I should
point out that it requires having the same number of polynomial
equations defining $P(\cT, \alpha)$ as there are variables.  As
mentioned, for $P(\cT, \alpha)$ we have one variable for each edge and
two for the cusp, as well as one equation for each tetrahedron and a
final equation for the cusp.  This is actually one more variable than
equation, but in fact one really works with the reduced Ptolemy
variety \cite[\S 4.1]{GaroufalidisGoernerZickert2015a}, which is
equivalent to setting one edge parameter to 1, giving us a system of
the required type.  All of the Ptolemy coordinates must be nonzero,
which puts $P(\cT, \alpha)$ into the ``sparse system'' setting of
Table 8.1 of \cite{SommeseVerscheldeWampler2005}, where the number of
points in the random starting system is determined by a mixed volume
computation via the Bernshte\u{\i}n theory.

\subsection{Numerical results}\label{sec:numerical}

For each approximate point of $P(\cT, \alpha)$, the next step was to
identify whether it gives an irreducible representation of $\pi_1(Y)$
whose image is moreover conjugate into $\SL{2}{\R} \leq
\SL{2}{\C}$. For this, I converted over to the tetrahedra shape
coordinates of $D(\cT)$ via
\cite[Equation~(4-6)]{GaroufalidisGoernerZickert2015a} and built the
corresponding holonomy representation
$\rho \maps \pi_1(M) \to \SL{2}{\C}$ as in
\cite[Lemma~3.5]{DunfieldGaroufalidis2012}. To see if $\rho$ factors
through to $\pi_1(Y)$, I checked whether $\rho(\alpha) \approx I$.  For
the manifolds in $\cY$, this procedure identified some 27.8~million
irreducible representations of their fundamental groups to
$\SL{2}{\C}$, an average of 90.5 per manifold. For each
representation, I checked whether it is conjugate into $\SU_2$
or $\SL{2}{\R}$; see Table~\ref{table:reps} for the resulting data.
The running time for each manifold was typically between 10 seconds
and a minute.

\begin{table}
  \begin{center}
    \begin{tabular}{lrrrr}
      \toprule
           &      $\SL2{\C}$ & $\SU_2$ & $\SL{2}\R$ &  $\PSLRtilde$ \\
      \midrule
      L-spaces      &  103.2 &   22.8 &         9.2 &   0.0 \\
      non-L-spaces  &   79.3 &   21.4 &         6.6 &   0.7 \\
      \bottomrule
    \end{tabular}
  \end{center}
  \caption{This table shows the average number of irreducible
    representations of $\pi_1(Y)$ for $Y \in \cY$ found in
    Section~\ref{sec:numerical} for various target groups, separated
    out by whether the manifold is an L-space or not.
    Representations to $\SU_2$ and $\SL{2}{\R}$ are also counted in
    the $\SL{2}{\C}$ column, so e.g.~the average number of Zariski
    dense representations to $\SL{2}{\C}$ for an L-space was 71.2.
  }
  \label{table:reps}
\end{table}

\subsection{The Euler class and representations to $\PSLRtilde$}
\label{sec:eulerreal}

For each representation $\rho$ from $\pi_1(Y)$ to $\SL{2}{\R}$, I
computed the Euler class $e(\rho) \in H^2(Y; \Z)$ of the induced
action on the circle, which is also the obstruction to lifting $\rho$
to $\rhotil \maps \pi_1(Y) \to \PSLRtilde$.  The Euler class vanished
for 113,721 (4.7\%) of the $\SL{2}{\R}$\hyp representations, providing
compelling numerical evidence that some 68,054 manifolds in $\cY$ are
orderable.  All of the apparently orderable manifolds are
non-L-spaces, and I will prove they are indeed orderable in the vast
majority (94.3\%) of cases in Section~\ref{sec:prove_order}.

\begin{remark}
  \label{rem:crazy}
  As with Remark~\ref{rem:tauteuler}, it is natural to ask to what
  extent $e(\rho)$ behaves like a random element of $H^2(Y; \Z)$.  If
  it was random, the expected number of $\rho$ with $e(\rho) = 0$ is
  64{,}530, whereas the observed number is 76.2\% more than that. That
  is, the Euler class was more likely to vanish than the size of
  $H^2(Y; \Z)$ would suggest.  However, if we look at just the
  L-spaces, the opposite is very much the case: if $e(\rho) = 0$ with
  probability $1/\abs{H^2(Y; \Z)}$, the expected number of $\rho$ with
  $e(\rho) = 0$ and where $Y$ is an L-space is 6{,}318 which is much
  larger than the zero such $\rho$ observed.  Indeed, the probability
  that a random $e(\rho)$ is nonzero whenever $Y$ is an L-space in
  this sample is less than $10^{-2{,}700}$.  This is quite good
  evidence for the conjecture that L-spaces are not orderable!

  Given that more than half the manifolds in $\cY$ are L-spaces, you
  might wonder why one expects only $6{,}318$ cases where
  $e(\rho) = 0$ for L-spaces and $58{,}212$ such cases for
  non-L-spaces. While Table~\ref{table:reps} shows that L-spaces
  average 39.4\% more representations to $\SL{2}{\R}$ than
  non-L-spaces, their homology was on average much larger, with a mean
  size of 257.5 versus 61.2.
\end{remark}

\section{Proving orderability}\label{sec:prove_order}

This section gives the proof of Theorem~\ref{thm:tilreps}. As in
Section~\ref{sec:numreps}, let $Y$ be a closed \3-manifold given as
the Dehn filling $M(\alpha)$ with $\cT$ a fixed ideal triangulation of
$M$.  The first step will be to adapt the interval analysis methods of
Section~\ref{sec:intervals} to certify a solution to the gluing
equations of $\cT$ where the associated holonomy representation has
image in $\PSL{2}{\R}$.

\subsection{Shapes giving $\PSL{2}{\R}$-representations}

Consider a solution to Thurston's gluing equations, that is, a point
$z \in D(\cT, \alpha) \subset \C^n$, with associated representation
$\rho \maps \pi_1(M) \to \PSL{2}{\C}$.  When can we conjugate $\rho$
so that the image is in $\PSL{2}{\R}$?  If all the shapes are real and
so $z \in \R^n$, we can choose the associated developing map
$\cTtil \to \H^3$ to have image in our preferred copy of
$\H^2 \subset \H^3$ and hence $\rho$ will preserve $\H^2$ as well.  In
this case, the image of $\rho$ is contained in the stabilizer of
$\H^2$ in $\PSL{2}{\C}$, which is the full isometry group of $\H^2$
and contains $\PSL{2}{\R}$ as its identity component.  For any 
generating set $\{\gamma_i\}$ of $\pi_1(M)$, when all
$\tr\big(\rho(\gamma_i)\big)$ are real (as opposed to some being
purely imaginary), it means $\rho$ has image in $\PSL{2}{\R}$.
However, it is possible to have some non-real shapes and yet $\rho$ is
still conjugate into $\PSL{2}{\R}$.  Indeed, for a
$\PSL{2}{\R}$-representation where some element of $\pi_1(\partial M)$
acts by an elliptic element, it is impossible for all the shapes to be
real, and this was the case for 69.0\% of the $\SL{2}{\R}$\hyp
representations in Table~\ref{table:reps} that lifted to
$\PSLRtilde$\hyp representations. That said, if
$\rho\left(\pi_1(\partial M)\right)$ contains a hyperbolic or
parabolic element then $\rho$ is conjugate into the stabilizer of
$\H^2$ if and only if $z \in \R^n$.

\subsection{Certifying real representations}
\label{sec:certreal}

To use the approach of \cite{HIKMOT2016} to certify $\PSL{2}{\R}$\hyp
representations, it is easiest to work only with real shapes.  As only
31\% of the candidate $\PSLRtilde$-representations identified in
Section~\ref{sec:eulerreal} had this property, I looked at other Dehn
filling descriptions of each $Y$ found by drilling out various short
curves.  Any irreducible $\rho \maps \pi_1(Y) \to \PSL{2}{\R}$ always
sends some element to a hyperbolic one, so drilling out such a loop is
likely to yield a triangulation where $\rho$ will be exhibited by real
shapes.  I succeed in finding an alternate description of $Y$ with a
representation to $\PSLRtilde$ coming from real shapes in 94.3\% of
cases, though this required working with triangulations with as many
as 13 tetrahedra.

Let $\IR$ denote the set of real intervals with rational endpoints.
For a point in $\R^n$ that appears to approximate a point in
$D(\cT, \alpha)$, the first step is to modify
Section~\ref{sec:intervals} to give a criterion for showing a nearby
$\z \in \IR^n$ is guaranteed to contain a point of
$z \in D(\cT, \alpha)$.  To do this, we simply replace $\IC$ with
$\IR$ and use the rectangular form of the gluing equations rather than
the logarithmic one.  The subtle point is that we need to check that
the associated $\rho$ has $\rho(\alpha) = I$ as opposed to
$\rho(\alpha)$ being parabolic.  In Section~\ref{sec:intervals}, this
follows geometrically because the shapes have positive imaginary
component and the log-holonomy of $\alpha$ is $2 \pi i$.  In the
current case, I instead computed the approximate
$\PSL{2}{\R}$-representation $\brho \maps \pi_1(M) \to \PGL{2}{\IR}$
analogous to what is built in Section~\ref{sec:approxholo} and checked
that $\tr^2 \brho(\beta) > 4$ for some slope $\beta \neq \alpha$.
This last condition forces the underlying $\rho(\beta)$ to be
hyperbolic, and so $\rho(\alpha)$ must be $I$ as $\rho(\alpha)$ and
$\rho(\beta)$ commute. (Aside: If $\rho(\alpha)$ is instead parabolic,
we can usually confirm this by showing $\brho(\alpha)$ does not
contain $I$. The remaining ambiguous cases, e.g.~when $\rho(\alpha)$
is trivial and $\rho(\beta)$ is parabolic, are probably best handled
by changing the Dehn filling description.)  Thus, at the end one has a
$\brho \maps \pi_1(Y) \to \PSL{2}{\IR}$ that is certified to contain
some actual representation $\rho \maps \pi_1(Y) \to \PSL{2}{\R}$,
which must be nontrivial since at least one $\rho(\beta)$ is
hyperbolic.  Moreover, we can lift $\rho$ to
$\rhotil \maps \pi_1(Y) \to \SL{2}{\R}$ when
$\pi_1(Y) = \spandef{S}{R_1, \ldots, R_n}$ if we can chose a lift
$\brhotil \maps S \to \SL{2}{\IR}$ such that
$-2 \notin \tr\left(\brhotil(R_i)\right)$ for all $i$.  Here, since I
started with shapes coming from an approximate point in the
$\SL{2}{\C}$ Ptolemy variety, I unsurprisingly always succeed in
finding such a lift.

\subsection{Lifting representations} \label{sec:liftreps}

For ease of notation, let $G = \SL{2}{\R}$, $\IG = \SL{2}{\IR}$, and
$\Gtil = \PSLRtilde$.  Also set $\Gamma = \pi_1(Y)$ with
$\spandef{S}{R_1, \ldots, R_n}$ a finite presentation for
$\Gamma$. Given $\brho \maps S \to \IG$ which is guaranteed to contain
a representation $\rho \maps \Gamma \to G$, we need a way to show
that $\rho$ lifts to $\Gtil$ under
\[
0 \to \Z \to \Gtil \to G \to 1
\]
by a computation using only $\brho$.  A \2-cocycle $c$ representing the
obstruction $e(\rho) \in H^2(\Gamma; \Z)$ to $\rho$ lifting can be
defined by choosing any lift $\rhotil \maps \FreeGroup(S) \to \Gtil$
of $\rho\vert_S$ and setting $c(R_i) = \rhotil(R_i)$ which is in the
center $\Z$ of $\Gtil$.  Thus, the key issue is how to work in $\Gtil$
starting from our approximate representation $\brho$.

\subsection{Working with $\Gtil$} As $\Gtil$ is nonlinear, i.e.~cannot
be embedded in any $\GL{n}{\R}$, it is not easy to work with
computationally, especially as we need to do so rigorously in the
context of interval arithmetic.  I used the following approach,
motivated by \cite[Chapter 2]{Bruggeman1994}, which you can consult for
additional details. In this section, we identify the hyperbolic plane
$\H^2$ with $\setdef{z \in \C}{\Im(z) > 0}$ and focus on the M\"obius
action of $G$ on $\H^2$.  The stabilizer of $i \in \H^2$ is
\[
  K = \setdef{\khat(t) =
   \twobytwomatrix{\cos t}{\sin t}{-\sin t}{\cos t}}%
 {t \in \R\big/2\pi\Z},
\]
where $\khat(t)$ rotates anticlockwise around the point $i$ through
angle $2t$.  The stabilizer of $\infty \in \partial \H^2$ is
\[
  P = \setdef{\phat(z) =
    \twobytwomatrix{\sqrt{y}}{x \big/ \sqrt{y}}{0}{1\big/\sqrt{y}}}{z = x + i y \in \H^2}
\]
where $\phat(z)$ takes $i$ to $z$. Any element
$g = \mysmallmatrix{a}{b}{c}{d}$ in $G$ can be uniquely written as
$\phat(z) \khat(t)$ where $z = g \cdot i$ and
$\khat(t) = \phat(z)^{-1} g$ or alternatively $t = \arg(d - i c)$.
The map $P \times K \to G$ taking $(p, k)$ to $p \cdot k$ is a
diffeomorphism, showing $G \cong \H^2 \times S^1$ as analytic
manifolds.  The universal covering group $\Gtil$ is thus topologically
$\H^2 \times \R$.  Taking $\Ptil \subset \Gtil$ to be the preimage of
$P$, the map $\Ptil \to P$ is a group isomorphism, and we define
$p(z) \in \Ptil$ to be the element mapping to $\phat(z)$.  The
preimage $\Ktil$ of $K$ is isomorphic to $(\R, +)$ with $\Ktil \to K$
having kernel $2 \pi \Z$; we denote the elements of $\Ktil$ by $p(t)$
for $t \in \R$ where $p(t) \cdot p(t') = p(t + t')$ and
$p(t) \mapsto \phat(t)$.  Thus each element of $\Gtil$ can be uniquely
expressed in the form $p(z) k(t)$ for some $z \in \H^2$ and
$t \in \R$.  I will record a $\gtil \in \Gtil$ by the pair
$(g, t) \in G \times \R$ where $\gtil \mapsto g$ and $t$ is the
$\Ktil$-coordinate of $\gtil$; concretely, this identifies $\Gtil$
with
\begin{equation}\label{eq:workGtil}
  \setdef{\Big(\mysmallmatrix{a}{b}{c}{d}, \ t\Big)
    \in G \times \R}{\frac{d - i c}{\sqrt{c^2 + d^2}} = e^{it}}
\end{equation}
If $g$ and $g'$ are in $G$, then setting $z' = g \cdot i$ it follows
from \cite[Page~28]{Bruggeman1994} that the group law on $\Gtil$ can
be written
\begin{equation}\label{eq:mulGtil}
  (g, t) \cdot (g', t') = \left(g \cdot g', \,  t + t' -
    \arg\big(e^{i t}(-z' \sin t + \cos t)\big) \right)
\end{equation}
where here $\arg$ takes values in $(-\pi, \pi]$; one can show the
input to $\arg$ is never in $(-\infty, 0] \subset \C$, so this term is
in fact continuous (indeed analytic) in $z'$ and $t$.  (In contrast,
the book \cite{Bruggeman1994} uses $(z, t)$ as coordinates on $\Gtil$,
which is more concise but the formula for the $z$ part of the group
law is then more complicated than just multiplying matrices.)

\subsection{An interval version of $\Gtil$}
\label{sec:intervalgtil}

Motivated by the description (\ref{eq:workGtil}) for $\Gtil$, define
\[
  \IGtil = \setdef{\Big(\mysmallmatrix{\a}{\b}{\c}{\d}, \t \Big) \in
    \big(M_{2}{\IR} \big) \times \IR}{%
    \parbox{0.55\textwidth}{There exist interval reps
      $a, b, c, d, t \in \R$ such that $a d - b c = 1$ and
      $(d - i c)/\sqrt{c^2 + d^2} = e^{i t}$}}
\]
where here $a \in \a, b \in \b$, etc.  We can view $\Gtil$ as a subset
of $\IGtil$, and given an element
$\g = \mysmallmatrix{\a}{\b}{\c}{\d}$ in $\SL{2}{\IR}$ known to
contain a $g \in \SL{2}{\R}$, then $(\g, \t)$ is in $\IGtil$ if we
take
\begin{equation}\label{eq:tset}
  \t = \arg\big((\d - i \c)/\sqrt{\c^2 + \d^2}\big)
\end{equation}
Using Equation \ref{eq:mulGtil}, we can multiply elements $(\g, \t)$ and
$(\g', \t')$ of $\IGtil$ together provided the element of $\IC$
\[
\bm{r} = e^{i \t}(-\z' \sin \t + \cos \t)\big) \mtext{where $\z' = \g' \cdot i$}
\]
is disjoint from $(-\infty, 0]$, which is something that can be easily
checked from the four corners of the rectangle $\bm{r}$.

\subsection{Proof of Theorem~\ref{thm:tilreps}}  

For 64{,}180 manifolds $Y$ in $\cY$, I found a Dehn surgery
description $Y = M(\alpha)$ where $M$ had an ideal triangulation $\cT$
with $n$ tetrahedra with a point $z_0 \in \R^n$ where the following
held.  First, the method of Section~\ref{sec:certreal} showed there
was a point in $D(\cT, \alpha)$ near $z_0$ such that the associated
$\rho \maps \pi_1(M) \to \PSL{2}{\R}$ factored through to $\pi_1(Y)$,
was nontrivial, and lifted to $G = \SL{2}{\R}$. Second, starting from
the approximate representation $\brho \maps \pi_1(Y) \to \IG$ and a
presentation $\pi_1(Y) = \spandef{S}{R_1, \ldots, R_n}$, I constructed
a lift $\brhotil \maps \FreeGroup(S) \to \IGtil$ of $\brho|_S$ by
setting $\brhotil(\gamma) = (\brho(\gamma), \t)$ where $\t$ is defined
by (\ref{eq:tset}).  For each relator $R_i$, one must have
$\brhotil(R_i) = (\g, \t)$ contains an element of the center of
$\Gtil$, that is $I \in \g$ and $\t \cap 2 \pi \Z \neq \emptyset$.
Provided there is a unique $m \in \Z$ with $2 \pi m \in \t$, it
follows that the 2-cocycle for $e(\rho)$ described in
Section~\ref{sec:liftreps} takes the value $m$ on $R_i$. I then
checked that this cocycle was a coboundary. This proved the existence
of a nontrivial representation $\rhotil \maps \pi_1(Y) \to \Gtil$ and
hence $Y$ is orderable.

As in Section~\ref{sec:praccon}, the amount of precision needed
varied, both to certify the initial solution and to avoid the branch
cut of $\arg$ when computing in $\IG$.  In 95\% of the cases, using
200 bits of precision sufficed, but the hardest 67 examples required
1{,}000 bits.  For each $Y$, the precise Dehn surgery description,
triangulation, the approximate shapes, and level of precision are
available at \cite{PaperData} along with the code used to check them
all.  This completes the proof of Theorem~\ref{thm:tilreps}.

\subsection{Proof of Theorem~\ref{thm:main}(\ref{item:euler0})}

Theorems~\ref{thm:tauteuler0} and \ref{thm:tilreps} respectively show
that at least 32{,}347 and 64{,}180 manifolds in $\cY$ are orderable.
From \cite{PaperData}, we can determine the overlap between the two
methods, and find that at least one of these theorems applies to some
80{,}236 of them. Comparing with the data behind
Theorem~\ref{thm:main}(\ref{item:Lsp}) confirms that these are all
non-L-spaces, proving Theorem~\ref{thm:main}(\ref{item:euler0}).

\section{Code and data}
\label{sec:code}

The full data and code for all these proofs is permanently archived on
the Harvard Dataverse \cite{PaperData}. One method I used to try to
ensure the correctness of the code is worth mentioning. The
results of Theorem~\ref{thm:main} about the manifolds $Y \in \cY$ came
from:
\begin{enumerate}
\item \label{item:cHF} Determining which $Y$ are L-spaces using
  Section~\ref{sec:compute-HF}.
\item \label{item:cNO}
  Showing $Y$ is not orderable using
  Section~\ref{sec:disorder}.
\item \label{item:cFOL} Finding taut foliations from foliar
  orientations (Sections~\ref{sec:foliar} and \ref{sec:persist}) and
  using these to show $Y$ is orderable (Section~\ref{sec:foleuler}).
\item Showing $Y$ is orderable using representations to $\PSLRtilde$
  (Section~\ref{sec:prove_order}).
\end{enumerate}
The only place where information from one part was used in another is
that 4\% of the foliations found in (\ref{item:cFOL}) were used in
(\ref{item:cHF}) to show the last stubborn 2.1\% of manifolds in $\cY$
were not L-spaces. Putting aside this 2.1\%, the above items represent
four completely independent programs.  The key safeguard is that all
four programs were run on the \emph{whole} of $\cY$, even though this
meant huge amounts of time was spent searching for foliar orientations
on manifolds already known to be L-spaces, searching for orders of
manifolds already known not to have any, etc.  Rather that being a
waste of computational resources, this represents compelling evidence
for the correctness of the code since no contradictory results were
obtained!  For example, suppose that the program in (\ref{item:cHF})
had a bug causing it to randomly output the wrong answer for 1 in
every 10,000 inputs.  We would then expect about 16 manifolds in $\cY$
for which a taut foliation would have been found in (\ref{item:cFOL})
but that were reported in (\ref{item:cHF}) as L-spaces, which is of
course impossible \cite{OSgenusbounds2004}.  Indeed, the probability
that no such ``L-space with a taut foliation'' was found with these
assumptions is less than $10^{-7}$.  Hence, one should expect that the
probability that the program in (\ref{item:cHF}) reports non-L-spaces
as L-spaces is well less than $10^{-4}$.

%% file: conjecture.bbl
\newcommand{\etalchar}[1]{$^{#1}$}
\begin{thebibliography}{HIKMOT}

\bibitem[Agol]{Agol2014}
I.~Agol.
\newblock {Virtual properties of 3-manifolds}.
\newblock In {\em Proceedings of the {I}nternational {C}ongress of
  {M}athematicians---{S}eoul 2014. {V}ol. 1}, pages 141--170. Kyung Moon Sa,
  Seoul, 2014.
\newblock \mathreviewsnumber{3728467}.

\bibitem[AL]{AgolLi2003}
I.~Agol and T.~Li.
\newblock \href{http://dx.doi.org/10.2140/gt.2003.7.287}{{An algorithm to
  detect laminar 3-manifolds}}.
\newblock {\em Geom. Topol.} {\bf 7} (2003), 287--309.
\newblock \href{http://arxiv.org/abs/arXiv:math/0201310}{{\tt
  arXiv:math/0201310}}, \mathreviewsnumber{1988287}.

\bibitem[BM]{BakerMoore2014}
K.~L. Baker and A.~H. Moore.
\newblock \href{https://doi.org/10.2969/jmsj/07017484}{{Montesinos knots,
  {H}opf plumbings, and {L}-space surgeries}}.
\newblock {\em J. Math. Soc. Japan} {\bf 70} (2018), 95--110.
\newblock \href{http://arxiv.org/abs/arXiv:1404.7585}{{\tt arXiv:1404.7585}},
  \mathreviewsnumber{3750269}.

\bibitem[BCP]{Magma}
W.~Bosma, J.~Cannon, and C.~Playoust.
\newblock \href{http://dx.doi.org/10.1006/jsco.1996.0125}{{The {M}agma algebra
  system. {I}. {T}he user language}}.
\newblock {\em J. Symbolic Comput.} {\bf 24} (1997), 235--265.
\newblock Computational algebra and number theory (London, 1993).
\newblock \mathreviewsnumber{1484478}.

\bibitem[Bow]{Bowden2016}
J.~Bowden.
\newblock \href{http://dx.doi.org/10.1007/s00039-016-0387-2}{{Approximating
  {$C^0$}-foliations by contact structures}}.
\newblock {\em Geom. Funct. Anal.} {\bf 26} (2016), 1255--1296.
\newblock \href{http://arxiv.org/abs/arXiv:1509.07709}{{\tt arXiv:1509.07709}},
  \mathreviewsnumber{3568032}.

\bibitem[BC]{BoyerClay2014}
S.~Boyer and A.~Clay.
\newblock \href{http://dx.doi.org/10.1016/j.aim.2017.01.026}{{Foliations,
  orders, representations, {L}-spaces and graph manifolds}}.
\newblock {\em Adv. Math.} {\bf 310} (2017), 159--234.
\newblock \href{http://arxiv.org/abs/arXiv:1401.7726}{{\tt arXiv:1401.7726}},
  \mathreviewsnumber{3620687}.

\bibitem[BGW]{BoyerGordonWatson2013}
S.~Boyer, C.~M. Gordon, and L.~Watson.
\newblock \href{http://dx.doi.org/10.1007/s00208-012-0852-7}{{On {L}-spaces and
  left-orderable fundamental groups}}.
\newblock {\em Math. Ann.} {\bf 356} (2013), 1213--1245.
\newblock \href{http://arxiv.org/abs/arXiv:1107.5016}{{\tt arXiv:1107.5016}},
  \mathreviewsnumber{3072799}.

\bibitem[BH]{BoyerHu2018}
S.~Boyer and Y.~Hu.
\newblock {Taut foliations in branched cyclic covers and left-orderable
  groups}.
\newblock Preprint 2018, 40 pages.
\newblock \href{http://arxiv.org/abs/arXiv:1711.04578}{{\tt arXiv:1711.04578}}.

\bibitem[BRW]{BoyerRolfsenWiest2005}
S.~Boyer, D.~Rolfsen, and B.~Wiest.
\newblock \href{http://aif.cedram.org/item?id=AIF_2005__55_1_243_0}{{Orderable
  3-manifold groups}}.
\newblock {\em Ann. Inst. Fourier (Grenoble)} {\bf 55} (2005), 243--288.
\newblock \href{http://arxiv.org/abs/arXiv:math/0211110}{{\tt
  arXiv:math/0211110}}, \mathreviewsnumber{2141698 (2006a:57001)}.

\bibitem[Bri]{Brittenham2001}
M.~Brittenham.
\newblock \href{https://doi.org/10.1142/S0218216501001372}{{Persistent
  laminations from {S}eifert surfaces}}.
\newblock {\em J. Knot Theory Ramifications} {\bf 10} (2001), 1155--1168.
\newblock \href{http://arxiv.org/abs/arXiv:math/9807139}{{\tt
  arXiv:math/9807139}}, \mathreviewsnumber{1871223}.

\bibitem[Bru]{Bruggeman1994}
R.~W. Bruggeman.
\newblock \href{https://doi.org/10.1007/978-3-0346-0336-2}{{\em Families of
  automorphic forms}}, volume~88 of {\em Monographs in Mathematics}.
\newblock Birkh\"{a}user Boston, Inc., Boston, MA, 1994.
\newblock \mathreviewsnumber{1306502}.

\bibitem[Bur]{Burton2014}
B.~A. Burton.
\newblock {The cusped hyperbolic census is complete}.
\newblock Preprint 2014, 32 pages.
\newblock \href{http://arxiv.org/abs/arXiv:1405.2695}{{\tt arXiv:1405.2695}}.

\bibitem[BBP{\etalchar{+}}]{Regina}
B.~A. Burton, R.~Budney, W.~Pettersson, et~al.
\newblock {Regina: Software for low-dimensional topology}.
\newblock \url{http://regina-normal.github.io/}, 1999--2018.

\bibitem[Cal1]{Calegari2000}
D.~Calegari.
\newblock \href{http://dx.doi.org/10.4310/CAG.2000.v8.n1.a5}{{Foliations
  transverse to triangulations of {$3$}-manifolds}}.
\newblock {\em Comm. Anal. Geom.} {\bf 8} (2000), 133--158.
\newblock \href{http://arxiv.org/abs/arXiv:math/9803109}{{\tt
  arXiv:math/9803109}}, \mathreviewsnumber{1730893}.

\bibitem[Cal2]{Calegari2007}
D.~Calegari.
\newblock {\em Foliations and the geometry of 3-manifolds}.
\newblock Oxford Mathematical Monographs. Oxford University Press, Oxford,
  2007.
\newblock \mathreviewsnumber{2327361}.

\bibitem[CaD]{CalegariDunfield2003}
D.~Calegari and N.~M. Dunfield.
\newblock \href{http://dx.doi.org/10.1007/s00222-002-0271-6}{{Laminations and
  groups of homeomorphisms of the circle}}.
\newblock {\em Invent. Math.} {\bf 152} (2003), 149--204.
\newblock \href{http://arxiv.org/abs/arXiv:math/0203192}{{\tt
  arXiv:math/0203192}}, \mathreviewsnumber{1965363 (2005a:57013)}.

\bibitem[CC]{CandelConlon2003}
A.~Candel and L.~Conlon.
\newblock \href{http://dx.doi.org/10.1090/gsm/060}{{\em Foliations. {II}}},
  volume~60 of {\em Graduate Studies in Mathematics}.
\newblock American Mathematical Society, Providence, RI, 2003.
\newblock \mathreviewsnumber{1994394}.

\bibitem[CR]{ClayRolfsen2015}
A.~Clay and D.~Rolfsen.
\newblock {\em Ordered groups and topology}, volume 176 of {\em Graduate
  Studies in Mathematics}.
\newblock American Mathematical Society, Providence, RI, 2016.
\newblock \href{http://arxiv.org/abs/arXiv:1511.05088}{{\tt arXiv:1511.05088}},
  \mathreviewsnumber{3560661}.

\bibitem[CuD]{CullerDunfield2018}
M.~Culler and N.~M. Dunfield.
\newblock \href{https://doi.org/10.2140/gt.2018.22.1405}{{Orderability and
  {D}ehn filling}}.
\newblock {\em Geom. Topol.} {\bf 22} (2018), 1405--1457.
\newblock \href{http://arxiv.org/abs/arXiv:1602.03793}{{\tt arXiv:1602.03793}},
  \mathreviewsnumber{3780437}.

\bibitem[CDGW]{SnapPy}
M.~Culler, N.~M. Dunfield, M.~Goerner, and J.~R. Weeks.
\newblock {Snap{P}y, a computer program for studying the topology and geometry
  of $3$-manifolds}.
\newblock Available at \url{http://snappy.computop.org}.

\bibitem[Del]{Delman1995}
C.~Delman.
\newblock \href{https://doi.org/10.1016/0166-8641(95)00085-U}{{Essential
  laminations and {D}ehn surgery on {$2$}-bridge knots}}.
\newblock {\em Topology Appl.} {\bf 63} (1995), 201--221.
\newblock \mathreviewsnumber{1334307}.

\bibitem[DR]{DelmanRobertsTBD}
C.~Delman and R.~Roberts.
\newblock {Alternating knots and Montesinos knots satisfy the L-space knot
  conjecture}.
\newblock GEAR talk at Univ.~of Illinois by Delman, 2017.
\newblock \url{https://youtu.be/1M8NvoiDGzc}

\bibitem[Dun1]{PaperData}
N.~Dunfield.
\newblock {{Code and data to accompany this paper}}.
\newblock Harvard Dataverse, 2019.
\newblock \url{https://doi.org/10.7910/DVN/LCYXPO}

\bibitem[Dun2]{DunfieldExceptional}
N.~M. Dunfield.
\newblock {A census of exceptional Dehn fillings}.
\newblock Preprint 2018, 15 pages.
\newblock \href{http://arxiv.org/abs/arXiv:1812.11940}{{\tt arXiv:1812.11940}}.

\bibitem[DG]{DunfieldGaroufalidis2012}
N.~M. Dunfield and S.~Garoufalidis.
\newblock
  \href{http://dx.doi.org/10.1090/S0002-9947-2012-05663-7}{{Incompressibility
  criteria for spun-normal surfaces}}.
\newblock {\em Trans. Amer. Math. Soc.} {\bf 364} (2012), 6109--6137.
\newblock \href{http://arxiv.org/abs/arXiv:1102.4588}{{\tt arXiv:1102.4588}},
  \mathreviewsnumber{2946944}.

\bibitem[DHL]{DunfieldHoffmanLicata2015}
N.~M. Dunfield, N.~R. Hoffman, and J.~E. Licata.
\newblock \href{http://dx.doi.org/10.4310/MRL.2015.v22.n6.a7}{{Asymmetric
  hyperbolic {$L$}-spaces, {H}eegaard genus, and {D}ehn filling}}.
\newblock {\em Math. Res. Lett.} {\bf 22} (2015), 1679--1698.
\newblock \href{http://arxiv.org/abs/arXiv:1407.7827}{{\tt arXiv:1407.7827}},
  \mathreviewsnumber{3507256}.

\bibitem[ECHL{\etalchar{+}}]{WordProcessing1992}
D.~B.~A. Epstein, J.~W. Cannon, D.~F. Holt, S.~V.~F. Levy, M.~S. Paterson, and
  W.~P. Thurston.
\newblock {\em Word processing in groups}.
\newblock Jones and Bartlett Publishers, Boston, MA, 1992.
\newblock \mathreviewsnumber{1161694}.

\bibitem[FWW]{FloydWeberWeeks2002}
W.~Floyd, B.~Weber, and J.~Weeks.
\newblock \href{http://projecteuclid.org/euclid.em/1057860318}{{The {A}chilles'
  heel of {$\rm O(3,1)$}?}}
\newblock {\em Experiment. Math.} {\bf 11} (2002), 91--97.
\newblock \mathreviewsnumber{1960304}.

\bibitem[Gab]{Gabai2000}
D.~Gabai.
\newblock \href{http://dx.doi.org/10.1007/s000140050114}{{Combinatorial volume
  preserving flows and taut foliations}}.
\newblock {\em Comment. Math. Helv.} {\bf 75} (2000), 109--124.
\newblock \mathreviewsnumber{1760497}.

\bibitem[GK]{GabaiKazez1990}
D.~Gabai and W.~H. Kazez.
\newblock \href{https://doi.org/10.1016/0166-8641(90)90018-W}{{Pseudo-{A}nosov
  maps and surgery on fibred {$2$}-bridge knots}}.
\newblock {\em Topology Appl.} {\bf 37} (1990), 93--100.
\newblock \mathreviewsnumber{1075377}.

\bibitem[Gao]{Gao2017}
X.~Gao.
\newblock \href{https://doi.org/10.2140/agt.2017.17.2511}{{Non-{L}-space
  integral homology 3-spheres with no nice orderings}}.
\newblock {\em Algebr. Geom. Topol.} {\bf 17} (2017), 2511--2522.
\newblock \href{http://arxiv.org/abs/arXiv:1609.07663}{{\tt arXiv:1609.07663}},
  \mathreviewsnumber{3686404}.

\bibitem[Gar]{Gardam2018}
G.~Gardam.
\newblock {Profinite rigidity in the SnapPea census}.
\newblock Preprint 2018, 16 pages.
\newblock \href{http://arxiv.org/abs/arXiv:1805.02697}{{\tt arXiv:1805.02697}}.

\bibitem[GGZ1]{GaroufalidisGoernerZickert2015b}
S.~Garoufalidis, M.~Goerner, and C.~K. Zickert.
\newblock \href{https://doi.org/10.2140/agt.2015.15.565}{{Gluing equations for
  {$\mathrm{PGL}(n,\Bbb{C})$}-representations of 3-manifolds}}.
\newblock {\em Algebr. Geom. Topol.} {\bf 15} (2015), 565--622.
\newblock \href{http://arxiv.org/abs/arXiv:1207.6711}{{\tt arXiv:1207.6711}},
  \mathreviewsnumber{3325748}.

\bibitem[GGZ2]{GaroufalidisGoernerZickert2015a}
S.~Garoufalidis, M.~Goerner, and C.~K. Zickert.
\newblock \href{https://doi.org/10.2140/agt.2015.15.371}{{The {P}tolemy field
  of 3-manifold representations}}.
\newblock {\em Algebr. Geom. Topol.} {\bf 15} (2015), 371--397.
\newblock \href{http://arxiv.org/abs/arXiv:1310.2497}{{\tt arXiv:1310.2497}},
  \mathreviewsnumber{3325740}.

\bibitem[GTZ]{GaroufalidisThurstonZickert2015}
S.~Garoufalidis, D.~P. Thurston, and C.~K. Zickert.
\newblock \href{https://doi.org/10.1215/00127094-3121185}{{The complex volume
  of {${\rm SL}(n,\Bbb{C})$}-representations of 3-manifolds}}.
\newblock {\em Duke Math. J.} {\bf 164} (2015), 2099--2160.
\newblock \href{http://arxiv.org/abs/arXiv:1111.2828}{{\tt arXiv:1111.2828}},
  \mathreviewsnumber{3385130}.

\bibitem[Ghys]{Ghys2001}
{\'E}.~Ghys.
\newblock
  \href{http://www.umpa.ens-lyon.fr/~ghys/articles/groups-acting-circle.pdf}{{Groups
  acting on the circle}}.
\newblock {\em Enseign. Math. (2)} {\bf 47} (2001), 329--407.
\newblock \mathreviewsnumber{1876932 (2003a:37032)}.

\bibitem[Goe]{GoernerPtolemyData}
M.~Goerner.
\newblock {Ptolemy data for SnapPy manifolds}.
\newblock \url{http://ptolemy.unhyperbolic.org/}

\bibitem[GZ]{GoernerZickert2018}
M.~Goerner and C.~K. Zickert.
\newblock \href{https://doi.org/10.1007/s00209-017-1970-4}{{Triangulation
  independent {P}tolemy varieties}}.
\newblock {\em Math. Z.} {\bf 289} (2018), 663--693.
\newblock \href{http://arxiv.org/abs/arXiv:1507.03238}{{\tt arXiv:1507.03238}},
  \mathreviewsnumber{3803807}.

\bibitem[GL]{GordonLidman2014}
C.~Gordon and T.~Lidman.
\newblock \href{http://dx.doi.org/10.1007/s40306-014-0091-y}{{Taut foliations,
  left-orderability, and cyclic branched covers}}.
\newblock {\em Acta Math. Vietnam.} {\bf 39} (2014), 599--635.
\newblock \href{http://arxiv.org/abs/arXiv:1406.6718}{{\tt arXiv:1406.6718}},
  \mathreviewsnumber{3292587}.

\bibitem[Han1]{Hanselman2016}
J.~Hanselman.
\newblock \href{https://doi.org/10.2140/agt.2016.16.3103}{{Bordered {H}eegaard
  {F}loer homology and graph manifolds}}.
\newblock {\em Algebr. Geom. Topol.} {\bf 16} (2016), 3103--3166.
\newblock \href{http://arxiv.org/abs/arXiv:1310.6696}{{\tt arXiv:1310.6696}},
  \mathreviewsnumber{3584255}.

\bibitem[Han2]{HFhatGraph}
J.~Hanselman.
\newblock {HFhat\_graph\_manifolds}, 2016.
\newblock \url{https://github.com/hanselman/HFhat_graph_manifolds}

\bibitem[HRRW]{HanselmanRasmussenRasmussenWatson2015}
J.~Hanselman, J.~Rasmussen, S.~D. Rasmussen, and L.~Watson.
\newblock {L-spaces, taut foliations, and graph manifolds}.
\newblock Preprint 2015, 9 pages.
\newblock \href{http://arxiv.org/abs/arXiv:1508.05911}{{\tt arXiv:1508.05911}}.

\bibitem[Hat]{Hatcher2002}
A.~Hatcher.
\newblock {\em Algebraic topology}.
\newblock Cambridge University Press, Cambridge, 2002.
\newblock \mathreviewsnumber{1867354}.

\bibitem[Hed]{Hedden2009}
M.~Hedden.
\newblock {On knot {F}loer homology and cabling. {II}}.
\newblock {\em Int. Math. Res. Not. IMRN} (2009), 2248--2274.
\newblock \href{http://arxiv.org/abs/arXiv:0806.2172}{{\tt arXiv:0806.2172}},
  \mathreviewsnumber{2511910}.

\bibitem[HK]{HirasawaKobayashi2001}
M.~Hirasawa and T.~Kobayashi.
\newblock \href{http://projecteuclid.org/euclid.ojm/1153492557}{{Pre-taut
  sutured manifolds and essential laminations}}.
\newblock {\em Osaka J. Math.} {\bf 38} (2001), 905--922.
\newblock \mathreviewsnumber{1864469}.

\bibitem[HW]{HodgsonWeeks1994}
C.~D. Hodgson and J.~R. Weeks.
\newblock \href{http://projecteuclid.org/euclid.em/1048515809}{{Symmetries,
  isometries and length spectra of closed hyperbolic three-manifolds}}.
\newblock {\em Experiment. Math.} {\bf 3} (1994), 261--274.
\newblock \mathreviewsnumber{1341719}.

\bibitem[HIKMOT]{HIKMOT2016}
N.~Hoffman, K.~Ichihara, M.~Kashiwagi, H.~Masai, S.~Oishi, and A.~Takayasu.
\newblock \href{http://dx.doi.org/10.1080/10586458.2015.1029599}{{Verified
  computations for hyperbolic 3-manifolds}}.
\newblock {\em Exp. Math.} {\bf 25} (2016), 66--78.
\newblock \href{http://arxiv.org/abs/arXiv:1310.3410}{{\tt arXiv:1310.3410}},
  \mathreviewsnumber{3424833}.

\bibitem[HW]{HoffmanWalsh2015}
N.~R. Hoffman and G.~S. Walsh.
\newblock \href{http://dx.doi.org/10.1090/bproc/20}{{The big {D}ehn surgery
  graph and the link of {$S^3$}}}.
\newblock {\em Proc. Amer. Math. Soc. Ser. B} {\bf 2} (2015), 17--34.
\newblock \href{http://arxiv.org/abs/arXiv:1311.3980}{{\tt arXiv:1311.3980}},
  \mathreviewsnumber{3422666}.

\bibitem[HEO]{HoltEickOBrien2005}
D.~F. Holt, B.~Eick, and E.~A. O'Brien.
\newblock \href{http://dx.doi.org/10.1201/9781420035216}{{\em Handbook of
  computational group theory}}.
\newblock Discrete Mathematics and its Applications (Boca Raton). Chapman \&
  Hall/CRC, Boca Raton, FL, 2005.
\newblock \mathreviewsnumber{2129747}.

\bibitem[HTW]{HosteThistlethwaiteWeeks1998}
J.~Hoste, M.~Thistlethwaite, and J.~Weeks.
\newblock \href{http://dx.doi.org/10.1007/BF03025227}{{The first 1,701,936
  knots}}.
\newblock {\em Math. Intelligencer} {\bf 20} (1998), 33--48.
\newblock \mathreviewsnumber{1646740}.

\bibitem[Joh]{Johansson2017}
F.~Johansson.
\newblock \href{http://dx.doi.org/0.1109/TC.2017.2690633}{{Arb: efficient
  arbitrary-precision midpoint-radius interval arithmetic}}.
\newblock {\em IEEE Transactions on Computers} {\bf 66} (2017), 1281--1292.
\newblock \url{http://arblib.org/}

\bibitem[Juh]{Juhasz2015}
A.~Juh\'asz.
\newblock \href{http://dx.doi.org/10.1142/9789814630627_0007}{{A survey of
  {H}eegaard {F}loer homology}}.
\newblock In {\em New ideas in low dimensional topology}, volume~56 of {\em
  Ser. Knots Everything}, pages 237--296. World Sci. Publ., Hackensack, NJ,
  2015.
\newblock \href{http://arxiv.org/abs/arXiv:1310.3418}{{\tt arXiv:1310.3418}},
  \mathreviewsnumber{3381327}.

\bibitem[KR]{KazezRoberts2017}
W.~H. Kazez and R.~Roberts.
\newblock \href{https://doi.org/10.2140/gt.2017.21.3601}{{{$C^0$}
  approximations of foliations}}.
\newblock {\em Geom. Topol.} {\bf 21} (2017), 3601--3657.
\newblock \href{http://arxiv.org/abs/arXiv:1509.08382}{{\tt arXiv:1509.08382}},
  \mathreviewsnumber{3693573}.

\bibitem[KMOS]{KronheimerMrowkaOzsvathSzabo2007}
P.~Kronheimer, T.~Mrowka, P.~Ozsv\'ath, and Z.~Szab\'o.
\newblock \href{http://dx.doi.org/10.4007/annals.2007.165.457}{{Monopoles and
  lens space surgeries}}.
\newblock {\em Ann. of Math. (2)} {\bf 165} (2007), 457--546.
\newblock \href{http://arxiv.org/abs/arXiv:math/0310164}{{\tt
  arXiv:math/0310164}}, \mathreviewsnumber{2299739}.

\bibitem[Li]{Li2002}
T.~Li.
\newblock \href{http://dx.doi.org/10.2140/gt.2002.6.153}{{Laminar branched
  surfaces in 3-manifolds}}.
\newblock {\em Geom. Topol.} {\bf 6} (2002), 153--194.
\newblock \href{http://arxiv.org/abs/arXiv:math/0204012}{{\tt
  arXiv:math/0204012}}, \mathreviewsnumber{1914567}.

\bibitem[LOT1]{LipshitzOzsvathThurston2011}
R.~Lipshitz, P.~S. Ozsv\'ath, and D.~P. Thurston.
\newblock \href{http://dx.doi.org/10.1073/pnas.1019060108}{{Tour of bordered
  {F}loer theory}}.
\newblock {\em Proc. Natl. Acad. Sci. USA} {\bf 108} (2011), 8085--8092.
\newblock \href{http://arxiv.org/abs/arXiv:1107.5621}{{\tt arXiv:1107.5621}},
  \mathreviewsnumber{2806643}.

\bibitem[LOT2]{LipshitzOzsvathThurston2014}
R.~Lipshitz, P.~S. Ozsv\'ath, and D.~P. Thurston.
\newblock {Bordered Heegaard Floer homology: Invariance and pairing}.
\newblock Preprint 2014, 283 pages.
\newblock \href{http://arxiv.org/abs/arXiv:0810.0687}{{\tt arXiv:0810.0687}}.

\bibitem[MT]{MatsuzakiTaniguchi1998}
K.~Matsuzaki and M.~Taniguchi.
\newblock {\em Hyperbolic manifolds and {K}leinian groups}.
\newblock Oxford Mathematical Monographs. The Clarendon Press, Oxford
  University Press, New York, 1998.
\newblock Oxford Science Publications.
\newblock \mathreviewsnumber{1638795}.

\bibitem[MKC]{MooreKearfottCloud2009}
R.~E. Moore, R.~B. Kearfott, and M.~J. Cloud.
\newblock \href{http://dx.doi.org/10.1137/1.9780898717716}{{\em Introduction to
  interval analysis}}.
\newblock Society for Industrial and Applied Mathematics (SIAM), Philadelphia,
  PA, 2009.
\newblock \mathreviewsnumber{2482682}.

\bibitem[Nai]{Naimi1997}
R.~Naimi.
\newblock \href{https://doi.org/10.2140/pjm.1997.180.153}{{Constructing
  essential laminations in {$2$}-bridge knot surgered {$3$}-manifolds}}.
\newblock {\em Pacific J. Math.} {\bf 180} (1997), 153--186.
\newblock \mathreviewsnumber{1474899}.

\bibitem[Neu]{Neumann1981}
W.~D. Neumann.
\newblock \href{https://doi.org/10.2307/1999331}{{A calculus for plumbing
  applied to the topology of complex surface singularities and degenerating
  complex curves}}.
\newblock {\em Trans. Amer. Math. Soc.} {\bf 268} (1981), 299--344.
\newblock \mathreviewsnumber{632532}.

\bibitem[OS1]{OSgenusbounds2004}
P.~Ozsv{\'a}th and Z.~Szab{\'o}.
\newblock \href{http://dx.doi.org/10.2140/gt.2004.8.311}{{Holomorphic disks and
  genus bounds}}.
\newblock {\em Geom. Topol.} {\bf 8} (2004), 311--334.
\newblock \href{http://arxiv.org/abs/arXiv:math/0311496}{{\tt
  arXiv:math/0311496}}, \mathreviewsnumber{2023281 (2004m:57024)}.

\bibitem[OS2]{OSLensSpace2005}
P.~Ozsv{\'a}th and Z.~Szab{\'o}.
\newblock \href{http://dx.doi.org/10.1016/j.top.2005.05.001}{{On knot {F}loer
  homology and lens space surgeries}}.
\newblock {\em Topology} {\bf 44} (2005), 1281--1300.
\newblock \href{http://arxiv.org/abs/arXiv:math/0303017}{{\tt
  arXiv:math/0303017}}, \mathreviewsnumber{2168576 (2006f:57034)}.

\bibitem[RR]{Rasmussens2015}
J.~Rasmussen and S.~D. Rasmussen.
\newblock \href{https://doi.org/10.1016/j.aim.2017.10.014}{{Floer simple
  manifolds and {L}-space intervals}}.
\newblock {\em Adv. Math.} {\bf 322} (2017), 738--805.
\newblock \href{http://arxiv.org/abs/arXiv:1508.05900}{{\tt arXiv:1508.05900}},
  \mathreviewsnumber{3720808}.

\bibitem[Ras]{Rasmussen2017}
S.~D. Rasmussen.
\newblock \href{https://doi.org/10.1112/S0010437X16008319}{{L-space intervals
  for graph manifolds and cables}}.
\newblock {\em Compos. Math.} {\bf 153} (2017), 1008--1049.
\newblock \href{http://arxiv.org/abs/arXiv:1511.04413}{{\tt arXiv:1511.04413}},
  \mathreviewsnumber{3705248}.

\bibitem[Ril]{Riley2013}
R.~Riley.
\newblock \href{http://dx.doi.org/10.1016/j.exmath.2013.01.003}{{A personal
  account of the discovery of hyperbolic structures on some knot complements}}.
\newblock {\em Expo. Math.} {\bf 31} (2013), 104--115.
\newblock \href{http://arxiv.org/abs/arXiv:1301.4601}{{\tt arXiv:1301.4601}},
  \mathreviewsnumber{3057120}.

\bibitem[Rum]{Rump2010}
S.~M. Rump.
\newblock \href{http://dx.doi.org/10.1017/S096249291000005X}{{Verification
  methods: rigorous results using floating-point arithmetic}}.
\newblock {\em Acta Numer.} {\bf 19} (2010), 287--449.
\newblock \mathreviewsnumber{2652784}.

\bibitem[SW]{SarkarWang2010}
S.~Sarkar and J.~Wang.
\newblock \href{http://dx.doi.org/10.4007/annals.2010.171.1213}{{An algorithm
  for computing some {H}eegaard {F}loer homologies}}.
\newblock {\em Ann. of Math. (2)} {\bf 171} (2010), 1213--1236.
\newblock \href{http://arxiv.org/abs/arXiv:math/0607777}{{\tt
  arXiv:math/0607777}}, \mathreviewsnumber{2630063}.

\bibitem[Seg]{Segerman2012}
H.~Segerman.
\newblock \href{https://doi.org/10.2140/agt.2012.12.2179}{{A generalisation of
  the deformation variety}}.
\newblock {\em Algebr. Geom. Topol.} {\bf 12} (2012), 2179--2244.
\newblock \href{http://arxiv.org/abs/arXiv:0904.1893}{{\tt arXiv:0904.1893}},
  \mathreviewsnumber{3020204}.

\bibitem[SVW]{SommeseVerscheldeWampler2005}
A.~J. Sommese, J.~Verschelde, and C.~W. Wampler.
\newblock \href{https://doi.org/10.1007/3-540-27357-3_8}{{Introduction to
  numerical algebraic geometry}}.
\newblock In {\em Solving polynomial equations}, volume~14 of {\em Algorithms
  Comput. Math.}, pages 301--335. Springer, Berlin, 2005.
\newblock \mathreviewsnumber{2161992}.

\bibitem[Soos]{CryptoMiniSat}
M.~Soos.
\newblock {CryptoMiniSat, Version 5.0}, 2016.
\newblock \url{https://github.com/msoos/cryptominisat}

\bibitem[Sage]{SageMath}
W.~Stein et~al.
\newblock {\em {S}ageMath, the {S}age {M}athematics {S}oftware {S}ystem
  ({V}ersion 8.3)}, 2018.
\newblock \url{http://www.sagemath.org}

\bibitem[Thu1]{ThurstonsNotes}
W.~P. Thurston.
\newblock {Geometry and topology of three-manifolds}.
\newblock Unpublished lecture notes, 1980.
\newblock \url{http://library.msri.org/books/gt3m/}

\bibitem[Thu2]{Thurston1997}
W.~P. Thurston.
\newblock {\em Three-dimensional geometry and topology. {V}ol. 1}, volume~35 of
  {\em Princeton Mathematical Series}.
\newblock Princeton University Press, Princeton, NJ, 1997.
\newblock Edited by Silvio Levy.
\newblock \mathreviewsnumber{1435975}.

\bibitem[Trn]{Trnkova2017}
M.~Trnkov\'a.
\newblock {Rigorous computations with an approximate Dirichlet domain}.
\newblock Preprint 2017, 10 pages.
\newblock \href{http://arxiv.org/abs/arXiv:1703.02595}{{\tt arXiv:1703.02595}}.

\bibitem[Ver]{Verschelde1999}
J.~Verschelde.
\newblock \href{http://doi.acm.org/10.1145/317275.317286}{{Algorithm 795:
  PHCpack: A General-purpose Solver for Polynomial Systems by Homotopy
  Continuation}}.
\newblock {\em ACM Trans. Math. Softw.} {\bf 25} (June 1999), 251--276.

\bibitem[V{\etalchar{+}}]{PHCPack}
J.~Verschelde et~al.
\newblock {PHCpack: solving polynomial systems via homotopy continuation}.
\newblock \url{http://math.uic.edu/~jan}, 1999--2015.

\bibitem[Wee]{Weeks2005}
J.~Weeks.
\newblock
  \href{http://dx.doi.org/10.1016/B978-044451452-3/50011-3}{{Computation of
  hyperbolic structures in knot theory}}.
\newblock In {\em Handbook of knot theory}, pages 461--480. Elsevier B. V.,
  Amsterdam, 2005.
\newblock \href{http://arxiv.org/abs/arXiv:math/0309407}{{\tt
  arXiv:math/0309407}}, \mathreviewsnumber{2179268}.

\bibitem[Zhan]{BFHPython}
B.~Zhan.
\newblock {\texttt{bfh\_python}: Computations in Bordered Heegaard Floer
  Homology}, 2014.
\newblock \url{https://github.com/bzhan/bfh_python}

\bibitem[Zic]{Zickert2016}
C.~K. Zickert.
\newblock \href{https://doi.org/10.1007/s00209-015-1608-3}{{Ptolemy
  coordinates, {D}ehn invariant and the {$A$}-polynomial}}.
\newblock {\em Math. Z.} {\bf 283} (2016), 515--537.
\newblock \href{http://arxiv.org/abs/arXiv:1405.0025}{{\tt arXiv:1405.0025}},
  \mathreviewsnumber{3489078}.

\end{thebibliography}
